\documentclass[a4paper,leqno,10pt]{article}
\usepackage{multicol,color}
\usepackage{amsmath}
\usepackage{latexsym,amsbsy,amssymb,amsfonts,bm,mathtools}
\usepackage{enumerate}
\usepackage{pifont,bbding,umoline}
\usepackage[latin1]{inputenc}
\def\disp{\displaystyle}

\def\Limsup{\mathop{{\rm Lim}\,{\rm sup}}}

\def\tto{\;{\lower 1pt\hbox{$\rightarrow$}}\kern-10pt
\hbox{\raise 2pt\hbox{$\rightarrow$}}\;}
\def\Hat{\widehat}
\def\hat{\widehat}
\def\Tilde{\widetilde}
\def\tilde{\widetilde}

\def\ra{\rangle}
\def\la{\langle}
\def\ve{\varepsilon}
\def\B{\mathbb{B}}
\def\h{\hfill\Box}
\def\R{\mathbb{R}}
\def\N{\mathbb{N}}

\def\t{\tau}
\def\ox{\bar{x}}
\def\op{\bar{p}}
\def\oy{\bar{y}}
\def\oz{\bar{z}}
\def\ov{\bar{v}}
\def\ou{\bar{u}}
\def\ow{\bar{w}}

\def\ob{\bar{b}}
\def\co{\mbox{\rm co}\,}
\def\gph{\mbox{\rm gph}\,}
\def\epi{\mbox{\rm epi}\,}

\def\dom{\mbox{\rm dom}\,}

\def\h{\hfill\triangle}
\def\dn{\downarrow}
\def\O{\Omega}

\def\ph{\varphi}
\def\emp{\emptyset}

\def\oR{\overline{\R}}

\def\lm{\lambda}
\def\gg{\gamma}
\def\dd{\delta}
\def\al{\alpha}
\def\bb{\beta}

\def\th{\theta}
\def\vt{\vartheta}
\newcounter{lk}

\begin{document}
\begin{center}
\vspace*{0.3in}{\bf OPTIMAL CONTROL OF THE SWEEPING PROCESS\\OVER POLYHEDRAL CONTROLLED SETS}\\[3ex]
G. COLOMBO\footnote{Department of Mathematics, University of Padova, Padua, Italy (colombo@math.unipd.it). Research of this author was partially
supported by GNAMPA of INdAM.},
R. HENRION\footnote{Weierstrass Institute for Applied Analysis and Stochastics, Berlin, Germany (henrion@wias-berlin.de).
Research of this author was partially supported by the DFG Research Center MATHEON.},
N. D. HOANG\footnote{Departamento de Matem\'atica, Universidad T\'echnica Federico Santa Mar\'ia, Valpara\'iso, Chile (hoang.nguyen@usm.cl).
Research of this author was partially supported by FONDECYT $\#$3140060 and Basal Projects at CMM, Universidad de Chile.} and
B. S. MORDUKHOVICH\footnote{Department of Mathematics, Wayne State University, Detroit, Michigan, USA (boris@math.wayne.edu).
Research of this author was partially supported by the USA National Science Foundation under grants DMS-1007132 and DMS-1512846.}
\end{center}
\small{\bf Abstract.} The paper addresses a new class of optimal control problems governed by the dissipative and discontinuous
differential inclusion of the sweeping/Moreau process while using controls to determine the best shape of moving convex polyhedra
in order to optimize the given Bolza-type functional, which depends on control and state variables as well as their velocities.
Besides the highly non-Lipschitzian nature of the unbounded differential inclusion of the controlled sweeping process,
the optimal control problems under consideration contain intrinsic state constraints of the inequality and equality types.
All of this creates serious challenges for deriving necessary optimality conditions. We develop here the method of discrete
approximations and combine it with advanced tools of first-order and second-order variational analysis and generalized differentiation.
This approach allows us to establish constructive necessary optimality conditions for local minimizers of the controlled sweeping
process expressed entirely in terms of the problem data under fairly unrestrictive assumptions. As a by-product of the developed approach,
we prove the strong $W^{1,2}$-convergence of optimal solutions of discrete approximations to a given local minimizer of the
continuous-time system and derive necessary optimality conditions for the discrete counterparts.
The established necessary optimality conditions for the sweeping process are illustrated by several examples.\\[1ex]
{\bf Abbreviated title.} Optimal control of the sweeping process

\newtheorem{Theorem}{Theorem}[section]
\newtheorem{Proposition}[Theorem]{Proposition}
\newtheorem{Remark}[Theorem]{Remark}
\newtheorem{Lemma}[Theorem]{Lemma}
\newtheorem{Corollary}[Theorem]{Corollary}
\newtheorem{Definition}[Theorem]{Definition}
\newtheorem{Example}[Theorem]{Example}
\renewcommand{\theequation}{\thesection.\arabic{equation}}
\normalsize
\def\proof{\normalfont\medskip
{\noindent\itshape Proof.\hspace*{6pt}\ignorespaces}}
\def\endproof{$\h$ \vspace*{0.1in}}

\section{Introduction and Problem Formulation}
\setcounter{equation}{0}

This paper is devoted to the study of the {\em sweeping process}, a class of models introduced by Jean-Jacques Moreau in the 1970s
to describe a number of quasistatic mechanical problems; see \cite{mor_frict,mor_jde,mor_unilat} and the book \cite{MM} for more details.
Besides the original motivations, models of this type have found significant applications to elastoplasticity \cite{HR}, hysteresis \cite{kre},
electric circuits \cite{ABB}, etc. For its own sake, the sweeping process theory has become an important area of nonlinear and
variational analysis with numerous mathematical achievements and challenging open questions; see, e.g., \cite{CT,KMM} and the references therein.

Mathematically the sweeping process is governed by the dissipative differential inclusion
\begin{equation}\label{sw}
\dot{x}(t)\in-N\big(x(t);C(t)\big)\;\mbox{ a.e. }\;t\in[0,T],
\end{equation}
which describes the movement of a point belonging to a continuous moving set $C(t)$ while its velocity belongs for a.e.\ $t$
to the negative normal cone to $C(t)$ at $x(t)$. The Cauchy problem $x(0)=x_0$ for the sweeping process \eqref{sw} enjoys a developed
well-posedness theory for convex and mildly nonconvex moving sets; see, e.g., \cite{CT}. Higher-order and state-dependent
sweeping processes have also been studied (to a much lesser extent) in the literature; see, e.g., the book \cite{b}.
Let us also mention the recent paper \cite{mrs}, which contains existence and well-posedness results obtained via advanced
tools of variational analysis and generalized differentiation for a broad class of evolution systems including the sweeping process.

Among the central issues of the sweeping process theory is establishing the existence and {\em uniqueness} of solutions
to the Cauchy problem for the sweeping differential inclusion \eqref{sw} under reasonable assumptions on the given moving set $C(t)$.
This tells us that it does not make any sense to optimize the sweeping process generated by the {\em given set} $C(t)$
from the standard viewpoint of optimal control theory well developed for Lipschitzian differential inclusions and the like;
see, e.g., \cite{m-book1,sm,v} and the references therein.

In our first paper on the sweeping process \cite{chhm} we suggested to take a new viewpoint on optimizing the dynamical system \eqref{sw}
by {\em controlling} the moving set $C(t)$ with the usage of control actions that change the {\em shape} of $C(t)$ and hence the right-hand
side of the sweeping differential inclusion \eqref{sw}. This idea was partly implemented in \cite{chhm} for the case when the sweeping process
was driven by a moving affine hyperplane whose normal direction and boundary were acting as control variables.
Furthermore, it was assumed in \cite{chhm} the independence of the running cost on time, control variables, and
control velocities as well as the uniform Lipschitzian continuity of feasible controls.
Apparently this first attempt was limited from both viewpoints of control theory and possible applications.
A much more realistic while significantly more challenging case appears when, along with general
running costs depending on {\em state, control}, and their {\em velocity} variables, controlled moving
sets are described by {\em convex polyhedra} governed by finitely many controls in normal directions
and polyhedron boundaries under inequality and equality constraints. Polyhedral descriptions of moving
sets in the (uncontrolled) sweeping process were largely explored, e.g., in \cite{HR,kre,KV},
where the reader can find interesting applications to particular models of elastoplasticity and hysteresis.

In the other line of development we mention the recent paper \cite{bk} and its subsequent extension \cite{ao}, which address a different class of optimal control
problems for an equivalent variational inequality description of the sweeping process of the rate-independent hysteresis type, where the convex moving set is fixed
while controls appear in an associated ordinary differential equation. Another recent paper \cite{MCF}, in the framework of $BV$
solutions of a sweeping process whose given moving set is lower semicontinuous with nonempty interior, concerns relaxation issues and dynamic programming.
Controls appear there via perturbations of the dynamics given by the normal cone while being the barycenter of a Borel finite measure.\vspace*{0.02in}

In this paper we study the following {\em optimal control problem} $(P)$ of the generalized Bolza type for the sweeping process \eqref{sw}
as well as some of its modifications. Given an extended-real-valued terminal cost function $\ph\colon\R^n\to\oR:=(-\infty,\infty]$
and a running cost $\ell\colon[0,T]\times\R^{2(n+nm+m)}\to\oR$, minimize the functional
\begin{equation}\label{eq:MP}
J[x,u,b]\colon=\varphi\big(x(T)\big)+\int_0^T\ell\big(t,x(t),u(t),b(t),\dot{x}(t),\dot{u}(t),\dot{b}(t)\big)\,dt
\end{equation}
over the controlled sweeping dynamics described by
\begin{eqnarray}\label{sw-con1}
\mathop x\limits^.(t)\in-N\big(x(t);C(t)\big)\;\mbox{ for a.e. }\;t\in[0,T],\quad x(0):=x_0\in C(0)
\end{eqnarray}
with the inequality and equality constraint defined by
\begin{eqnarray}\label{sw-con2}
C(t):=\big\{x\in\R^n\big|\;\la u_i(t),x\ra\le b_i(t),\;i=1,\ldots,m\big\}
\end{eqnarray}
\begin{eqnarray}\label{sw-con3}
\mbox{with }\;\|u_i(t)\|=1\;\mbox{ for all }\;t\in[0,T],\;i=1,\ldots,m,
\end{eqnarray}
where the controls actions $u(\cdot)=\big(u_1(\cdot),\ldots,u_m(\cdot)\big)$ and $b(\cdot)=\big(b_1(\cdot),\ldots,b_m(\cdot)\big)$ are
{\em absolutely continuous} on $[0,T]$, the final time $T$ is fixed, and the absolutely continuous trajectories $x(\cdot)$
of the differential inclusion are understood in the standard sense of Carath\'eodory. This class of problems contains several
novel features, which either have never been investigated or have been studied insufficiently in control theory;
see more discussions below. Recall now that the {\em normal cone} to a {\em convex} set $C\subset\R^n$ is defined by
\begin{eqnarray}\label{nor}
N(x;C):=\big\{v\in\R^n\big|\;\la v,y-x\ra\le 0,\;\;y\in C\big\}\;\mbox{ if }\;x\in C\;\mbox{ and }\;N(x;C):=\emp\;\mbox{ if }\;x\notin C.
\end{eqnarray}
Hence the sweeping process with a given moving set \eqref{sw} can be considered as an evolution variational inequality,
or a differential variational inequality in the terminology of \cite{ps}. On the other hand, the control model
\eqref{sw-con1}--\eqref{sw-con3} relates rather to evolutionary {\em quasi-variational inequalities} with
controlled parameters, which should be determined to optimize the dynamical process.

The {\em main goal} of this paper is to derive {\em necessary optimality conditions} entirely in terms of the problem
data for the so-called {\em intermediate local minimizers} of $(P)$ and its modifications that occupy an intermediate
position between the standard notions of weak and strong local minima in variational and control problems; see Section~3
for more discussions. Our approach is based on developing an appropriate version of the {\em method of discrete approximations},
which largely follows the scheme of \cite{m95,m-book1} implemented therein for the case of Lipschitzian and uniformly bounded
differential inclusions, while now requiring a novel extension to the case of totally non-Lipschitzian and unbounded differential
inclusions in \eqref{sw-con1}. Some results on discrete approximations of feasible trajectories of \eqref{eq:MP} and the convergence
of optimal solutions to appropriate discretizations of the continuous-time system in \eqref{eq:MP}--\eqref{sw-con3} have been recently obtained in
our preceding paper \cite{chhm2}. However, they do not provide enough information for passing to the limit in necessary optimality
conditions for discrete approximations and thus establishing in this way necessary optimality conditions for local minimizers of
the original continuous-time systems governed by the controlled sweeping process.

In this paper we are going to proceed further in this direction by improving the previous discrete approximation results to make
it possible deriving necessary optimality conditions for the continuous-time systems by passing to the limit from those for their
discrete approximations that are proved to satisfy the desired well-posedness and convergence properties.
The realization of this approach requires overcoming significant difficulties, which have {\em never been addressed}
earlier in control theory from this or any other method of deriving necessary optimality conditions even for more
simple problems with smooth data. Besides the aforementioned totally {\em non-Lipschitzian} and {\em unbounded} nature
of the sweeping process, serious challenges come, in particular, from the intrinsic presence of {\em state constraints}
of the {\em inequality} and {\em equality} types combined with the {\em quasi-variational inequality} structure of the
controlled sweeping process. Indeed, we show in Section~3 that
problem $(P)$ can be rewritten in the more conventional form of the generalized Bolza problem for a non-Lipschitzian
and unbounded differential inclusion with a {\em fixed} right-hand side, where the relations in \eqref{sw-con2} and
\eqref{sw-con3} are treated as state constraints of the inequality and equality type, respectively.
It is worth mentioning that, in contrast to the inequality state constraints well studied for standard control
systems and Lipschitzian differential inclusions (see, e.g., \cite{aa,v} and the references therein),
the {\em equality} state constraints have been just very recently addressed in \cite{ak} for smooth
control systems under regularity assumptions formulated via full ranks of the corresponding Jacobians
of the constraint functions. Needless to say that neither the results nor the approach of \cite{ak} can be applied in our setting.

A crucial ingredient of our approach within the method of discrete approximations is applying advanced nonconvex
tools of first-order and second-order of {\em variational analysis} and {\em generalized differentiation},
which are required even in the case of smooth terminal and running costs in the convex cone setting of \eqref{sw}.
This allows us not only to establish the desired strong convergence of discrete approximations and then to pass to
the liming in the necessary optimality conditions obtained for discrete problems, but also to derive necessary optimality
conditions for the continuous-time control system entirely {\em in terms of the problem data} and the
{\em given} local optimal solution to the controlled sweeping process.

In fact, our major necessary optimality conditions are derived for a certain {\em parametric perturbation}
$(P^\tau)$ of the original problem $(P)$ with the control constraints in \eqref{sw-con3} replaced by
\begin{eqnarray}\label{delta}
\|u_i(t)\|=1\;\mbox{ on }\;[\tau,T-\tau]\;\mbox{ and }\;\frac{1}{2}\le\|u_i(t)\|\le\frac{3}{2}\;\mbox{ on }\;[0,\tau)\cup(T-\tau,T],\;i=1,\ldots,m,
\end{eqnarray}
where the time endpoint perturbation parameter $\tau>0$ is arbitrarily small, and so $(P^\tau)$ is not much
different from $(P)$. The purpose of the equality constraint relaxation on the small intervals adjacent
to the time endpoints is to avoid {\em degeneracy} of necessary optimality conditions, which otherwise may
hold for all the feasible solutions under some choice of nontrivial dual elements.
Such a {\em degeneracy phenomenon} for necessary optimality conditions of the Pontryagin Maximum Principle (PMP)
type has been discovered and well investigated in control theory with {\em inequality} state constraints;
in particular, for Lipschitzian and compact-valued differential inclusions as in \cite{aa,rv,v}.
Our case is significantly different from the previous studies
in both directions of the problem setting and the results obtained. We derive {\em nondegenerate} necessary
optimality conditions for intermediate local minimizers of $(P^\tau)$ with $\tau>0$ while the passage to the
limit therein as $\tau\dn 0$ leads us to the conditions that generally degenerate. As examples show, even the
degenerate optimality conditions obtained in this way for $(P)$ can be useful to find optimal controls,
but anyway we treat as our main result the more trustworthy ones established for $(P^\tau)$.
Of course, there is no difference between problems $(P)$ and $(P^\tau)$ if the equality constraints \eqref{sw-con3} are not imposed.

Although the obtained necessary optimality conditions for the controlled sweeping process are constructively expressed via the problem data and turn out to be efficient as illustrated by various examples presented in this paper, they are rather complicated and contain measures, which is not surprising for state-constrained systems. However, a crucial advantage of the method of discrete approximations and its strong convergence established below is that we can stop, with {\em any prescribed accuracy}, at a suitable step of discretization and treat the corresponding optimal solution to the discrete problem satisfying the (much simpler) discrete optimality conditions as an {\em approximate/suboptimal} solution to the continuous-time one.\vspace*{0.02in}

The rest of the paper is organized as follows. The main result of Section~2 shows that the class of {\em absolutely continuous} controls $(u(\cdot),b(\cdot))$ and the corresponding absolutely continuous trajectories $x(\cdot)$ is a right choice for {\em feasible solutions} to $(P)$ and $(P^\tau)$, since such a control pair satisfying the polyhedral constraints \eqref{sw-con2} ensures the existence of an absolutely continuous solution to the Cauchy problem \eqref{sw-con1} under an appropriate constraint qualification, which is also used in deriving necessary optimality conditions.

The major aim of Section~3 is to construct well-posed discrete approximations of the optimal control problem $(P^\tau)$ for any $\tau\in[0,T]$, with
$P^0:=P$, such that they admit optimal solutions whose piecewise linear extensions on $[0,T]$ converges to the given intermediate local minimizer
$(\ox^\tau(\cdot),\ou^\tau(\cdot),\ob^\tau(\cdot))$ of $(P^\tau)$ {\em strongly in $W^{1,2}[0,T]$} with some {\em additional properties}
allowing us to derive nondegenerate necessary optimality conditions for $(\ox^\tau(\cdot),\ou^\tau(\cdot),\ob^\tau(\cdot))$ as $\tau>0$
by passing to the limit from discrete approximations. This essentially distinguishes our new results in this direction from those obtained
in the preceding paper \cite{chhm2} devoted to discrete approximations of the control sweeping process. As a crucial step of this procedure,
we justify the strong $W^{1,2}$-approximation with the desired additional properties for {\em any} feasible solution to the controlled
sweeping differential inclusion \eqref{sw-con1} without taking into account the cost functional \eqref{eq:MP}.

Since optimal control problems $(P^\tau)$ for $\tau\in[0,T]$ and its discrete counterparts are {\em intrinsically nonsmooth}
due to the sweeping dynamics \eqref{sw-con1} and its finite-difference approximations, we need to employ suitable
constructions of {\em generalized differentiation} satisfying extensive calculus rules to obtain necessary optimality
conditions first for discrete-time and then for continuous-time systems. Section~4 is devoted to the description of
such constructions and the explicit calculations of the major {\em second-order} one---the coderivative of the normal
cone mapping---entirely in {\em terms of the given data} of the controlled sweeping process.
These second-order calculations are certainly of their own interest while playing a crucial role in
the efficient realization of our approach to deriving necessary optimality conditions.

Section~5 presents {\em necessary optimality conditions} for the {\em discrete-time} optimal control problems
appearing in the discrete approximation procedure for $(P^\tau)$, $\tau\in[0,T]$, developed in Section~3.
These conditions are obtained by reducing the discrete-time problems to nonsmooth mathematical programs with
many geometric and functional constraints with the usage of generalized differential calculus and the
second-order calculations from Section~4. The conditions obtained are expressed explicitly via the problem data.

Section~6 is a {\em culmination} of the paper. It contains the formulation and proof of the main
{\em nondegenerate necessary conditions} for intermediate local minimizers of the sweeping control problem
$(P^\tau)$ whenever $\tau\in(0,T)$, their limiting versions as $\tau\dn 0$, and those for some special cases.
The proof of the main result is rather involved and significantly depends on the major results obtained in the previous sections.

The concluding Section~7 contains some applications to problems of {\em quasistatic elastoplasticity with hardening}
and also present several {\em examples} showing the strength and illustrating specific features of the necessary
optimality conditions obtained for the controlled sweeping process.

The notation of this paper is standard in variational analysis and optimal control; see, e.g., \cite{m-book,v}.
Recall that $\B$ stands for the closed unit ball of the space in question, $B(x,r):=x+r\B$, and $\N:=\{1,2,\ldots\}$.

\section{Feasible Solutions to the Controlled Sweeping Process}\label{sec:bdd}
\setcounter{equation}{0}

To begin our study, we want to make sure that the choice of absolutely continuous controls in problems $(P^\tau)$ as
$\tau\in[0,T]$ is appropriate from the viewpoint of {\em feasibility}, i.e., such a choice of $(u(\cdot),b(\cdot))$
in \eqref{sw-con2} ensures the existence of a solution $x(\cdot)$ to the Cauchy problem in \eqref{sw-con1},
which is at least {\em absolutely continuous} on $[0,T]$. Observe that the unbounded polyhedral moving set $C(t)$
generated by such a pair $(u(\cdot),b(\cdot))$ in \eqref{sw-con2} is not absolutely continuous on $[0,T]$ in the
Hausdorff sense for set-valued mappings (not even talking about Lipschitz continuity), and hence we cannot deduce
the existence of an absolutely continuous trajectory $x(\cdot)$ of \eqref{sw-con1} from known existence theorem
for the sweeping process; see, e.g., \cite{CT}.

Consider first the general sweeping process \eqref{sw-con1} generated by an arbitrary closed and convex moving set
$C(t)$ in $\R^n$, which is assumed to be nonempty for all $t\in[0,T]$. Denote by $v(t):=\pi_{C(t)}(0)$ the unique
{\em projection} of the origin onto $C(t)$ and define the shifted set $K(t):=C(t)-v(t)$.
The following result has been recently proved in \cite[Theorem~2.1]{chhm2} while being our starting point in this section.

\begin{Lemma}{\bf (existence of absolutely continuous sweeping trajectories for general moving sets).}\label{prop2}
Let the projection $v\colon[0,T]\to\R^n$ be absolutely continuous on $[0,T]$.
We assume that for any positive numbers $r,\ve$ there is a number $\delta=\delta(r,\ve)>0$ satisfying the estimate
\begin{equation}\label{tau}
\sum\limits_{i=1}^{l}\mathop{\rm max}\limits_{z\in K(\alpha_i)\cap r\B}{\rm dist}\big(z;K(\beta_i)\big)\le\ve.
\end{equation}
for every collection of mutually disjoint subintervals
\[
\big\{[\alpha_i,\beta_i]\big|\;i=1,\ldots,l\big\}\;\mbox{ of }\;[0,T]\;\mbox{ with }\;\disp\sum\limits_{i=1}^{l}|\beta_i-\alpha_i|\le\delta.
\]
Then there exists an absolutely continuous solution of the Cauchy problem in \eqref{sw-con1}.
\end{Lemma}

We now use this result to establish an existence theorem for absolutely continuous trajectories of \eqref{sw-con2}
generated by absolutely continuous controls in the {\em polyhedral} description \eqref{sw-con2} under the
{\em Linear Independence Constraint Qualification} (LICQ). Note that this qualification condition was missed
in the statement of \cite[Corollary~2.2]{chhm2}, where the proof was given only in the case of $m=1$ in \eqref{sw-con2}.
We are very grateful to Alexander Tolstonogov for observing that an additional condition is needed for the validity of
the latter result for $m>1$ and that the case of $m=1$ follows from his more general recent existence theorem in \cite{sansanych}.
The following result new for the case of $m>1$ is what we needed to justify {\em well-posedness} of the absolutely
continuous framework for feasible solutions to problems $(P^\tau)$ as $\tau\in[0,T]$.

\begin{Theorem} {\bf (existence of absolutely continuous sweeping trajectories for controlled polyhedra).}\label{ex-polyh}
Let the controls $(u(\cdot)=(u_1(\cdot),\ldots,u_m(\cdot))$ and $b(\cdot)=(b_1(\cdot),\ldots,b_m(\cdot))$ be
absolutely continuous on $[0,T]$, let the inequality system in \eqref{sw-con2} be consistent,
$($i.e., $C(t)\ne\emp))$ for all $t\in[0,T]$, and let
\begin{equation}\label{LICQ}
\text{the vectors $\big\{u_i(t)\big|\;i\in I_x(t)\big\}$ be linearly independent whenever }\;x\in C(t),\;t\in[0,T],
\end{equation}
where $I_x(t):=\{i\in\{1,\ldots,m\}\;\mbox{ with }\;\langle u_i (t),x\rangle=b_i(t)\}$.
Then the corresponding Cauchy problem in \eqref{sw-con1} admits a unique absolutely continuous solution $x(\cdot)$ on $[0,T]$.
\end{Theorem}
{\bf Proof.} Let us show first that the (well-defined) projection $v(\cdot)$ is absolutely continuous on $[0,T]$.
Indeed, observe that for each $t\in[0,T]$ the vector $v(t)$ solves a positive-definite parametric program under
the imposed LICQ \eqref{LICQ}. Then all the assumptions of Robinson's stability theorem from \cite[Theorems~2.1 and 4.1]{rob1} are satisfied.
Taking into account that $f$ appearing in the stability condition \cite[formula (2.4)]{rob1} is the Lagrangian of our
quadratic program, $x$ is our $v(t)$, $C$ is the first orthant, and the parameter $p$ therein is our $t$,
we conclude from the aforementioned stability condition that the modulus of continuity of $v(\cdot)$ is proportional
to the modulus of continuity of the problem data with respect to $t$. This readily justifies the claimed absolute continuity of the projection function.

To deduce the existence theorem for the sweeping process under consideration from Lemma~\ref{prop2},
it remains to verify the validity of condition (\ref{tau}) in this case.
Recall that a multifunction $\Gamma:\left[0,T\right]\rightrightarrows\mathbb{R}^{n}$ is absolutely continuous if for any $\varepsilon>0$ there
is some $\delta\left(\varepsilon\right)>0$ such that the implication
\[
\sum\nolimits_{i}\left\vert\beta_{i}-\alpha_{i}\right\vert<\delta\left(\varepsilon\right)
\Longrightarrow\sum\nolimits_{i}d_{H}\left(\Gamma\left(\alpha_{i}\right),
\Gamma\left(\beta_{i}\right)\right)<\varepsilon
\]
holds for any finite collection of mutually disjoint intervals $\left[\alpha_{i},\beta_{i}\right]\subset\left[0,T\right]$,
where $d_{H}$ refers to the Hausdorff distance. For $r>0$ and $t\in\left[0,T\right]$ define
$\Gamma_{r}(t):=K(t)\cap r\mathbb{B}$. Then it follows from the definition of $K(\cdot)$ that
$\Gamma_{r}(t)=\cap _{j=1}^{m}\Gamma _{r}^{(j)}(t)$ with the notation
\[
\Gamma_{r}^{(j)}(t):=\left(C_{j}(t)-v(t)\right)\cap r\mathbb{B}\;\mbox{ and }\;
C_{j}(t):=\big\{x\in\mathbb{R}^{n}\big|\left\langle u_{j}(t),x\right\rangle\le b_{j}(t)\big\},\quad j=1,\ldots,m.
\]
We claim now the validity of the estimate
\begin{eqnarray}\label{beh1}
\mathrm{dist}\,(z,\Gamma _{r}^{(j)}(t))\le r\left\Vert u_{j}(t)-u_{j}(s)\right\Vert+\left\vert b(t)-b(s)\right\vert+\left\vert
\left\langle u_{j}(t),v(t)\right\rangle-\left\langle u_{j}(s),v(s)\right\rangle\right\vert
\end{eqnarray}
for all $s,t\in[0,T]$, $r>0$, $j\in\left\{1,\ldots,m\right\}$ and $z\in\Gamma_{r}^{(j)}(s)$.
Note that (\ref{beh1}) is trivially satisfied for $z\in\Gamma_{r}^{(j)}(t)$, and hence we
may assume that $z\notin\Gamma_{r}^{(j)}(t)$. Now fix such $s,t,r,j,z$ and show first that
\begin{eqnarray}\label{beh2}
\mathrm{dist}\,(z,\Gamma_{r}^{(j)}(t))=\gg:=\left\langle u_{j}(t),z+v(t)\right\rangle-b_{j}(t)\ge 0.
\end{eqnarray}
Indeed, \eqref{beh2} follows from the easily verifiable fact that $\tilde{z}:=z-\gg u_{j}(t)$ is a unique minimizer of the program
\begin{equation}\label{prog1}
\min_{y}\big\{\left\Vert z-y\right\Vert^{2}\big|\;y\in C_{j}(t)-v(t)\big\}
\end{equation}
by $\left\Vert u_{j}(t)\right\Vert=1$. Observing that $v(t)\in C(t)\subset C_{j}(t)$ and $\gg\ge 0$ gives us
$\gg\le 2\left\langle u_{j}(t),z\right\rangle$, and so
\[
\left\Vert \tilde{z}\right\Vert^{2}=\left\Vert z\right\Vert^{2}-2\gg\left\langle u_{j}(t),z\right\rangle+\gg^{2}=
\left\Vert z\right\Vert^{2}+\gg\left(\gg-2\left\langle u_{j}(t),z\right\rangle\right)\le\left\Vert z\right\Vert^{2}\le r^{2},
\]
which means that $\tilde{z}\in r\mathbb{B}$ and thus $\tilde{z}$ solves not just (\ref{prog1}) but also the program
\[
\min_{y}\big\{\left\Vert z-y\right\Vert^{2}\big|\;y\in\Gamma_{r}^{(j)}(t)\big\}
\]
over the smaller constraint set. This verifies (\ref{beh2}). Further, it follows from $z\in\Gamma _{r}^{(j)}(s)$ that
\begin{eqnarray*}
\gg &=&\left\langle u_{j}(t)-u_{j}(s),z\right\rangle+\left\langle u_{j}(s),z\right\rangle+\left\langle u_{j}(t),v(t)\right\rangle -b_{j}(t)\\
&\le&\left\langle u_{j}(t)-u_{j}(s),z\right\rangle+b_{j}(s)-\left\langle u_{j}(s),v(s)\right\rangle+\left\langle u_{j}(t),v(t)\right\rangle-b_{j}(t),
\end{eqnarray*}
which implies (\ref{beh1}) by taking into account that $\left\Vert z\right\Vert\le r$. Interchanging the roles of $s$ and $t$ in (\ref{beh1}) gives us
\begin{eqnarray*}
d_{H}(\Gamma_{r}^{(j)}(s),\Gamma_{r}^{(j)}(t))\le r\left\Vert u_{j}(t)-u_{j}(s)\right\Vert+\left\vert b(t)-b(s)\right\vert+\left\vert
\left\langle u_{j}(t),v(t)\right\rangle-\left\langle u_{j}(s),v(s)\right\rangle\right\vert
\end{eqnarray*}
for the same $s,t,r,j,z$. On the other hand, the functions $ru_{j}(\cdot)$, $b(\cdot)$, and $\left\langle u_{j}(\cdot),v(\cdot)\right\rangle$ are absolutely
continuous on $[0,T]$ since $u_{j}(\cdot)$ and $b(\cdot)$ were assumed while $v$ was shown to be such. This tells us that the multifunction
$\Gamma_{r}^{(j)}(\cdot)$ is absolutely continuous for any $r>0$ and $\,j\in\left\{1,\ldots,m\right\}$. We want to derive from here that $\Gamma_{r}(\cdot)$ is
absolutely continuous for any $r>0$, which would follow from this property of the intersection mapping $\Gamma_{r}^{(1)}(\cdot)\cap\Gamma_{r}^{(2)}(\cdot)$.
Since this intersection is bounded, the desired fact follows from
\begin{equation} \label{beh3}
\Gamma_{r}^{(1)}\left(t\right)\cap\mathrm{int}\,\Gamma_{r}^{(2)}\left(t\right)\ne\emp\;\mbox{ whenever }\;t\in\left[0,T\right]
\end{equation}
by \cite[Proposition on p.\ 274]{jjmoreau75}. To verify (\ref{beh3}), observe first that the assumed LICQ and nonemptiness of $C(t)$ yields the existence of a
(time dependent) Slater point $\tilde{x}(t)$ in the description of $C(t)$:
\[
\left\langle u_{j}(t),\tilde{x}(t)\right\rangle<b_{j}(t)\;\mbox{ for all }\;t\in\left[0,T\right],\;j=1,\ldots,m.
\]
Recall the inclusion $v(t)\in C_{j}(t)$ for all $j=1,\ldots,m$ on $[0,T]$, which means that
\[
\left\langle u_{j}(t),v(t)\right\rangle\le b_{j}(t)\;\mbox{ whenever }\;t\in\left[0,T\right],\;j=1,\ldots,m.
\]
Define $x_{\nu}(t):=\nu\tilde{x}(t)+\left(1-\nu\right)v(t)$ for $t\in[0,T]$ and $\nu\in\left[0,1\right]$ and deduce that
$\left\langle u_{j}(t),x_{\nu}(t)\right\rangle<b_{j}(t)$ and $\nu\in\mathrm{int}\,r\mathbb{B}$ as $t\in\left[0,T\right]$, $j=1,2$, $\nu\in(0,\min\{r,1\})$.
It yields $x_{\nu}(t)\in\Gamma_{r}^{(1)}\left(t\right)\cap\mathrm{int}\,\Gamma _{r}^{(2)}\left(t\right)$ justifying (\ref{beh3}).

It is proved therefore that $\Gamma_{r}(\cdot)$ is absolutely continuous on $[0,T]$ for any $r>0$, which means that whenever $r,\varepsilon>0$
there is $\delta\left(r,\varepsilon\right)>0$ such that the implication
\[
\sum\nolimits_{i}\left\vert\beta _{i}-\alpha_{i}\right\vert<\delta
\left(r,\varepsilon\right)\Longrightarrow\sum\nolimits_{i}d_{H}\left(\Gamma_{r}\left(\alpha_{i}
\right),\Gamma_{r}\left(\beta_{i}\right)\right)<\varepsilon
\]
holds for any finite collection of mutually disjoint intervals $\left[\alpha_{i},\beta_{i}\right]\subset\left[0,T\right]$. This verifies that
\begin{eqnarray*}
\sum\limits_{i=1}^{l}\max_{z\in K\left(\alpha_{i}\right)\cap r\mathbb{B}}\mathrm{dist}\,(z,K\left(\beta_{i}\right))\le\sum\limits_{i=1}^{l}
\max_{z\in K\left(\alpha_{i}\right)\cap r\mathbb{B}}\mathrm{dist}\,(z,K\left(\beta_{i}\right)\cap r\mathbb{B})\le\sum\limits_{i=1}^{l}d_{H}
\left(\Gamma_{r}\left(\alpha_{i}\right)\Gamma_{r}\left(\beta_{i}\right)\right)<\varepsilon
\end{eqnarray*}
and thus completes the existence part of the proof of the theorem. Uniqueness follows from a well known argument based on the convexity of the moving sets and
Gronwall's lemma. $\h$\vspace*{0.02in}

The next example demonstrates that just the consistency condition $C(t)\ne\emp$ is not sufficient for the
existence of absolutely continuous trajectories in \eqref{sw-con1}, \eqref{sw-con2} generated by
${\cal C}^\infty$ controls $(u(\cdot),b(\cdot))$ for $m=3$.

\begin{Example}{\bf (LICQ is essential for the existence of absolutely continuous sweeping trajectories).}
{\rm Consider the controlled sweeping system \eqref{sw-con1}, \eqref{sw-con2} in $\R^2$ generated by the controls
\[
u_1(t):=e_1,\;u_2(t):=-e_1,\;u_3(t):=\big(-\cos t,-\sin t\big),\;b_1(t)=1, b_2(t):=-1, b_3(t):=-\cos t-\sin t
\]
on $[0,\pi]$, which are obviously ${\cal C}^\infty$ functions on this intervals while LICQ \eqref{LICQ} fails. Then we have
\[
C(t)=\begin{cases}
\{1\}\times\R&\text{for }\;t=0,\\
\{1\}\times[1,\infty)&\text{for }\;0<t<\pi
\end{cases}\quad\mbox{and}\quad
v(t)=\begin{cases}
(1,0)&\text{for }\;t=0,\\
(1,1)&\text{for }\;0<t<\pi.
\end{cases}
\]
Due to the discontinuity of $v(t)$, the assumptions of both Lemma~\ref{prop2} and Theorem~\ref{ex-polyh}
are not satisfied although $C(t)\ne\emp$ for all $t\in[0,\pi]$. Observe further that the corresponding
Cauchy problem \eqref{sw-con1} with $x(0)=(1,0)$ cannot have absolutely continuous solutions, since even
the requirement $x(t)\in C(t)$ on $[0,T]$ is not met for continuous functions $x(t)$ by the discontinuity of $C(t)$.}
\end{Example}

\begin{Remark}{\bf (Lipschitzian sweeping trajectories).}\label{polyhedral}
{\rm It follows from \cite[Proposition~2.3]{chhm2} that the Lipschitz continuity assumptions on
$v(\cdot)$ and $K(\cdot)\cap r\B$ in addition to those imposed in Lemma~\ref{prop2} ensure
the existence of a {\em unique Lipschitzian} solution to \eqref{sw-con1} with an explicit estimate
of its Lipschitz constant. Then the arguments in the proof of Theorem~\ref{ex-polyh} allow us to
conclude under these assumptions that any Lipschitzian control pair $(u(\cdot),b(\cdot))$ in \eqref{sw-con2}
generates a unique solution to \eqref{sw-con1} with the aforementioned properties.
Furthermore, if all the controls $u(\cdot)$ and $b(\cdot)$ in \eqref{sw-con2} are {\em uniformly}
Lipschitzian on $[0,T]$ with the moduli $L_u$ and $L_b$, respectively, then we can prove the existence
of a constant $M>0$ dependent only on $T$, $x_0$, $L_u$, and $L_b$ such that the sweeping process in \eqref{sw-con1}
is equivalent to the {\em bounded} differential inclusion
\begin{eqnarray*}
-\dot{x}(t)\in N\big(x(t);C(t)\big)\cap M\B\;\mbox{ a.e. }\;t\in[0,T]\;\mbox{ with }\;x(0)=x_0\in C(0).
\end{eqnarray*}
Such a boundedness reduction was first established in \cite{thibault} for the sweeping process with an absolutely continuous moving set $C(t)$, which is not
the case here, and then was further developed in \cite{chhm} for \eqref{sw-con1}, \eqref{sw-con2} with $m=1$ under uniform Lipschitzian assumptions.
In this paper, in contrast to \cite{chhm}, we prefer not to require the uniform Lipschitzness of feasible controls that led us in \cite{chhm} to incomplete results on discrete approximations and necessary optimality conditions.}
\end{Remark}

\section{Well-Posed Discrete Approximations}\label{sect:discr}
\setcounter{equation}{0}

In this section we start developing a {\em discretization approach} to the study of the optimal control problems
$(P^\tau)$, $\tau\in[0,T]$, which will finally result in deriving robust necessary optimality conditions for their
local minimizers. As the first step of this device, we construct here well-posed discrete approximations for an
arbitrary {\em feasible} solution to $(P^\tau)$ and establish an appropriate {\em strong convergence} of such approximations.

To begin with, let us represent the controlled sweeping differential inclusion \eqref{sw-con1} under the polyhedral
constraints \eqref{sw-con2} in an equivalent differential inclusion form for extended trajectories.
Consider the vectors $u:=(u_1,\ldots,u_m)\in\R^{nm}$, $b:=(b_1,\ldots,b_m)\in\R^m$, and
$z:=(x,u,b)\in\R^n\times\R^{nm}\times\R^m$ and define the set-valued mapping $F\colon\R^n\times\R^{nm}\times\R^m\tto\R^n$ by
\begin{equation}\label{F}
F(z):=-N\big(x;C(u,b)\big)\;\mbox{ with }\;C(u,b):=\big\{x\in\R^n|\;\la u_i,x\ra\le b_i,\;i=1,\ldots,m\big\}.
\end{equation}
Then the sweeping differential inclusion in \eqref{sw-con1} under constraints \eqref{sw-con2} can be rewritten as
\begin{equation}\label{G}
\mathop z\limits^.(t)\in G\big(z(t)\big):=F\big(z(t)\big)\times\R^{nm}\times\R^m\;\mbox{ a.e. }\;t\in[0,T],
\end{equation}
where the initial condition $z(0)=(x_0,u(0),b(0))$ is such that $\la u_i(0),x_0\ra\le b_i(0)$ for all $i=1,\ldots,m$. We can see that form \eqref{G} treats controls $(u(\cdot),b(\cdot))$ and the corresponding sweeping trajectories $x(\cdot)$ symmetrically. Theorem~\ref{ex-polyh} tells us that any choice of absolutely continuous functions $u(\cdot)$ and $b(\cdot)$ satisfying LICQ \eqref{LICQ} generates a feasible solution $z(\cdot)$ for the extended differential inclusion \eqref{G}.

Now the constrained sweeping system \eqref{sw-con1}, \eqref{sw-con2} is written in the conventional form of the theory of differential inclusions with fixed right-hand sides as in \cite{m-book1,v} while they do not possess major properties under which this theory has been developed. Indeed, the right-hand side of
\eqref{G} is described by a {\em highly irregular} (discontinuous and unbounded) set-valued mapping. Further, it follows from definition \eqref{nor} of the normal cone that \eqref{G} implicitly contains the {\em state constraints} on the trajectories $z(t)$ given by
\begin{equation}\label{inequastate}
x(t)\in C\big(u(t),b(t)\big)\;\mbox{ for all }\;t\in[0,T].
\end{equation}
Moreover, besides the more conventional state constraints of the inequality type \eqref{sw-con2}, each optimal control problem $(P^\tau)$ as $\tau\in[0,T]$ contains those (nonsmooth and irregular) of the {\em equality} type given by \eqref{delta}, which has not been investigated in optimal control theory even in much simpler settings; see the discussion in Section~1. All of this emphasizes serious challenges we face to study these problems.

The next theorem establishes the aforementioned {\em strong $W^{1,2}$-approximation} of a given feasible solution to \eqref{G} subject to the state constraints in \eqref{delta} (those in \eqref{inequastate} are contained in \eqref{G}) by a sequence of feasible solutions to its discrete counterparts. The underlying difference of this theorem from the previous one in \cite[Theorem~3.1]{chhm2} is deriving, under additional assumptions, new approximation properties that play a crucial role in the subsequent passage to the limit from optimality conditions for discrete-time control problems. We check then that the additional assumptions imposed are not actually restrictive.\vspace*{0.03in}

In what follows the symbol $j_\tau(k)$ stands for the smallest natural number $j$ such that $t_j^k\ge\tau$ whenever $\tau\in[0,T]$, while $j^\tau(k)$ signifies the largest $j\in\N$ with $t^k_j\le T-\tau$.

\begin{Theorem}{\bf (strong discrete approximation of feasible solutions).}\label{feasible} Let $\oz(\cdot)=(\ox(\cdot),\ou(\cdot),\bar b(\cdot))\in W^{1,2}[0,T]:=
W^{1,2}([0,T];\R^{n+nm+m})$ be an arbitrary feasible solution to problem $(P^\tau)$ with any fixed parameter $\tau\in[0,T]$ and
define the uniform discrete partitions of $[0,T]$ by setting
\begin{equation}\label{partition}
\Delta_k:=\big\{0=t^k_0<t^k_1<\ldots<t^k_k=T\big\}\;\mbox{ with }\;h_k:=t^k_{j+1}-t^k_j\dn 0,\;j=0,\ldots,k-1,\;\mbox{ as }\;k\to\infty.
\end{equation}
Assume in addition that $\oz(\cdot)$ satisfies the following properties at the mesh points $($all these properties are automatically satisfied if, e.g., $\oz(\cdot)\in W^{2,\infty}[0,T])$: the differential inclusion \eqref{G} holds for $\oz(\cdot)$ at all
$t^k_j$, $j=0,\ldots ,k-1$ for each $k\in\N$, we have
\begin{equation}\label{intconv}
\sum^{k-1}_{j=0}h_k\Big\|\frac{\ox(t^k_{j+1})-\ox(t^k_j)}{h_k}-\dot\ox(t^k_j)\Big\|^2\to 0\;\mbox{ as }\;k\to\infty,
\end{equation}
and there exists a constant $M>0$ independent of $k$ such that
\begin{equation}\label{discretbound1}
\sum^{k-1}_{j=0}\Big\|\frac{\ox(t^k_{j+1})-\ox(t^k_j)}{h_k}-\dot\ox(t^k_j)\Big\|\le M,\quad
\Big\|\frac{\ou(t^k_{1})-\ou(t^k_0)}{h_k}\Big\|\le M,\quad \Big\|\frac{\ob(t^k_{1})-\ob(t^k_0)}{h_k}\Big\|\le M,
\end{equation}
\begin{equation}\label{discretebound}
\sum^{k-2}_{j=0}\Big\|\frac{\ou(t^k_{j+2})-\ou(t^k_{j+1})}{h_k}-\frac{\ou(t^k_{j+1})-\ou(t^k_{j})}{h_k}
\Big\|\le M,\quad\sum^{k-2}_{j=0}\Big\|\frac{\ob(t^k_{j+2})-\ob(t^k_{j+1})}{h_k}-\frac{\ob(t^k_{j+1})-\ob(t^k_{j})}{h_k}
\Big\|\le M
\end{equation}
whenever $k\in\N$. Then there is a sequence of piecewise linear functions $z^k(t):=(x^k(t),u^k(t),b^k(t))$ on $[0,T]$ with $\big(x^k(0),u^k(0),b^k(0)\big)=\big(x_0,\ou(0),\ob(0)\big)$,
\begin{eqnarray}\label{contr-disc}
\left\{\begin{array}{ll}
\|u_i^k(t^k_j)\|=1\;\mbox{ for }\;j=j_\tau(k),\ldots,j^\tau(k),\\
\disp\frac{1}{2}\le\|u^k_i(t^k_j)\|\le\frac{3}{2}\;\mbox{ for }\;j\le j_\tau(k)-1\;\mbox{ and }\;j\ge j^\tau(k)+1,\quad i=1,\ldots,m,
\end{array}\right.
\end{eqnarray}
satisfying for all $k\in\N$ the extended finite-difference inclusions
\begin{equation}\label{discrete-inclusion1}
x^k(t)=x^k(t_j)+(t-t_j)v^k_j\;\mbox{ whenever }\;t_j^k\le t\le t_{j+1}^k\;\mbox{ with }\;v^k_j\in F\big(z^k(t_j^k)\big),\quad j=0,\ldots,k-1,
\end{equation}
and such that the functions $z^k(\cdot)$ converge to $\oz(\cdot)$ in the norm topology of $W^{1,2}[0,T]$, i.e.,
\begin{equation}\label{zk}
z^k(t)\to\oz(t)\;\;uniformly\;\;on\;\;[0,T]\;\;and\;\;\int\limits_0^T\|\dot{z}^k(t)-\dot\oz(t)\|^2\,dt\to 0\;\mbox{ as }\;k\to\infty.
\end{equation}
Moreover, there exists a constant $\tilde{M}\ge M$ depending only on the total variation of $\ou(\cdot)$ on $[0,T]$ so that for every $k\in\N$ we have the estimates
\begin{equation}\label{initialboundedub}
\Big\|\frac{u^k(t^k_{1})-u^k(t^k_0)}{h_k}\Big\|\le \tilde{M},\quad\Big\|\frac{b^k(t^k_{1})-b^k(t^k_0)}{h_k}\Big\|\le \tilde{M},
\end{equation}
\begin{equation}\label{uniformbounded}
{\rm var}\big(\dot u^k;[0,T]\big)\le \tilde{M},\;\mbox{ and }\;{\rm var}\big(\dot b^k;[0,T]\big)\le \tilde{M}.
\end{equation}
Finally, the sequence $\{x^k(\cdot)\}$ has uniformly bounded variations on $[0,T]$ being in
addition uniformly Lipschitzian on $[0,T]$ if $\ox(\cdot)$ is Lipschitz continuous on this interval.
\end{Theorem}
{\bf Proof.} Following the proof of \cite[Theorem~3.1]{chhm2}, for any $k\in\N$ we construct the continuous-time functions $y^k(t):=\big(y^k_1(t),y^k_2(t),y^k_3(t)\big)$ as the piecewise linear extensions on $[0,T]$ of the discrete-time triples
\begin{equation*}
\big(y^k_1(t^k_j),y^k_2(t^k_j),y^k_3(t^k_j)\big):=\big(\ox(t^k_j),\ou(t^k_j),\ob(t^k_j)\big),\quad j=0,\ldots,k,
\end{equation*}
and denote by $w^k(t)=\big(w_1^{k}(t),w_2^{k}(t),w_3^{k}(t)\big):=\dot{y}^k(t)$ their derivatives at non-mesh points.
The assumptions in \eqref{discretebound} yield ${\rm var}(w_i^k;[0,T])\le M$ for $i=2,3$ whenever $k\in\N$. It follows from the above that
\[
y^k(\cdot)\to\oz(\cdot)\;\mbox{ uniformly on }\;[0,T]\;\mbox{ and }\;w^k(\cdot)\to\dot{\oz}(\cdot)\;\mbox{ strongly in }\;L^2[0,T]\;\mbox{ as }\;k\to\infty.
\]
Defining $u^k(t):=y^k_2(t)$ on $[0,T]$, we get $u^k(0)=y^k_2(0)=\ou(0)$ and deduce from \eqref{delta}
and the constructions above that the constraints on $u^k(\cdot)$ in \eqref{contr-disc} and \eqref{initialboundedub} hold.
Fix $k\in\N$, denote $t_j:=t^k_j$ as $j=1,\ldots,k-1$, and construct the desired functions $b^k(t),x^k(t)$ on $[0,T]$ by induction.
To proceed, put $\big(x^k(0),b^k(0)\big)=\big(x_0,\ob(0)\big)$, suppose that the value of $x^k(t_j)$ is known, and define $b^k_i(t)$ at the mesh points so that
\[
\begin{split}
b_i^k(t_j):=\big\la x^k(t_j),u_i^k(t_j)\big\ra&\mbox{if }\;y_{3i}^k(t_j)=\big\la{y_1^k(t_j),y_{2i}^k(t_j)}\big\ra,\\
b_i^k(t_j)>\big\la x^k(t_j),u_i^k(t_j)\big\ra&\mbox{if }\;y_{3i}^k(t_j)>\big\la{y_1^k(t_j),y_{2i}^k(t_j)}\big\ra.
\end{split}
\]
We can clearly arrange $b_i^k(t_j)-\ob_i(t_j)=b_i^k(t_j)-y_{3i}^k(t_j)=\la x^k(t_j),u_i^k(t_j)\ra-\la{y_1^k(t_j),y_{2i}^k(t_j)}\ra=\la x^k(t_j)-y_1^k(t_j),\ou_{i}(t_j)\ra$. Using the projection $v_j^k=\pi_{F(x^k(t_j),u^k(t_j),b^k(t_j))}(w_{1j}^k)$, define next $x^k(t)$ on $(t_j,t_{j+1}]$ by
\eqref{discrete-inclusion1}, then construct $b^k_i(t_{i+1})$ as above and extend it linearly to $[t_j,t_{j+1}]$. Observe that our construction yields
\begin{equation}\label{bv}
\begin{split}
&\disp\Big|\frac{b_i^k(t_{j+1})-b_i^k(t_j)}{h_k}-\frac{\ob_i(t_{j+1})-\ob_i(t_j)}{h_k}\Big|=
\disp\Big|\frac{b_i^k(t_{j+1})-\ob(t_{j+1})}{h_k}-\frac{b_i^k(t_{j})-\ob(t_j)}{h_k}\Big|\\
&\qquad\qquad =\disp\Big|\Big\la\frac{x^k(t_{j+1})-\ox(t_{j+1})}{h_k},\ou_i(t_{j+1})\Big\ra-
\Big\la\frac{x^k(t_{j})-\ox(t_{j})}{h_k},\ou_i(t_{j})\Big\ra\Big|\\
&\qquad\qquad \le\disp\frac{3}{2}\Big\|\frac{x^k(t_{j+1})-x^k(t_j)}{h_k}-\frac{\ox(t_{j+1})-\ox(t_j)}{h_k}
\Big\|\disp+\Big\|\frac{x^k(t_{j})-\ox(t_{j})}{h_k}
\Big\|\cdot\big\|\ou_i(t_{j+1})-\ou_i(t_{j})\big\|
\end{split}
\end{equation}
for all $j=0,\ldots,k-1$. We also have the equalities
\begin{equation}\label{bv2}
\Big\|\frac{x^k(t_{j+1})-x^k(t_j)}{h_k}-\frac{\ox(t_{j+1})-\ox(t_j)}{h_k}\Big\|=
\Big\|v^k_j-\frac{\ox(t_{j+1})-\ox(t_j)}{h_k}\Big\|=\|v^k_j-w_{1j}^k\|.
\end{equation}
Furthermore $F\big(x^k(t_j),u^k(t_j),b^k(t_j)\big)=F(\ox(t_j),\ou(t_j),\ob(t_j))$, $v_j^k=\pi_{F(\ox(t_j),\ou(t_j),\ob(t_j))}(w_{1j}^k)$, and
\begin{equation}\label{v-w}
\begin{split}
\sum^{k-1}_{j=0}\|v^k_j-w^k_{1j}\|&=\sum^{k-1}_{j=0}{\rm dist}\Big(\frac{\ox(t_{j+1})-\ox(t_j)}{h_k};F(\ox(t_j),\ou(t_j),\ob(t_j))\Big)\\
&\le\sum^{k-1}_{j=0}\Big\|\frac{\ox(t_{j+1})-\ox(t_j)}{h_k}-\dot\ox(t_j)\Big\|\le M.
\end{split}
\end{equation}
Employing this together with \eqref{bv2} gives us
\begin{equation}\label{liplip}
\sum^{k-1}_{j=0}\Big\|\frac{x^k(t_{j+1})-x^k(t_j)}{h_k}-\frac{\ox(t_{j+1})-\ox(t_j)}{h_k}\Big\|\le M,
\end{equation}
which readily implies the estimates
\begin{equation}\label{bar-k}
\Big\|\frac{x^k(t_{j})-\ox(t_j)}{h_k}\Big\|\le
\sum^{j-1}_{i=0}\Big\|\frac{x^k(t_{i+1})-\ox(t_{i+1})}{h_k}-\frac{x^k(t_{i})-\ox(t_i)}{h_k}\Big\|+\Big\|\frac{x^k(t_{0})-\ox(t_0)}{h_k}\Big\|\le M
\end{equation}
for every $j=0,\ldots,k-1$ as $k\in\N$. As a result, it follows from \eqref{bv}, \eqref{bv2}, \eqref{v-w}, and \eqref{bar-k} that
\[
\sum^{k-1}_{j=0}\Big|\frac{b_i^k(t_{j+1})-b_i^k(t_j)}{h_k}-\frac{\ob_i^k(t_{j+1})-\ob_i^k(t_j)}{h_k}\Big|\le M+M{\rm var}\big(\ou;[0,T]\big).
\]
Combining the latter with \eqref{discretebound}, we arrive at
\[
\sum^{k-2}_{i=0}\Big|\frac{b^k(t_{j+2})-b^k(t_{j+1})}{h_k}-\frac{b^k(t_{j+1})-b^k(t_{j})}{h_k}\Big|\le 3M+2M{\rm var}\big(\ou;[0,T]\big),
\]
with verifies the validity of \eqref{uniformbounded} for $\dot{b}^k$. Observe simultaneously the fulfillment of the estimate for $b^k$
in \eqref{initialboundedub}, which follows from \eqref{bar-k} and the representation $b_i^k(t_1)-\ob_i(t_1)=\la x^k(t_1)-y_1^k(t_1),\ou_{i}(t_1)\ra$.

Next we justify the $W^{1,2}$-convergence in \eqref{zk} for which it suffices in fact to check the
$L^2$-convergence of $\dot{x}^k$ and $\dot{b}^k$. To verify the former one, observe by \eqref{intconv} that
\begin{equation*}
\int^T_0\|\dot x^k(t)-w^k_1(t)\|^2dt=\sum^{k-1}_{j=0}h_k\|v_{j}^k-w^k_j\|^2\le\sum^{k-1}_{j=0}h_k\Big\|\frac{\ox(t^k_{j+1})-\ox(t^k_j)}{h_k}-\dot \ox(t^k_j)\Big\|^2\to 0
\end{equation*}
as $k\to\infty$. The claimed convergence for $\dot{b}^k$ follows from \eqref{bv} and
\begin{equation*}
\int^T_0\|\dot b^k(t)-w^k_3(t)\|^2dt\le\frac92\int^T_0\|\dot x^k(t)-w^k_1(t)\|^2dt+2M^2\sum^{k-1}_{j=0}h_k\|\ou_i(t_{j+1})-\ou_i(t_{j})\|^2
\end{equation*}
by the absolute continuity of $\ou(\cdot)$ on $[0,T]$. The last statement of the theorem on $\{x^k(\cdot)\}$ follows
immediately from \eqref{discretbound1}, \eqref{liplip}, and the fact that $\ox(\cdot)$ has bounded
variation on $[0,T]$. To complete the proof of the theorem, it remains to observe that the validity of all the
assumptions in \eqref{intconv}--\eqref{discretebound} for the case of $\oz(\cdot)\in W^{2,\infty}[0,T])$ is a direct
consequence of the definitions. $\h$\vspace*{0.03in}

It is not hard to check that all the assumptions of Theorem~\ref{feasible} are satisfied if $\oz(\cdot)$ is
{\em piecewise ${\cal C}^1[0,T]$ with $\dot{\oz}(\cdot)\in BV([0,T])$}. In this case, more general than
$\oz(\cdot)\in W^{2,\infty}[0,T]$, the derivatives appearing in \eqref{G} and in the subsequent formulas are the right derivatives.\vspace*{0.05in}

Our next goal is to construct a well-posed discrete approximation for a {\em given} local optimal solution of the control problem
$(P^\tau)$ as $\tau\in[0,T]$,
which satisfies the assumptions of Theorem~\ref{feasible}. We consider a rather broad class of local
minimizers introduced and first studied for differential inclusions in \cite{m95} under the name of
{\em intermediate local minimizers} (i.l.m.). This notion obviously covers strong local minimizers
(corresponding to $\al=0$ in the definition below) and occupies an intermediate position between weak
and strong minimizers in dynamic optimization and optimal control; see \cite{m95} and \cite[Chapter~6]{m-book1}
for more details. Note that a related class of local minimizers in optimal control problems for differential
inclusions has been studied later under the name of $W^{1,1}$-minimizers; see, e.g., \cite{v}.
We now present an adaptation of the i.l.m.\ notion to the case of the sweeping control problems $(P^\tau)$ under consideration.

\begin{Definition}{\bf (intermediate local minimizers for the controlled sweeping process).}\label{ilm}
Fix any $\tau\in[0,T]$ and consider a feasible solution
$\oz(\cdot)=\big(\ox(\cdot),\ou(\cdot),\ob(\cdot)\big)\in W^{1,2}[0,T]$ to $(P^\tau)$.
We say that $\oz(\cdot)$ is an {\sc intermediate local minimizer} for this problem if there are numbers $\al\ge 0$ and $\ve>0$ such that $J[\oz]\le J[z]$
for any feasible solution $z(\cdot)=\big(x(\cdot),u(\cdot),b(\cdot)\big)$ to $(P^\tau)$ satisfying
\begin{equation*}
\big\|\big(x(t),u(t),b(t)\big)-\big(\ox(t),\ou(t),\ob(t)\big)\big\|<\ve\;\mbox{ as }\;t\in[0,T],\;\mbox{ and}
\end{equation*}
\begin{equation}\label{interm}
\al\int\limits_{{0}}^{{T}}\Big(\Big\|\dot x(t)-\dot\ox(t)\Big\|^2+\Big\|\dot u(t)-\dot\ou(t)\Big\|^2+\Big\|\dot b(t)-\dot\ob(t)\Big\|^2\Big)\,dt<\ve.
\end{equation}
\end{Definition}

It is easy to see that the general setting of $\al\ge 0$ in \eqref{interm} reduces to the cases when either
$\al=1$ or $\al=0$, and that from the viewpoint of necessary optimality conditions it suffices to examine
only the case of intermediate local minimizers with $\al=1$, which we do in what follows.\vspace*{0.03in}

Given a number $\tau\in[0,T]$ and an i.l.m.\ $\oz(\cdot)$ for $(P^\tau)$, construct now the family of
{\em discrete approximation problems} $(P_k^\tau)$, $k\in\N$, having optimal solutions that converge to
$\oz(\cdot)$ in some strong sense, which eventually allows us to derive eventually necessary optimality
conditions for $\oz(\cdot)$ by passing to the limit from those in discrete approximations.
Suppose that for the given i.l.m.\ $\oz(\cdot)$ and the discrete mesh $\Delta_k$ in \eqref{partition} all the assumptions (and hence conclusions) of
Theorem~\ref{feasible} are satisfied and then define each problem $(P_k^\tau)$ by:
\begin{eqnarray}\label{cost-pk}
\begin{array}{ll}
\mbox{minimize}\;\;
J_k[z^k]:=&\disp\varphi(x_k^k)+\disp{h_k}\sum\limits_{j=0}^{k-1}\disp{\ell\Big(t^k_j,x_j^k,u_j^k,b_j^k,\frac{x_{j+1}^k-x_j^k}{h_k},
\disp\frac{u_{j+1}^k-u_j^k}{h_k},\frac{b_{j+1}^k-b_j^k}{h_k}\Big)}\\
&+\sum\limits_{j=0}^{k-1}\int\limits_{{t^k_j}}^{{t^k_{j+1}}}\disp\Big(\Big\|\frac{x_{j+1}^k-x_j^k}{h_k}-\dot\ox(t)\Big\|^2+
\Big\|\disp\frac{u_{j+1}^k-u_j^k}{h_k}-\dot{\ou}(t)\Big\|^2+\Big\|\frac{b_{j+1}^k-b_j^k}{h_k}-\dot{\ob}(t)\Big\|^2\Big)\,dt\\
&\disp+\,{\rm dist}^2\Big(\Big\|\frac{u_{1}^k-u_{0}^k}{h_k}\Big\|,(-\infty,\tilde{M}]\Big)+
{\rm dist}^2\disp\Big(\sum_{j=0}^{k-2}\Big\|\frac{b_{1}^k-b_{0}^k}{h_k}\Big\|,(-\infty,\tilde{M}]\Big)\\
&\disp+\,{\rm dist}^2\Big(\sum_{j=0}^{k-2}\Big\|\frac{u_{j+2}^k-2u_{j+1}^k+u_{j}^k}{h_k}\Big\|,
(-\infty,\tilde{M}]\Big)\\
&\disp +{\rm dist}^2\disp\Big(\sum_{j=0}^{k-2}\Big\|\frac{b_{j+2}^k-2b_{j+1}^k+b_{j}^k}{h_k}\Big\|,(-\infty,\tilde{M}]\Big)
\end{array}
\end{eqnarray}
over collections $z^k:=(x_0^k,\ldots,x_k^k,u_0^k,\ldots,u_k^k,b_0^k,\ldots,b_k^k)$ with $u_j^k:=(u_{j1}^k,\ldots,u_{jm}^k)\in\R^{nm}$ for every index $j=0,\ldots,k$
subject to the constraints in \eqref{contr-disc} together with
\begin{equation}\label{discrete-inclusion}
x_{j+1}^k\in x_j^k+h_k F(x_j^k,u_j^k,b_j^k)\;\mbox{ for }\;j=0,\ldots,k-1\;\mbox{ with }\;\big(x^k_0,u^k_0,b^k_0\big)=\big(x_0,\ou(0),\ob(0)\big),
\end{equation}
\begin{equation}\label{state-constraint}
\la u_{ki}^k,x_{ik}^k\ra\le b_{ki}^k\;\mbox{ for }\;i=1,\ldots,m,
\end{equation}
\begin{equation}\label{nei}
\big\|(x^k_j,u^k_j,b^k_j)-\big(\ox(t^k_j),\ou(t^k_j),\ob(t^k_j)\big)\|\le\ve/2\;\mbox{ for }\;j=0,\ldots,k,
\end{equation}
\begin{equation}\label{neighborhood1}
\sum\limits_{j=0}^{k-1}\int\limits_{{t^k_j}}^{{t^k_{j+1}}}\disp\Big(\Big\|\frac{x_{j+1}^k-x_j^k}{h_k}\disp-
\dot\ox(t)\Big\|^2+\Big\|\disp\frac{u_{j+1}^k-u_j^k}{h_k}-\dot\ou(t)\Big\|^2+\Big\|\frac{b_{j+1}^k-b_j^k}{h_k}\disp
-\dot\ob(t)\Big\|^2\Big)\,dt\le\frac{\ve}{2},
\end{equation}
\begin{equation}\label{initialubconstraint}
\;\;\;\Big\|\frac{u^k_{1}-u^k_0}{h_k}\Big\|\le\tilde{M}+1,\;\;\;\Big\|\frac{b^k_{1}-b^k_0}{h_k}\Big\|\le \tilde{M}+1,
\end{equation}
\begin{equation}\label{bvconstraint}
\sum_{j=0}^{k-2}\Big\|\frac{u_{j+2}^k-2u_{j+1}^k+u_{j}^k}{h_k}\Big\|\le \tilde{M}+1,\;\mbox{ and }\;
\sum_{j=0}^{k-2}\Big\|\frac{b_{j+2}^k-2b_{j+1}^k+b_{j}^k}{h_k}\Big\|\le \tilde{M}+1,
\end{equation}
where $\tilde{M}$ is taken from in Theorem~\ref{feasible} while $\ve>0$ is taken from Definition~\ref{ilm} with $\al=1$.
Note that the index $j$ plays the role of the discrete time in $(P_k^{\tau})$ and that inclusions \eqref{discrete-inclusion}
correspond to those in \eqref{discrete-inclusion1} at the mesh points of $\Delta_k$. Observe that, in contrast to part
\eqref{contr-disc} of the state constraints, the other part \eqref{state-constraint} needs to be imposed only at the endpoints
$(x^k_k,u^k_k,b^k_k)$ while the counterparts of \eqref{state-constraint} at $(x^k_j,u^k_j,b^k_j)$ for $j=0,\ldots,k-1$ follows
from \eqref{discrete-inclusion} due to the structure of $F$ in \eqref{F} by the definition of the normal cone in \eqref{nor}.
It is important to emphasize that the set of {\em feasible solutions} to each problem $(P^\tau_k)$ with $\tau\in[0,T]$
and $k\in\N$ sufficiently large is {\em nonempty} by Theorem~\ref{feasible}.\vspace*{0.03in}

To employ and justify the method of discrete approximations in deriving necessary optimality conditions for the control
sweeping process, we need to make sure that for all $\tau\in[0,T]$ and all $k\in\N$ sufficiently large each problem $(P_k^{\tau})$ admits an {\em optimal solution}. Despite the finite-dimensionality, this issue is nontrivial for $(P^\tau_k)$ due to the possible {\em nonclosedness} of the feasible solution set to this problem because of the dynamic constraints (\ref{discrete-inclusion}) generated by the normal cone to the moving set in (\ref{F}); see \cite[Example~4.5]{chhm2}. To overcome such a possibility, we employed in \cite{chhm2} the {\em Positive Linear Independence Constraint Qualification} (PLICQ) for the given i.l.m.\ $\oz(\cdot)$ in the original problem $(P)$ formulated as follows:
\begin{eqnarray}\label{PLICQ}
\Big[\disp\sum_{i\in I(\ox(t),\ou(t),\ob(t))}\al_i\ou_i(t)=0,\;\al_i\ge 0\Big]\Longrightarrow\al_i=0\;
\mbox{ as }\;i\in I\big(\ox(t),\ou(t),\ob(t)\big)\;\mbox{ on } \;[0,T],
\end{eqnarray}
where the collection of the active constraint indices $i\in I(\ox(t),\ou(t),\ob(t))$ is defined by
\begin{equation}\label{active_constr}
I(x,u,b):=\big\{i\in\{1,\ldots,m\}\big|\;\langle u_i,x\rangle=b_i\big\}.
\end{equation}
This condition, which permits the linearly dependence of active gradients, is obviously weaker than the classical
{\em LICQ} for $\oz(\cdot)$ on $[0,T]$ that has already be used in Theorem~\ref{ex-polyh}. It is worth mentioning that in our polyhedral setting
\eqref{sw-con2} under the additional assumption of $C(t)\ne\emp$ it follows that (\ref{PLICQ}) corresponds to Slater's condition, while we keep the term PLICQ here since a more general framework of moving sets described by smooth inequalities will be considered in Section \ref{sect:basics} and further developments.

In \cite[Theorem~4.4]{chhm2} we proved the existence of optimal solutions to problem $(P_k)$ with $P_k:=P^0_k$
for all large numbers $k\in\N$ under the validity of PLICQ \eqref{PLICQ} by using the normal cone/subdifferential
structure of the constraints in \eqref{discrete-inclusion} and employing Attouch's theorem on the subdifferential
convergence (see, e.g., \cite[Theorem~12.35]{a}) as well as the extremal principle of variational analysis from
\cite[Theorem~2.8]{m-book1}. The given proof holds without any change in the case of problems $(P^\tau_k)$ with
$\tau>0$ by taking into account the existence of feasible solutions to $(P^\tau_k)$ justified in Theorem~\ref{feasible}.
This brings us to the following result.

\begin{Proposition}{\bf (existence of discrete optimal solutions).}\label{ex-disc} Let the cost functions $\ph$ and
$\ell(t,\cdot,\cdot)$  be lower semicontinuous around the given i.l.m.\ $\oz(\cdot)$ satisfying PLICQ \eqref{PLICQ}
whenever $t\in[0,T]$. Then for all $\tau\in[0,T]$ and all $k\in\N$ sufficiently large there exist optimal solutions
to the discrete problems $(P_k^{\tau})$.
\end{Proposition}

To proceed next with establishing the concluding result of this section on the desired strong convergence of
optimal solutions for $(P^\tau_k)$ to the given local minimizer $\oz(\cdot)$, we need to impose one more requirement
on $\oz(\cdot)$. Fix any quadruple $(t,x,u,b)$ and denote $\Hat{\ell}_F(t,x,u,b,v,w,\nu)$ the {\em convexification}
of the integrand in \eqref{eq:MP} on the set $F(x,u,b)$ from \eqref{F} with respect to the {\em velocity} variables $(v,w,\nu)$,
i.e., the largest convex and lower semicontinuous (l.s.c.) function majorized by ${\ell}(t,x,u,b,\cdot,\cdot,\cdot)$ on this set.
Then for all $\tau\ge$ define the {\em relaxation} $(R^\tau)$ of problem $(P^\tau)$ as follows:
\begin{equation}\label{R}
{\rm{minimize }}\;\;\Hat J[z]:=\varphi\big(x(T)\big)+\int_0^T\Hat{\ell}_{F}\big(t,x(t),u(t),b(t),\dot{x}(t),\dot{u}(t),\dot{b}(t)\big)\,dt
\end{equation}
over the triples $z(t)=(x(t),u(t),b(t))$ of absolutely continuous functions on $[0,T]$ satisfying the constraints in
\eqref{sw-con1}--\eqref{sw-con3} and \eqref{delta}. It follows from the construction of $\Hat{\ell}_F$ and the convexity
of the set on the right-hand side of \eqref{sw-con1} that the relaxed problem $(R^\tau)$ reduces to the original one
$(P^\tau)$ if the integrand $\ell$ in \eqref{eq:MP} is convex and l.s.c.\ with respect to the velocity variables $(v,w,\nu)$.
In the general case we say that $\oz(\cdot)$ is a {\em relaxed intermediate local minimizer} (r.i.l.m.) for $(P^{\tau})$ if
it is an i.l.m.\ for this problem with $J[\oz]=\Hat J[\oz]$.

A remarkable phenomenon well-recognized in control theory for continuous-time systems reveals that in many nonconvex
settings the value of the cost functional does not change under the integrand convexification with respect to velocity variables.
It is known as ``hidden convexity" being related to Bogolyubov-type relaxation results and Lyapunov's convexity theorem for
integrals of set-valued mappings; see, e.g., \cite{AC,m-book1,v}. To the best of our knowledge, the most general
Bogolyubov-type theorem is obtained in \cite{dfm07} for optimal control problems governed by differential inclusions
satisfying the so-called ``modified one-sided Lipschitzian" condition with respect to state variables. However, the latter condition does not hold for the sweeping inclusion \eqref{sw-con1}. Thus we cannot drop so far the relaxation property of intermediate local minimizers in the following theorem, which is crucial for justifying the method of discrete approximations to derive necessary optimality conditions for the
sweeping control problems under consideration.

\begin{Theorem}{\bf (strong convergence of discrete optimal solutions).}\label{conver}
Given an arbitrary number $\tau\in[0,T]$, let $\oz(\cdot)=(\ox(\cdot),\ou(\cdot),\ob(\cdot))$ be a r.i.l.m.\ for problem $(P^{\tau})$
satisfying the assumptions of Theorem~{\rm\ref{feasible}} and Proposition~{\rm\ref{ex-disc}} and such that $\ph$ is
continuous around $\ox(T)$, $\ell(t,\cdot,\cdot)$ is continuous around $(\oz(t),\dot\oz(t))$ uniformly on $[0,T]$ while
$\ell(\cdot,z,\dot z)$ is a.e.\ continuous on $[0,T]$ being uniformly majorized by a summable function near the given local minimizer.
Then any sequence of piecewise linear extensions of the optimal solutions $\oz^k=(\ox^k,\ou^k,\ob^k)$ to the discrete
problems $(P_k^{\tau})$ converges to $\oz(\cdot)$ in the norm topology of $W^{1,2}[0,T]$.
Furthermore, we have the following estimates:
\begin{equation}\label{initialubconstraint1}
\limsup_{k\to\infty}\Big\|\frac{\ou^k_{1}-\ou^k_0}{h_k}\Big\|\le \tilde{M},\quad\limsup_{k\to\infty}\Big\|\frac{\ob^k_{1}-\ob^k_0}{h_k}\Big\|\le\tilde{M}
\;\mbox{ for all }\;k\in\N,
\end{equation}
\begin{equation}\label{bvbounded}
\limsup_{k\to\infty}\sum_{j=0}^{k-2}\Big\|\frac{\ou_{j+2}^k-2\ou_{j+1}^k+\ou_{j}^k}{h_k}\Big\|\le \tilde{M},\;\mbox{ and }\;
\limsup_{k\to\infty}\sum_{j=0}^{k-2}\Big\|\frac{\ob_{j+2}^k-2\ob_{j+1}^k+\ob_{j}^k}{h_k}\Big\|\le \tilde{M},
\end{equation}
where the constant $\tilde{M}\ge M$ is taken from \eqref{uniformbounded}.
\end{Theorem}
{\bf Proof.} Fix any sequence of the (well-defined by Proposition~\ref{ex-disc}) extended optimal solutions $\oz^k(\cdot)$
to problems $(P^\tau_k)$  and observe that strong $W^{1,2}$-convergence to $\oz(\cdot)$ on $[0,T]$ as well as the properties
\eqref{initialubconstraint1} and \eqref{bvbounded} follow directly from the equality
\begin{align}
\mathop{\lim}\limits_{k\to\infty}\Bigg(\int_0^T\Big(\Big\|\dot\ox^k(t)\disp-\dot\ox(t)\Big\|^2+\Big\|\dot\ou^k(t)\disp-\dot\ou(t)\Big\|^2
+\Big\|\dot\ob^k(t)\disp-\dot\ob(t)\Big\|^2\Big)\,dt\qquad\qquad\qquad\qquad\qquad&\notag\\
+\,{\rm dist}^2\Big(\Big\|\frac{u_{1}^k-u_{0}^k}{h_k}\Big\|,(-\infty,\tilde{M}]\Big)+{\rm dist}^2
        \Big(\sum_{j=0}^{k-2}\Big\|\frac{b_{1}^k-b_{0}^k}{h_k}\Big\|,(-\infty,\tilde{M}]\Big)\qquad\label{convergence}\\
+\,{\rm dist}^2\Big(\sum_{j=0}^{k-2}\Big\|\frac{\ou_{j+2}^k-2\ou_{j+1}^k+\ou_{j}^k}{h_k}\Big\|,(-\infty,\tilde{M}]\Big)
+{\rm dist}^2\Big(\sum_{j=0}^{k-2}\Big\|\frac{\ob_{j+2}^k-2\ob_{j+1}^k+\ob_{j}^k}{h_k}\Big\|,(-\infty,\tilde{M}]\Big)\Bigg)
&=0\notag
\end{align}
due to the initial conditions $\big(\ox^k(0),\ou^k(0),\ob^k(0)\big)=\big(\ox(0),\ou(0),\ob(0)\big)$ as $k\in\N$.
To justify \eqref{convergence}, suppose the contrary, i.e., the limit along a subsequence therein (no relabeling)
equals to some $\gg>0$. By the weak compactness of the unit ball in $L^2[0,T]:=L^2([0,T];\R^n\times\R^{nm}\times\R^m)$,
find $(v(\cdot),w(\cdot),\nu(\cdot))\in L^2[0,T]$ and (if necessary) another subsequence of $\{\oz^k(\cdot)\}$ so that
\[
\big(\dot{\ox}^k(\cdot),\dot{\ou}^k(\cdot),\dot{\ob}^k(\cdot)\big)\to\big(v(\cdot),w(\cdot),\nu(\cdot)\big)\;\mbox{ weakly in }\;L^2[0,T].
\]
Next we define the absolutely continuous triple $\Tilde z(\cdot):=(\Tilde x(\cdot),\Tilde u(\cdot),\Tilde b(\cdot))\colon[0,T]\to\R^{n+nm+m}$ by
\[
\Tilde z(t):=\big(x_0,\ou(0),\ob(0)\big)+\int\limits_0^t\big(v(s),w(s),\nu(s)\big)\,ds,\quad t\in[0,T],
\]
for which $\dot{\Tilde z}(t)=(v(t),w(t),\nu(t))$ a.e.\ on $[0,T]$. Applying Mazur's weak closure theorem gives us a sequence of convex combinations of
$(\dot{\ox}^k(\cdot),\dot{\ou}^k(\cdot),\dot{\ob}^k(\cdot))$ converging to $(v(\cdot),w(\cdot),\nu(\cdot))$ strongly in $L^2[0,T]$ and thus a.e.\ on $[0,T]$ along a subsequence. Then passing to the limit as $k\to\infty$ in the discrete inclusions \eqref{discrete-inclusion}
with the convex normal structure of the mapping $F$ from \eqref{F} and employing the aforementioned
Attouch's theorem tell us that $\Tilde x(\cdot)$ satisfies the sweeping inclusion \eqref{sw-con1} with
the set $C(t)$ in \eqref{sw-con2} defined via $\Tilde u(\cdot)$ and $\Tilde b(\cdot)$. The validity of the $\tau$-constraints in \eqref{delta} for $\Tilde u(\cdot)$ follows, whenever $\tau\in[0,T]$, from the uniform convergence on $[0,T]$ of the designated sequence of convex
combinations of $\ou^k(\cdot)$ to the limiting control function $\Tilde u(\cdot)$. It also follows from the strong
$L^2$-convergence of the above convex combinations of $(\dot{\ox}^k(\cdot),\dot{\ou}^k(\cdot),\dot{\ob}^k(\cdot))$
that the limiting triple $\Tilde z(\cdot)$ belongs to the prescribed $\ve$-neighborhood (in $W^{1,2}$) of the i.l.m.\ $\oz(\cdot)$ from Definition~\ref{ilm}.

It remains to pass to the limit in the discrete cost functional \eqref{cost-pk} along the optimal triple
$\oz^k(\cdot)$ for $(P^\tau_k)$ as $k\to\infty$. We can directly deduce from the construction of $\Hat{\ell}_F$ and its convexity in velocities that
\[
\int_0^T\Hat{\ell}_{F}\big(t,\Tilde x(t),\Tilde u(t),\Tilde b(t),\dot{\Tilde x}(t),\dot{\Tilde u}(t),\dot{\Tilde b}(t)\big)\,dt
\le \liminf_{k\to\infty}\disp{h_k}\sum\limits_{j=0}^{k-1}
\disp{\ell\Big(t^k_j,\ox_j^k,\ou_j^k,\ob_j^k,\frac{\ox_{j+1}^k-\ox_j^k}{h_k},\frac{\ou_{j+1}^k-\ou_j^k}{h_k}},\frac{\ob_{j+1}^k-\ob_j^k}{h_k}\Big).
\]
By the structure of \eqref{cost-pk}, the lower semicontinuity of the total variation, and the choice of $\gg$ above we get
\begin{equation}\label{ine1}
\disp\Hat J[\Tilde z]+\gamma=\ph\big(\Tilde{x}(T)\big)+\int_0^T\Hat
\ell_F\big(t,\Tilde{x}(t),\Tilde u(t),\Tilde b(t),\dot{\Tilde x}(t),\dot{\Tilde u}(t),\dot{\Tilde b}(t)\big)\,dt+\gamma\le\liminf_{k\to\infty}J_k[\oz^k]
\end{equation}
by using the Lebesgue dominated convergence theorem due to the assumptions made.
On the other hand, applying Theorem~\ref{feasible} to the local minimizer $\oz(\cdot)$
under consideration gives us a sequence $\{z^k(\cdot)\}$ of the feasible solutions to $(P_k^{\tau})$
that approximates $\oz(\cdot)$ in the norm topology of $W^{1,2}[0,T]$. Since $\oz_k(\cdot)$ is an
optimal solution to problem $(P_k^{\tau})$ while $z^k(\cdot)$ is feasible to it for each $k$, we have
\begin{equation}\label{appr}
J_k[\oz^k]\le J_k[z^k]\;\mbox{ whenever }\;k\in\N.
\end{equation}
It now follows from the structure of the cost functional \eqref{cost-pk} in $(P_k^{\tau})$ with $\tilde{M}\ge M$, the strong $W^{1,2}$-convergence in
Theorem~\ref{feasible}, and the assumed continuity of $\ph$ and $\ell$ that $J_k[z^k]\to J[\oz]$ as $k\to\infty$. 
Thus by taking (\ref{appr}) into account we obtain
\begin{equation}\label{ine2}
\limsup_{k\to\infty}J_k[\oz^k]\le J[\oz].
\end{equation}
The obtained relationships (\ref{ine1}) and (\ref{ine2}) together with the assumption on $\gg>0$
imply that $\widehat J[\Tilde{z}]<\widehat J[\oz]$ contradicting therefore the choice of $\oz(\cdot)$
as a r.i.l.m.\ for $(P^\tau)$. Hence $\gg=0$, which shows that \eqref{convergence} holds and thus completes the proof of the theorem. $\h$

\section{Generalized Differentiation and Second-Order Calculations}\label{sect:basics}
\setcounter{equation}{0}

After establishing well-posedness of the discrete approximation problems $(P^\tau_k)$ and the desired
strong convergence of their optimal solutions to the given r.i.l.m.\ $\oz(\cdot)$ for the sweeping control
problem $(P^\tau)$ with any fixed $\t\in[0,T]$, our further strategy is as follows: obtain necessary
optimality conditions for finite-dimensional discrete-time problems $(P^\tau_k)$ whenever $k\in\N$ and
then justify the possibility of passing to the limit as $k\to\infty$ in the obtained discrete
relationships as to derive necessary optimality conditions for $\oz(\cdot)$ in $(P^\t)$. Since problems $(P^\tau_k)$ and $(P^\tau)$ are always {\em nonsmooth} due to the dynamic constraints
independently on the smoothness of the cost functions $\ph$ and $\ell$ in \eqref{eq:MP}, we have to
employ appropriate generalized differential constructions of variational analysis enjoying comprehensive
calculus and robustness properties. In our setting not only {\em first-order} but also {\em second-order} generalized
differentiation is needed.

The main results of this section give upper estimates as well as precise formulas for calculating the coderivative of the normal cone mapping to moving
sets as in \eqref{F}, which is a second-order object playing a decisive role in the subsequent results of this paper.
We begin with some basic definitions from generalized differentiation while referring the reader to \cite{m-book,RW}
for more details on the first-order constructions and to \cite{m-book} and the papers mentioned below for
the second-order ones and their equivalent descriptions.

Recall that, for a set-valued mapping/multifunction $F\colon\R^n\tto\R^m$, the symbol
\begin{equation}\label{pk}
\disp\Limsup_{x\to\ox}F(x):=\big\{y\in\R^m\big|\;\exists\;\mbox{ sequences }
\;x_k\to\ox,\;y_k\to y\;\mbox{ with }\;y_k\in F(x_k)\;\mbox{ for all }\;k\in\N\big\}
\end{equation}
signifies the (Kuratowski-Painlev\'e) {\em outer limit} of $F$ at $\ox$.
Given a subset $\O\subset\R^n$ locally closed around $\ox\in\O$, the {\em normal cone} to
$\O$ at $\ox$ (known also as the limiting/basic/Mordukhovich one) is defined by
\begin{equation}\label{nc}
N(\ox;\O)=N_\O(\ox):=\Limsup_{x\to\ox}\big\{\mbox{cone}\big[x-\Pi(x;\O)\big]\big\}
\end{equation}
via the outer limit \eqref{pk}, where $\Pi(x;\O)$ stands for the Euclidean projection of $x$ onto $\O$,
and where `cone' denotes the conic hull of the set. When $\O$ is convex, the normal cone
\eqref{nc} reduces to the classical one of convex analysis, while in general the cone \eqref{nc}
is nonconvex even for simple sets $\O$, e.g., for $\O:=\{(x_1,x_2)\in\R^2|\;x_2=|x_1|\}$.
Nevertheless, the normal cone and associated subdifferential and coderivative constructions
for functions and multifunctions enjoy {\em full calculi} based on variational principles; see \cite{m-book,RW}.

Given a set-valued mapping $F\colon\R^n\times\R^m$ whose graph
$$
\gph F:=\big\{(x,y)\in\R^n\times\R^m\big|\;y\in F(x)\big\}
$$
is locally closed around $(\ox,\oy)$, the {\em coderivative} of $F$ at $(\ox,\oy)$ is defined by
\begin{equation}\label{cod}
D^*F(\ox,\oy)(u):=\big\{v\in\R^n\big|\;(v,-u)\in N\big((\ox,\oy);\gph F\big)\big\},\quad u\in\R^m,
\end{equation}
where $\oy=F(\ox)$ is omitted if $F$ is single-valued. When $F\colon\R^n\to\R^m$ is smooth around $\ox$, we have
$$
D^*F(\ox)(u)=\big\{\nabla F(\ox)^*u\big\}\;\mbox{ for all }\;u\in\R^n,
$$
with $A^*$ standing for adjoint operator/matrix transposition of the Jacobian $A=\nabla F(\ox)$.

For a l.s.c.\ extended-real-valued function $\ph\colon\R^n\to\oR$ with the domain and epigraph
$$
\dom\ph:=\big\{x\in\R^n\big|\;\ph(x)<\infty\big\}\;\mbox{ and }\;\epi\ph:=\big\{(x,\mu)\in\R^{n+1}\big|\;\mu\ge\ph(x)\big\}
$$
its (first-order) {\em subdifferential} at $\ox\in\dom\ph$ is generated by \eqref{nc} as
\begin{equation}\label{subdifferential}
\partial\varphi(\ox):=\big\{v\in\R^m\big|\;(v,-1)\in N\big((\ox,\ph(\ox);\epi\ph\big)\big\}.
\end{equation}

Our main objects here are the {\em second-order} generalized differential constructions defined by the scheme of \cite{m-book} as follows.
Given $\ov\in\partial\ph(\ox)$ from \eqref{subdifferential}, the {\em second-order subdifferential}
(or {\em generalized Hessian}) of $\ph$ at $\ox$ relative to $\ov$ is the mapping $\partial^2\ph(\ox,\ov)\colon\R^n\tto\R^n$ with the values
\begin{equation}\label{second-order-sub}
\partial^2\varphi(\ox,\ov)(u):=(D^*\partial\varphi)(\ox,\ov)(u),\quad u\in\R^n.
\end{equation}
Having an extended real-valued function $\ph\colon\R^n\times\R^d\to\oR$ of two variables $(x,w)\in\R^n\times\R^d$
and its partial (in $x$) first-order subdifferential mapping
\begin{equation*}
\partial_x\ph(x,w):=\big\{{\rm{set\;of\;subgradients}}\;v\;{\rm{of}}\;\ph_w:=\ph(\cdot,w)\;{\rm{at}}\;x\big\},\quad (x,w)\in\dom\ph,
\end{equation*}
define the {\em partial second-order subdifferential} of $\ph$ in $x$ at $(\ox,\ow)$ relative to $\ov\in\partial_x\ph(\ox,\ow)$ by
\begin{equation}\label{par2}
\partial^2_x\ph(\ox,\ow,\ov)(u):=(D^*\partial_x\ph)(\ox,\ow,\ov)(u),\quad u\in\R^n.
\end{equation}
Note that for ${\cal C}^2$-smooth functions $\ph$ the constructions in \eqref{second-order-sub} and \eqref{par2} reduce,
respectively, to
\begin{equation*}
\partial^2\ph(\ox)(u)=\big\{\nabla^2_{xx}\ph(\ox)u\big\},\quad\partial^2\ph(\ox,\ow)(u)=\big\{\big(\nabla^2_{xx}\ph(\ox,\ow),
\nabla^2_{xw}\ph(\ox,\ow)\big)\big\},\quad u\in\R^n,
\end{equation*}
expressed in terms of the classical (symmetric) Hessian matrices.
The partial second-order construction \eqref{par2} has been studied in \cite{BR}
under the name of ``extended partial second-order subdifferential" with the notation
$\Tilde\partial^2_x\ph$. Since no other partial second-order constructions are used here,
we drop both the word ``extended" and the tilde-notation for \eqref{par2}. Our goal is to
{\em estimate} and {\em calculate} this second-order construction for the special class of
functions arising in the controlled sweeping process \eqref{sw-con1}.

To proceed further, consider the smooth parametric inequality system
\begin{equation}\label{ctild}
S(w):=\{x\in\R^n\big|\;g(x,w)\in\R^m_-\big\},
\end{equation}
where $\R^m_-$ is the nonpositive orthant of $\R^m$ and  $g\colon\R^n\times\R^d\to\R^m$ is an arbitrary
${\cal C}^2$-smooth vector function. Associate with \eqref{ctild} the {\em normal cone mapping} ${\cal N}\colon\R^n\times\R^d\tto\R^n$ defined by
\begin{equation}\label{nm}
{\cal N}(x,w):=N\big(x;S(w)\big)\;\mbox{ for }\;x\in S(w)
\end{equation}
via the normal cone \eqref{nc} to the moving set $S(w)$, which we denote as $N_{S(w)}(x)$ for convenience.
It is easy to see that the mapping ${\cal N}$ in \eqref{nm} admits the {\em composite representation}
\begin{equation*}
{\cal N}(x,w)=\partial_x\ph(x,w)\;\mbox{ with }\;\ph(x,w):=\big(\dd_{\R^m_-}\circ g\big)(x,w)
\end{equation*}
by using the indicator function $\dd_{\R^m_-}$ of the orthant $\R^m_-$. Thus we get by definition \eqref{par2} that
\begin{equation}\label{nm1}
\partial^2_x\ph(\ox,\ow,\ov)(u)=D^*{\cal N}(\ox,\ow,\ov)(u)\;\mbox{ for any }\;\ov\in{\cal N}(\ox,\ow)\;\mbox{ and }\;u\in\R^n.
\end{equation}
Dealing with the moving set \eqref{ctild}, we use in what follows the coderivative form
\eqref{nm1} of the partial second-order subdifferential of the function $\ph$ in question.\vspace*{0.03in}

The main issue is to evaluate this construction entirely via the given data of \eqref{ctild}.
The next lemma based on the {\em second-order chain rules} from \cite{MO07,BR} plays an important
role in the subsequent calculations. Recall that the mapping $M\colon\R^s\tto\R^q$ is {\em calm} at
$(\bar s,\bar q)\in\gph M$ if there are numbers $\mu\ge 0$ and $\eta>0$ with
\begin{equation}\label{calm1}
M(s)\cap(\bar q+\eta\B)\subset M(\bar s)+\mu\|s-\bar s\|\B\;\mbox{ whenever }\;s\in\bar s+\eta\B.
\end{equation}
\begin{Lemma}{\bf (coderivative of the normal cone mapping for smooth inequality systems).}\label{morout}
Let $(\ox,\ow)\in\R^n\times\R^d$ be such that $g(\ox,\ow)\in\R^m_-$, let $\ov\in{\cal N}(\ox,\ow)$, and let
$$
I(\ox,\ow):=\big\{i\in\{1,\ldots,m\}\big|\;g_i(\ox,\ow)=0\big\}
$$
be the collection of active indices for \eqref{ctild} at $(\ox,\ow)$. The following assertions hold:

{\bf (i)} Assume that the partial gradients $\{\nabla_x g_i(\ox,\ow)|\;i\in I(\ox,\ow)\}$ are positively linearly
independent and that the mapping $\vartheta\mapsto\{(x,w,p)|\;(g(x,w),p)+\vartheta\in\gph N_{\R^m_-}\}$ is calm
at $(0,\ox,\ow,p)$ for all $p\in N_{\R^m_-}(g(\ox,\ow))$ with $\nabla_x g(\ox,\ow)^*p=\ov$. Then for all $u\in\R^n$ we have the upper estimate
\begin{eqnarray*}
&D^*{\cal N}(\bar{x},\bar{w},\bar{v})(u)\subset&\\&
\bigcup\limits_{\begin{subarray}{l}p\in N_{\mathbb{R}_{-}^{m}}(g(\bar{x},\bar{w}))\\\nabla_{x}g(\bar{x},\bar{w})^*p=\bar{v}\end{subarray}}
\left\{\left[\begin{array}{c}
\nabla_{xx}^{2}\langle p,g\rangle(\bar{x},\bar{w})\\
\nabla_{xw}^{2}\langle p,g\rangle(\bar{x},\bar{w})
\end{array}\right]u+\nabla g(\bar{x},\bar{w})^*D^*N_{\mathbb{R}_{-}^{m}}\big(g(\bar{x},\bar{w}),p\big)
\big(\nabla_{x}g(\bar{x},\bar{w})u\big)\right\}.&\nonumber
\end{eqnarray*}

{\bf (ii)} Assume that the partial gradients $\{\nabla_x g_i(\ox,\ow)|\;i\in I(\ox,\ow)\}$ are linearly independent,
and let the vector $\bar p\in\R^m$ be uniquely defined by
\begin{equation*}
\bar p\in N_{\mathbb{R}_{-}^{m}}\big(g(\bar{x},\bar{w})\big),\quad\nabla_{x}g(\bar{x},\bar{w})^*\bar p=\bar{v}.
\end{equation*}
Then for all $u\in\R^n$ we have the precise coderivative formula
\begin{eqnarray*}
D^*{\cal N}(\bar{x},\bar{w},\bar{v})(u)=\left[\begin{array}{c}
\nabla_{xx}^{2}\langle\bar p,g\rangle(\bar{x},\bar{w})\\
\nabla_{xw}^{2}\langle\bar p,g\rangle(\bar{x},\bar{w})
\end{array}\right]u+\nabla g(\bar{x},\bar{w})^*D^*N_{\mathbb{R}_{-}^{m}}\big(g(\bar{x},\bar{w}),\bar p\big)\big(\nabla_{x}g(\bar{x},\bar{w})u\big).
\end{eqnarray*}
\end{Lemma}
{\bf Proof.} We derive (i) from the second-order chain rule of the inclusion type established in \cite[Corollary~3.2(b)]{MO07},
where the first-order qualification condition follows from the positive linear independence of the active constraint gradients in \eqref{ctild}.
Assertion (ii) is a direct consequence of the precise (equality type) second-order chain rule obtained
in \cite[Theorem~3.1]{BR} under the full rank condition, which is ensured here by the assumed linear independence
of the active constraint gradients. $\h$\vspace*{0.05in}

Next we apply these results to the case of {\em bilinear} vector function $g(x,w)$ in \eqref{ctild}, which covers
our controlled sweeping setting in \eqref{sw-con1}. Define
\begin{equation}\label{aff}
g(x,w):=Ax-b\;\mbox{ with }\;w:=(A,b)\;\mbox{ for }\;x\in\R^n\;\mbox{ and }\;b\in\R^m,
\end{equation}
where $A$ is an $m\times n$-matrix, and {\em both} $A$ and $b$ are variable. In this case system \eqref{ctild} is written as
\begin{equation}\label{ctild2}
S(A,b):=\big\{x\in\mathbb{R}^{n}\big|\;Ax\le b\big\}.
\end{equation}
Taking into account that the values of $S(\cdot,\cdot)$ are polyhedral sets, we refer to \eqref{ctild2} as to the
{\em polyhedral system}. Note that the graph of $S$ {\em may not} be a convex polyhedron in $\R^n\times\R^{nm}\times\R^m$.
For any fixed $(\bar A,\ob)$ the active index set from Lemma~\ref{morout} reduces to
$$
I(\ox,\bar A,\ob):=\big\{i\in\{1,\ldots,m\}\big|\;\bar A_i\ox=\ob_i\big\},
$$
and we label $\{\bar A_i|\;i\in I(\ox,\bar A,\ob)\}$ as {\em active rows}.
Based on Lemma~\ref{morout} and the affine structure of \eqref{aff}, we arrive at the next lemma,
which relates the coderivative of ${\cal N}$ with that of $N_{\mathbb{R}_{-}^{m}}$.

\begin{Lemma}{\bf (coderivative of the normal cone mapping for polyhedral systems).}\label{coderest}
Let $(\ox,\bar A,\bar b)\in\R^n\times\R^{nm}\times\R^m$ be such that $\bar A\ox\le\ob$, and let
$\ov\in{\cal N}(\ox,\bar A,\ob)$ for the corresponding normal cone mapping \eqref{nm} generated by the polyhedral
system \eqref{ctild2}. Assume that the active rows $\{\bar A_i|\;i\in I(\ox,\bar A,\ob)\}$ are positively linearly independent.
Then we have the upper estimate
\begin{eqnarray*}
&&D^*{\cal N}(\ox,\bar{A},\bar{b},\bar{v})(u)\subset\\
&&\bigcup\left\{\left.\left(\begin{array}{c}
\bar{A}^*q\\\hline
p_{1}u+q_{1}\bar{x}\\
\vdots\\p_{m}u+q_{m}\bar{x}\\\hline
-q
\end{array}
\right)\right\vert{p\in N_{\R^m_-}(\bar A\ox-\ob),\;\bar A^*p=\ov,\;q\in D^*N_{\mathbb{R}_{-}^{m}}(\bar{A}
\bar{x}-\bar{b},p)(\bar{A}u)}\right\},\quad u\in\R^n,
\end{eqnarray*}
for all $u\in\R^n$. If moreover the active rows $\{\bar A_i|\;i\in I(\ox,\bar A,\ob)\}$ are linearly independent,
then we have the precise formula for the coderivative calculation
\begin{equation*}
D^*{\cal N}(\ox,\bar{A},\bar{b},\bar{v})(u)=\bigcup\left\{\left.\left(\begin{array}{c}
\bar{A}^*q\\\hline
\bar p_{1}u+q_{1}\bar{x}\\
\vdots\\\bar p_{m}u+q_{m}\bar{x}\\\hline
-q
\end{array}
\right)\right\vert{q\in D^*N_{\mathbb{R}_{-}^{m}}(\bar{A}
\bar{x}-\bar{b},\bar p)(\bar{A}u)}\right\},\quad u\in\R^n,
\end{equation*}
where the vector $\bar p\in N_{{\mathbb{R}_{-}^{m}}(\bar{A}\bar{x}-\bar{b})}$ is uniquely determined by ${\bar{A}}^*{\bar p=\bar{v}}$.
\end{Lemma}
{\bf Proof.} Applying Lemma~\ref{morout} to $g(x,w)$ from \eqref{aff} and using $\ox^T$ for the corresponding vector row yield
\[
\nabla g\left(\ox,\bar{A},\bar{b}\right)=\left(\begin{array}{c|c|c}
\bar{A}&
\begin{array}{ccc}
\bar{x}^T&0&0\\
0&\ddots&0\\
0&0&\bar{x}^T
\end{array}
&-I
\end{array}
\right),
\]
\begin{eqnarray*}
\nabla_{xx}^{2}\langle p,g\rangle=0,\quad\nabla_{x,(A,b)}^{2}\langle p,g\rangle=\big(p_1I\;
\big|\;\ldots\;|\;p_mI\;|\;0\big)^*
\end{eqnarray*}
for any fixed $(\ox,\ow)$ with $\ow=(\bar A,\ob)$ and any $p\in\mathbb{R}^{m}$. Observe that the mapping
\[
M(\vartheta):=\big\{(x,w,p)\big|\;(Ax-b,p)+\vartheta\in\gph N_{\mathbb{R}_{-}^{m}}\big\},\quad\vartheta=(\vartheta_1,\vartheta_2)\in\R^{2m},
\]
is automatically calm at $(0,\ox,\bar A,\ob,p)$ for any $p$ as required in Lemma~\ref{morout}.
This is a consequence of the polyhedrality of $M$ by the classical Robinson theorem from \cite{rob81}.
Thus the asserted formulas follow immediately from Lemma~\ref{morout} and the Jacobian and Hessian calculations given above. $\h$\vspace*{0.03in}

Now we are ready to derive from Lemma~\ref{coderest} the desired results for evaluating the coderivative
$D^*{\cal N}$ of \eqref{nm} entirely via the given data of \eqref{ctild2} by using the calculations of
$D^*N_{\mathbb{R}_{-}^{m}}$ available in the literature. Consider the mapping ${\cal F}(x,A,b):=-{\cal N}(x,A,b)$,
which actually appears in the sweeping inclusion.

\begin{Theorem}{\bf (coderivative of the normal cone mapping via the given data).}\label{cod_pol}
In the setting of Lemma~{\rm\ref{coderest}}, suppose that the active rows $\{\bar{A}_{i}|\;i\in I(\ox,\bar A,\ob)\}$
are positively linearly independent. For all $u\in\R^n$ and $p\in\R^m$, respectively, define the sets
\begin{eqnarray*}
P(u)&:=&\{p\in N_{\R_{-}^m}(\bar{A}\bar{x}-\bar{b})
\mid\bar{A}^*p=-\bar{v}\}\;\mbox{ if }\;u\in\bigcap_{\{i|\,p_{i}>0\}}\bar{A}_{i}^{\perp}\;\mbox{ and }\;P(u):=\emp\quad\mbox{otherwise},\\
Q(p)&:=&\left\{q\in\mathbb{R}^{m}\left\vert\begin{array}{ll}
q_i=0&\mbox{if }\;\bar{A}_i\bar{x}<\bar{b}_i\;\mbox{ or if }\;\bar{A}_{i}\bar{x}=\bar{b}_{i},\;p_{i}=0,\;\bar{A}_{i}u>0
\\
q_i\ge 0&\mbox{if }\;\bar{A}_{i}\bar{x}=\bar{b}_{i},\;p_{i}=0,\;\bar{A}_{i}u<0
\end{array}\right.\right\}.
\end{eqnarray*}
Then for all $u\in\mathbb{R}^{n}$ we have the upper estimate
\begin{equation}\label{cod_inclusion}
D^*{\cal F}(\ox,\bar A,\ob,\ov)(u)\subset\bigcup\limits_{\begin{subarray}{l}p\in P(u)\\q\in Q(p)
\end{subarray}}\left\{\left(\begin{array}{c}
\bar{A}^*q\\\hline
q_1\bar{x}-p_{1}u\\
\vdots\\
q_m\bar{x}-p_mu\\\hline
-q
\end{array}
\right)\right\}.
\end{equation}
If furthermore the active rows $\{\bar{A}_i|\;i\in I(\ox,\bar A,\ob)\}$ are linearly independent, then either
\begin{equation}\label{cod_inclusion_equa}
D^*{\cal F}(\ox,\bar A,\ob,\ov)(u)=\bigcup_{q\in Q(\op)}\left\{\left(\begin{array}{c}
\bar{A}^*q\\\hline
q_1\bar{x}-\op_{1}u\\
\vdots\\
q_m\bar{x}-\op_mu\\\hline
-q
\end{array}
\right)\right\}\quad\mbox{if}\quad{u\in}\bigcap\limits_{\left\{i|\op_{i}>0\right\}}\left[\bar A_{i}\right]^{\perp},
\end{equation}
or $D^*{\cal F}(\ox,\bar A,\ob,\ov)(u)=\emp$ otherwise. Here the vector $\op\in N_{\mathbb{R}_{-}^{m}(\bar{A}\bar{x}-\bar{b})}$ is uniquely
defined by ${\bar{A}}^*\op=-\bar{v}$.
\end{Theorem}
{\bf Proof}. Observe that $D^*{\cal F}(\ox,\bar A,\ob,\ov)(u)=D^*{\cal N}(\ox,\bar{A},\bar{b},-\bar{v})(-u)$.
The claimed results follow from Lemma~\ref{coderest} by substituting therein the precise coderivative calculation
\begin{eqnarray*}
D^*N_{\R^{m}_-}(\alpha,\beta)(\gg)=\left\{\begin{array}{ll}
\emp&\mbox{ if }\;\bb_i\gg_i\ne 0\;\mbox{ for some }\;i,\\
\big\{\eta\in\R^m\big|\;\eta_i=0\mbox{ if}\;i\in I_1\;\mbox{ and }\;\eta_i\ge 0\;\mbox{for}\;i\in I_2\big\}&\mbox{otherwise}
\end{array}\right.
\end{eqnarray*}
given in \cite[p.\ 1215]{heoutsur}, where the index subsets of $\{1,\ldots,m\}$ are defined by
$$
I_1:=\big\{i\big|\;\al_i<0\}\cup\big\{i\big|\;\al_i=\bb_i=0,\;\gg_i<0\big\},\quad I_2:=\big\{i\big|\;\al_i=\bb_i=0,\;\gg_i> 0\big\}.
$$
This verifies both coderivative formulas \eqref{cod_inclusion} and \eqref{cod_inclusion_equa} of the theorem. $\h$\vspace*{0.05in}

\section{Necessary Optimality Conditions for Discrete Approximations}\label{sect:nec_cond}
\setcounter{equation}{0}

The aim of this section is to obtain necessary conditions for optimal solutions of the discrete approximation problems $(P_k^{\tau})$
for any fixed $\t\in[0,T]$ and $k\in\N$. First we derive optimality conditions for a generalized version of $(P_k^{\tau})$,
where the dynamic constraints \eqref{discrete-inclusion} are described by an arbitrary closed-graph mapping $F$.
Then, by using the coderivative calculations of Section~\ref{sect:basics}, we arrive at optimality conditions expressed
entirely via the problem data of $(P_k^{\tau})$ with $F$ given in the particular normal cone form \eqref{F} of the
sweeping process under consideration. Our standing assumptions in this and next sections are that the cost functions $\ph$
and $\ell(t,\cdot,\cdot,\cdot,\cdot,\cdot)$ are {\em locally Lipschitzian} around the points in question.
Note that the subdifferential \eqref{subdifferential} of the running cost $\ell$ is taken with respect to its all but $t$ variables.
In what follows we drop indicating the time-dependence
of $\ell$ for brevity and use the notation
\begin{eqnarray}\label{notat}
\big[c,q\big]:=\big(c_1q_1,\ldots,c_m q_m\big)\in\R^{nm}\;\mbox{ and }\;{\rm rep}_m(x):=\big(x,\ldots,x\big)\in\R^{nm}
\end{eqnarray}
for vectors $c=(c_1,\ldots,c_m)\in\R^m$, $x\in\R^n$, and $q=(q_1,\ldots,q_m)\in\R^{nm}$ with $q_i\in\R^n$ as $i=1,\ldots,m$.

\begin{Theorem}\label{discrete} {\bf (necessary optimality conditions for general discrete inclusions).}
For fixed $k\in\N$ and $\t\in[0,T]$ let $\oz^k=(x_0,
\ox^k_1\ldots,\ox_k^k,\ou_0^k,\ldots,\ou_{k}^{k},\ob_0^k,\ldots,\ob_{k}^{k})$ be an optimal solution to the
discrete problem \eqref{cost-pk}--\eqref{bvconstraint}
written in the format of $(P^\t_k)$ but with the discrete inclusion \eqref{discrete-inclusion} governed by a general closed-graph mapping $F$. For each
$j=0,\ldots,k-1$ we denote
\begin{equation}\label{th}
\Big(\theta_{j}^{xk},\theta_{j}^{uk},\theta_{j}^{bk}\Big):=2\int\limits_{{t_j^k}}^{{t_{j+1}^k}}\Big(\frac{\ox_{j+1}^k-\ox_j^k}{h_k}
\disp-\dot\ox(t),\frac{\ou_{j+1}^k-\ou_j^k}{h_k}-\dot\ou(t),\frac{\ob_{j+1}^k-\ob_j^k}{h_k}-\dot\ob(t)\Big)dt.
\end{equation}
Then there exist dual elements $\lm^k\ge 0$, $\alpha^k\in\R_+^{m}$, $\xi^k=(\xi^k_0,\ldots,\xi^k_k)\in\R^{(k+1)m}$, $p_j^k=(p^{xk}_j,p^{uk}_j,p^{bk}_j)
\in\mathbb{R}^{n+nm+m}$ as $j=0,\ldots,k$ and subgradient vectors
\begin{equation}\label{subl}
\big(w^{xk}_{j},w^{uk}_{j},w^{bk}_{j},v^{xk}_{j},v^{uk}_{j},v^{bk}_{j}\big)\in\partial\ell\left(\oz^k_j,\frac{\oz_{j+1}^k-\oz_j^k}{h_k}\right),\quad j=0,
\ldots, k-1,
\end{equation}
such that the following conditions are satisfied:
\begin{equation}\label{nontriv}
\lm^k+\disp\|\alpha^k\|+\|\xi^k\|+\sum_{j=0}^{k-1}\|p^{xk}_j\|+\|p^{uk}_0\|+\|p^{bk}_0\|\ne 0,
\end{equation}
\begin{equation}\label{complementary_d}
\alpha^k_{i}\big(\langle\bar{u}^k_{ki},\bar{x}^k_k\rangle-\bar{b}^k_{ki}\big)=0,\quad i=1,\ldots,m,
\end{equation}
\begin{equation}\label{ksi}
\xi^k_{ji}\in N\big(\|\ou^k_{ji}\|;[1/2,3/2]\big)\;\mbox{ for }\;j=0,\ldots,j_{\tau}(k)-1\;\mbox{ and }\;j=j^{\tau}(k)+1,\ldots,k,\quad i=1,\ldots,m,
\end{equation}
\begin{align}\label{transversality_end_d}
-p^{xk}_k&\in\lambda^k\partial\varphi(\ox_k^k)+\sum_{i=1}^m\alpha_i^k\bar{u}^k_{ki};\quad p^{uk}_k=
-\big[\alpha^k,{\rm rep}_m(\bar{x}^k_k)\big]-2\big[\xi_k^k,\bar{u}_{k}^k\big],
\quad p^{bk}_k=\alpha^k,
\end{align}
\begin{equation}\label{psiu}
p^{uk}_{j+1}=\lambda^k(v^{uk}_j+h_k^{-1}\theta^{uk}_j),\quad p^{bk}_{j+1}=\lambda^k(v^{bk}_j+h_k^{-1}\theta^{bk}_j),\;j=0,\ldots,k-1,
\end{equation}
\begin{equation}\label{euler1}
\begin{split}
&\left(\frac{p_{j+1}^{xk}-p_j^{xk}}{h_k}-\lambda^k w_j^{xk},\frac{p_{j+1}^{uk}-p_j^{uk}}{h_k}-
\lambda^k w_j^{uk},\frac{p_{j+1}^{bk}-p_j^{bk}}{h_k}-\lambda^k w_j^{bk},p_{j+1}^{xk}-\lambda^k\Big(v_j^{xk}+\frac{1}{h_k}\theta_j^{xk}\Big)\right)\\
&\qquad\in\left(0,\frac{2}{h_k}\big[\xi_j^k,\bar{u}_{j}^k\big],0,0\right)+
N\Big(\Big(\ox_j^k,\ou_j^k,\ob_j^k,\frac{\ox_{j+1}^k-\ox_j^k}{h_k}\Big);\gph F\Big),\quad j=0,\ldots,k-1.
\end{split}
\end{equation}
\end{Theorem}
{\bf Proof.} Throughout the proof we omit indicating the (fixed) upper index `$k$' from the statement of this theorem;
the dependence of the result on $k$  will be needed in Section~\ref{nec-con}. Let
\[
y:=(x_0,\ldots,x_k,u_0,\ldots,u_{k},b_0,\ldots,b_{k},X_0,\ldots,X_{k-1},U_0,\ldots,U_{k-1},B_0,\ldots,B_{k-1}),
\]
where $x_0$ is fixed. Take $\ve>0$ from the construction of $(P^\t_k)$ and define the mathematical program $(MP)$:
\begin{equation*}
\begin{split}
\mathrm{minimize}\quad\varphi_0[y]&:=\varphi(x_k)+{h_k}\sum\limits_{j=0}^{k-1}{\ell(x_j,u_j,b_j,X_j,U_j,B_j)}+
\sum\limits_{j=0}^{k-1}{\int\limits_{{t_j}}^{{t_{j+1}}}{\disp\Big\|(X_j,U_j,B_j)-\dot\oz(t)\Big\|^2dt}}\\
&\;\quad+\,{\rm dist}^2\Big(\Big\|\frac{u_{1}-u_{0}}{h^k}\Big\|,(-\infty,\tilde{M}]\Big)+
{\rm dist}^2 \Big(\Big\|\frac{b_{1}-b_{0}}{h^k}\Big\|,(-\infty,\tilde{M}]\Big),\\
&\;\quad+{\rm dist}^2\Big(\sum_{j=0}^{k-2}\Big\|U_{j+1}-U_{j}\Big\|,(-\infty,\tilde{M}]\Big)+
{\rm dist}^2\Big(\sum_{j=0}^{k-2}\Big\|B_{j+1}-B_{j}\Big\|,(-\infty,\tilde{M}]\Big)
\end{split}
\end{equation*}
subject to equality, inequality, and geometric constraints
\allowdisplaybreaks
\begin{align*}
f^x_j(y)&:=x_{j+1}-x_j-h_kX_j=0 \ \;{\rm for}\ \;j=0,\ldots,k-1,\\
f^u_j(y)&:=u_{j+1}-u_j-h_kU_j=0 \ \;{\rm for}\ \;j=0,\ldots,k-1,\\
f^b_j(y)&:=b_{j+1}-b_j-h_kB_j=0 \ \;{\rm for}\ \;j=0,\ldots,k-1,\\
g_i(y)&:=\langle u_{ki},x_k\rangle -b_{ki}\le 0 \ \;{\rm for}\ \;i=1,\ldots ,m,\\
d_{ji}(y)&:=\|u_{ji}\|^2-1=0 \ \;{\rm for}\ \;j=j_{\tau}(k),\ldots,j^{\tau}(k),\quad i=1,\ldots,m,\\
y \in \O_{ji}&:=\big\{y\big|\;1/2\le\|u_{ji}\|\le 3/2 \big\}\ \;{\rm for}\ \;
j=0,\ldots,j_{\tau}(k)-1\;\mbox{ and }\;j=j^{\tau}(k)+1,\ldots,k,,\quad i=1,\ldots ,m,\\
\phi_j(y)&:=\big\|(x_j,u_j,b_j)-\oz(t_j)\big\|-\ve/2\le 0\;\;\mathrm{for}\;\;j=0,\ldots,k,\\
\phi_{k+1}(y)&:=\sum\limits_{j=0}^{k-1}\int\limits_{{t^k_j}}^{{t^k_{j+1}}}\disp\Big(\Big\|(X_j,U_j,B_j)-\dot\oz(t)\Big\|^2\Big)\,dt-\frac{\ve}{2}\le 0,\\
\phi_{k+2}(y)&:=\sum_{j=0}^{k-2}\Big\|U_{j+1}-U_{j}\Big\|\le \tilde{M}+1,\\
\phi_{k+3}(y)&:=\sum_{j=0}^{k-2}\Big\|B_{j+1}-B_{j}\Big\|\le \tilde{M}+1,\\
\phi_{k+4}(y)&:=\Big\|u_{1}-u_{0}\Big\|\le(\tilde{M}+1)h^k,\\
\phi_{k+5}(y)&:=\Big\|b_{1}-b_{0}\Big\|\le(\tilde{M}+1)h^k,\\
y\in\Xi_j&:=\big\{{y\big|\;X_j\in F(x_j,u_j,b_j)}\big\}\ \;{\rm for}\ \;j=0,\ldots,k-1,\\
y\in\Xi_k&:=\big\{{y\big|\;x_0\mbox{ is fixed},(u_0,b_0)=\big(\ou(0),\ob(0)\big)}\big\}.
\end{align*}
It is easy to see that $(MP)$ and $(P_k^{\tau})$ with an arbitrary mapping $F$ are equivalent.
Thus $\bar{y}:=(\bar{z},\bar{Z})$ is an optimal solution to (MP), where $\bar{z}:=\bar{z}^k$ is
the solution of $(P_k^{\tau})$ fixed in the theorem, and where
$\bar{Z}:=(\bar{X}_0,\ldots,\bar{X}_{k-1},\bar{U}_0,\ldots,\bar{U}_{k-1},\bar{B}_0,\ldots,\bar{B}_{k-1})$.
Necessary optimality conditions for this type of mathematical programs in terms the first-order
generalized differential constructions constructions of Section~\ref{sect:basics} are well known; see, e.g.,
\cite[Theorem~5.24]{m-book1}. Furthermore, it follows from Theorem~\ref{conver} above that all
the inequality constraints in $(MP)$ associated with functions $\phi_j\,\,(j=0,\ldots,k+5)$ are
{\em inactive} for all $k$ sufficiently large, and so the corresponding multipliers do not appear in
the optimality conditions. Taking this into account, we find $\lambda\ge 0$, $\alpha\in \R^m_+$,
$\xi=(\xi_0,\ldots,\xi_k)\in\R^{m(k+1)}$, $p_j=(p_j^x,p_j^u,p_j^b)\in\R^{n+nm+m}$ as $j=1,\ldots,k$, and
\begin{equation*}
y_j^*=\big(x^*_{0j},\ldots,x^*_{kj},u^*_{0j},\ldots,u^*_{kj},b^*_{0j},\ldots,b^*_{kj}, X^*_{0j},\ldots,X^*_{(k-1)j}, U^*_{0j},\ldots,U^*_{(k-1)j},
B^*_{0j},\ldots,B^*_{(k-1)j}\big)
\end{equation*}
for $j=0,\ldots,k$, which are not all zero and satisfy the conditions \eqref{ksi} together with
\begin{equation}\label{Lagr2}
y_j^*\in N(\oy;\Xi_j)\;\mbox{ for }\;j=0,\ldots,k,
\end{equation}
\begin{equation}\label{Lagr3}
-y_0^*-\ldots-y_{k}^*\in\lambda\partial\varphi_0(\oy)+\sum_{i=1}^m\alpha_i\nabla g_i(\oy)+
\sum_{j=0}^{k}\sum_{i=1}^m\xi_{ji}\nabla d_{ji}(\oy)+\sum_{j=0}^{k-1}\big(\nabla f_j(\oy)\big)^*p_{j+1},
\end{equation}
\begin{equation}\label{complementary}
\alpha_{i}g_{i}(\oy)=0\;\mbox{ for }\;i=1,\ldots,m.
\end{equation}
It follows from the definition of $\Xi_j$ and from $f_j^x(\bar{y})=0$ that the inclusions in \eqref{Lagr2} are equivalent to
\begin{equation}\label{inclusion}
(x^*_{jj},u^*_{jj},b^*_{jj},X^*_{jj})\in N\Big(\Big(\ox_j,\ou_j,\ob_j,\frac{\ox_{j+1}-\ox_j}{h_k}\Big);\gph F\Big),\quad j=0,\ldots,k-1.
\end{equation}
Note that every other component of $y_j^*$, which does not appear in (\ref{inclusion}), is zero.
Similarly $(x^*_{0k},u^*_{0k},b^*_{0k})$ can be the only nonzero component of $y_k^*$. Therefore we have
\begin{eqnarray}\label{entfalt}
-y_0^*-\ldots-y_{k}^*&=&\big(-x_{00}^*-x_{0k}^*,-x_{11}^*,\ldots ,-x_{k-1,k-1}^*,0,-u_{00}^* -u_{0k}^*,
\ldots ,-u_{k-1,k-1}^*,0,\\&&-b_{00}^*-b_{0k}^*,\ldots ,-b_{k-1,k-1}^*,0,-X_{00}^*,\ldots,-X_{k-1,k-1}^*,0,\ldots,0\big).\notag
\end{eqnarray}
Let us now calculate the three sums on the right-hand side of \eqref{Lagr3}.
For notational convenience we just specify the nonzero components, which are indexed according
to the partition of the vector $y$ introduced at the beginning of this proof. This gives us the equalities
\begin{eqnarray*}
\left(\sum_{i=1}^m\alpha_i\nabla g_i(\oy)\right)_{(x_k,u_k,b_k)}&=&
\left(\sum_{i=1}^m\alpha_i\bar{u}_{ki},\big[\alpha,{\rm rep}_m(\bar{x}_k)\big],-\alpha\right),\\
\left(\sum_{j=0}^{k}\sum_{i=1}^m\xi_{ji}\nabla d_{ji}(\oy)\right)_{u_j}&=&2\big[\xi_j,\bar{u}_j\big],
\quad j=0,\ldots,k,\\\left(\sum_{j=0}^{k-1}\big(\nabla f_j(\oy)\big)^*p_{j+1}\right)_{(x_j,u_j,b_j)}&=&
\left\{\begin{array}{cl}
-p_1&\mbox{if }j=0\\p_j-p_{j+1}&\mbox{if }j=1,\ldots,k-1\\p_k&\mbox{if }\;j=k
\end{array}\right.,\quad j=0,\ldots,k,\\
\left(\sum_{j=0}^{k-1}\big(\nabla f_j(\oy)\big)^*p_{j+1}\right)_{(X,U,B)}&=&
-h_kp.
\end{eqnarray*}
Introducing the auxiliary Lipschitzian functions (where the first ones are actually smooth around $\oy$)
\[
\rho_j(y):=\int\limits_{t_j}^{t_{j+1}}{\disp\Big\|(X_j,U_j,B_j)-\dot\oz(t)\Big\|^2dt},\quad j=0,\ldots,k-1,
\]
\begin{align*}
\sigma(y)&:={\rm dist}^2\Big(\Big\|\frac{u_{1}-u_{0}}{h^k}\Big\|,(-\infty,\tilde{M}]\Big)+{\rm dist}^2
\Big(\Big\|\frac{b_{1}-b_{0}}{h^k}\Big\|,(-\infty,\tilde{M}]\Big),\nonumber\\
&\qquad+{\rm dist}^2 \Big(\sum_{j=0}^{k-2}\Big\|U_{j+1}-U_{j}\Big\|,(-\infty,\tilde{M}]\Big)+{\rm dist}^2
\Big(\sum_{j=0}^{k-2}\Big\|B_{j+1}-B_{j}\Big\|,(-\infty,\tilde{M}]\Big)
\end{align*}
and then employing the subdifferential sum rule from \cite[Theorem~2.33(c)]{m-book}, we arrive at the inclusion
\[
\partial\varphi_0(\oy)\subset\partial\varphi(\bar{x}^k)+h_k\sum_{j=0}^{k-1}\partial\ell(\bar{x}_j,\bar{u}_j,\bar{b}_j,\bar{X}_j,\bar{U}_j,\bar{B}_j)
+\sum_{j=0}^{k-1}\nabla\rho_j(\bar{y})+\partial\sigma(\oy).
\]
Since the function ${\rm dist}^2(x;(-\infty,\tilde{M}])$ has the null derivative at all $x\le \tilde{M}$,
it implies together with \eqref{bvbounded} that $\partial\sigma(\oy)=\{0\}$.
Furthermore, the nonzero part of $\nabla\rho_j(\bar{y})$ is given by
$\nabla_{(X_j,U_j,B_j)}\rho_j(\bar{y})=(\theta_{j}^x,\theta_{j}^u,\theta_{j}^b)$ with the triple from \eqref{th},
and thus any element of the set $\lambda\partial\varphi_0(\oy)$ can be represented in the form
\[
\lambda\big(h_kw^x,\vartheta,h_kw^u,0,h_kw^b,0,h_kv^x+\theta^x,h_kv^u+\theta^u,h_kv^b+\theta^b\big),
\]
where $\vartheta\in\partial\varphi(\bar{x}^k)$ and the components of $(w^x,w^u,w^b,v^x,v^u,v^b)$ satisfy \eqref{subl}.
Combining this with the gradient expressions above, we deduce from \eqref{inclusion} via \eqref{entfalt} the componentwise relationship:
\begin{eqnarray}
-x_{00}^*-x_{0k}^*&=&\lambda h_kw_0^x-p_1^x\label{syst1}\\
-x_{jj}^*&=&\lambda h_kw_j^x+p_j^x-p_{j+1}^x,\quad j=1,\ldots k-1,\label{syst2}\\
0&=&\lambda\vartheta+\sum_{i=1}^m\alpha_i\bar{u}_{ki}+p_k^x,\label{syst3}\\
-u_{00}^*-u_{0k}^*&=&\lambda h_kw_0^u+2\big[\xi_0,\bar{u}_0\big]-p_1^u,\label{syst4}\\
-u_{jj}^*&=&\lambda h_kw_j^u+2[\xi_j,\bar{u}_j]+p_j^u-p_{j+1}^u,j=1,\ldots k-1,\label{syst5}\\
0&=&\big[\alpha,{\rm rep}_m(\bar{x}_k)\big]+2\big[\xi_k,\bar{u}_k\big]+p_k^u,\label{syst6}\\
-b_{00}^*-b_{0k}^*&=&\lambda h_kw_0^b-p_1^b,\label{syst7}\\
-b_{jj}^*&=&\lambda h_kw_j^b+p_j^b-p_{j+1}^b,\quad j=1,\ldots k-1,\label{syst8}\\
0&=&-\alpha+p_k^b,\label{syst9}\\
-X_{jj}^*&=&\lambda(h_kv_j^x+\theta_j^x)-h_kp_{j+1}^x,\quad j=0,\ldots k-1,\label{syst10}\\
0&=&\lambda(h_kv_j^u+\theta_j^u)-h_kp_{j+1}^u,\quad j=0,\ldots k-1,\label{syst11}\\
0&=&\lambda(h_kv_j^b+\theta_j^b)-h_kp_{j+1}^b,\quad j=0,\ldots k-1.\label{syst12}
\end{eqnarray}

Now let us derive from the obtained relationships the necessary optimality conditions of the theorem with
$p_0:=(x_{0k}^*,u_{0k}^*,b_{0k}^*)$. We have already got \eqref{ksi}. Observe now that \eqref{complementary_d}
is obviously implied by \eqref{complementary}, the conditions in \eqref{transversality_end_d} follow from
\eqref{syst3}, \eqref{syst6}, and \eqref{syst9} while those in \eqref{psiu} are a consequence of
\eqref{syst11} and \eqref{syst12}. Arguing by contradiction, suppose that the nontriviality condition
\eqref{nontriv} fails. Then it follows from \eqref{syst3} that $p_k^x=0$ as well. Since $x_{0k}^*=p_0^x=0$,
we deduce from \eqref{syst1},\eqref{syst2}, and \eqref{syst10} that $x_{jj}^*=0$ and $X_{jj}^*=0$ for $j=0,\ldots,k-1$.
Furthermore, \eqref{syst11} and \eqref{syst12} yield that $p_j^u=0$ and $p_j^b=0$ for $j=1,\ldots,k$ which in turn
implies by \eqref{syst4}, \eqref{syst5}, \eqref{syst7}, and \eqref{syst8}, that also $u_{jj}^*=0$ and $b_{jj}^*=0$
for $j=0,\ldots,k-1$. As already mentioned,
the components of $y_j^*$ different from $(x^*_{jj},u^*_{jj},b^*_{jj},X^*_{jj})$ are zero for $j=0,\ldots,k-1$,
and hence $y_j^*=0$ for $j=0,\ldots,k-1$. We similarly conclude that $y_k^*=0$ due to $x^*_{0k}=p_0^x=0$.
Getting all this together contradicts the nontriviality conditions in the mathematical program $(MP)$
formulated above and thus verifies the claimed nontriviality \eqref{nontriv}.

It remains to justify of the validity of the discrete-time adjoint conditions in \eqref{euler1},
which give us a discrete-time version of the {\em extended Euler-Lagrange inclusion} \cite{m95}
for the discrete optimal control problems under consideration. To get \eqref{euler1}, we substitute
the expressions in \eqref{syst1}, \eqref{syst2}, \eqref{syst4}, \eqref{syst5}, \eqref{syst7}, \eqref{syst8},
and \eqref{syst10} into the left-hand side of \eqref{inclusion} and
deduce from it that the vector
\begin{eqnarray*}
\big(p_{j+1}^x-p_j^x-\lambda h_kw_j^x,p_{j+1}^u-p_j^u-\lambda h_kw_j^u,p_{j+1}^b-p_j^b-\lambda
h_kw_j^b,h_kp_{j+1}^x-\lambda(h_kv_j^x+\theta_j^x)\big)-\big(0,2\big[\xi_j,\bar{u}_j\big],0,0\big)
\end{eqnarray*}
for each $j=0,\ldots,k-1$ belongs to the normal cone $N\big((\ox_j,\ou_j,\ob_j,\frac{\ox_{j+1}-\ox_j}{h_k});\gph F\big)$.
Dividing the obtained inclusions by $h_k>0$, we arrive at \eqref{euler1} and thus complete the proof of the theorem.$\h$\vspace*{0.05in}

The next theorem directly addresses the discrete approximation problems $(P^\t_k)$ for the controlled sweeping process, where the mapping $F$ in
\eqref{discrete-inclusion} in given in the particular form \eqref{F}. In this main case of our interest we
are able to derive, based on the second-order calculations of Section~\ref{sect:basics}, effective
necessary optimality conditions for $(P^\t_k)$ expressed entirely via the problem data.

\begin{Theorem}{\bf (necessary optimality conditions for the discretized sweeping process).}\label{R1}
In the setting of Theorem~{\rm\ref{discrete}} consider the discretized sweeping control problem $(P^\t_k)$,
where now the mapping $F$ in \eqref{discrete-inclusion} is defined by \eqref{F}.
Assume that for the given r.i.l.m.\ $\oz(\cdot)=(\ox(\cdot,\ou(\cdot),\ob(\cdot))$ of $(P^\t)$,
which is included in $(P^\t_k)$, the PLICQ property \eqref{PLICQ} holds and that all the
components of $\ou(t)$ are not zero on $[0,T]$. Then, in addition to the dual elements
$\lambda^k\ge 0$, $\xi^k=(\xi^k_0,\ldots,\xi^k_k)\in\R^{(k+1)m}$, $p_j^k=(p^{xk}_j,p^{uk}_j,p^{bk}_j)\in\mathbb{R}^{n+nm+m}$
as $j=0,\ldots,k$ and subgradients $w^{xk}_{j},w^{uk}_{j},w^{bk}_{j},v^{xk}_{j},v^{uk}_{j},v^{bk}_{j}$
from \eqref{subl} satisfying the relationships in \eqref{ksi} and \eqref{psiu} of
Theorem~{\rm\ref{discrete}} with the enhanced nontriviality condition
\begin{equation}\label{nontriv1}
\lm^k +\|p^{uk}_0\|+\|p^{bk}_0\|\ne 0,
\end{equation}
there exist vectors $\eta_j^k\in\mathbb{R}^m_+$ as $j=0,\ldots,k-1$ and $\gamma_j^k\in\mathbb{R}^m$ as
$j=0,\ldots,k-1$ such that we have the primal and dual/adjoint dynamic relationships
\begin{eqnarray}\label{etajk}
\disp\frac{\ox_{j+1}^k-\ox_j^k}{h_k}=-\sum_{i=1}^m\eta_{ji}^k\bar{u}^k_{ji},
\end{eqnarray}
\begin{eqnarray}\label{p_iteration}
\disp\frac{p_{j+1}^{xk}-p_j^{xk}}{h_k}-\lambda^k w_j^{xk}=\sum_{i=1}^m\gamma_{ji}^k\bar{u}_{ji}^k,
\end{eqnarray}
\begin{equation}\label{ee}
\frac{p_{j+1}^{uk}-p_j^{uk}}{h_k}-\lambda^k w_j^{uk}-\disp\frac{2}{h_k}\big[\xi_j^k,\bar{u}_{j}^k\big]
=\big[\gamma^k_j,{\rm rep}_m(\ox^k_j)\big]-\Big[\eta^k_j,{\rm rep}_m\left(\lambda^k\big(v_j^{xk}+\frac{1}{h_k}\theta_j^{xk}\big)-p_{j+1}^{xk}\right)\Big],
\end{equation}
\begin{eqnarray}\label{vstar}
\gamma^k_j=\disp\lambda^kw_j^{bk}-\frac{p^{bk}_{j+1}-p^{bk}_j}{h_k}
\end{eqnarray}
for all $j=0,\ldots,k-1$ together with the implications
\begin{eqnarray}\label{etajk-def}
\langle\ou^k_{ji},\ox^k_j\ra<\ob^k_{ji}\;\Longrightarrow&\eta^k_{ji}=0\;\mbox{ when }\;j=0,\ldots,k,\quad i=1,\ldots,m,
\end{eqnarray}
\begin{eqnarray}\label{orthogonal}
\;\;\;\qquad\eta_{ji}^k>0\;\Longrightarrow&\left\langle\bar{u}^k_{ji},\disp\lambda^k\left(v_j^{xk}+\frac{1}{h_k}\theta_j^{xk}\right)-p_{j+1}^{xk}\right\rangle=0\;\mbox{ when }\; j=0,\ldots,k-1,\;\;i=1,\ldots,m,
\end{eqnarray}
\begin{equation}\label{measuresign}
\big[\langle\ou^k_{ji},\ox^k_j\rangle<\ob^k_{ji}\;\mbox{ for all }\;i=1,\ldots,m\big]\Longrightarrow\gamma^k_j=0,\quad j=0,\ldots,k-1,
\end{equation}
as well as the transversality conditions
\begin{equation}\label{pxkk}
-p^{xk}_k\in\lambda^k\partial\varphi(\ox_k^k)+\sum\limits_{i=1}^m p^{bk}_{ki}\ou_{ki}^k,
\end{equation}
\begin{equation}\label{pk-1}
p^{uk}_{k}=-\big[p^{bk}_k,{\rm rep}_m(\ox^k_k)\big]-2\big[\xi^k_k,\ou^k_k\big],
\end{equation}
\begin{equation}\label{alphak1}
p^{bk}_{ki}\ge 0,\;\mbox{ and }\;\langle\ou^k_{ki},\ox^k_k\rangle<\ob^k_{ki}\Longrightarrow p^{bk}_{ki}=0\;\;\mbox{ for all }\;i=1,\ldots,m.
\end{equation}
\end{Theorem}
{\bf Proof.} In terms of the coderivative construction  \eqref{cod} we can rewrite \eqref{euler1} in the equivalent form
\begin{equation}\label{inclu2}
\begin{split}
&\left(\frac{p_{j+1}^{xk}-p_j^{xk}}{h_k}-\lambda^k w_j^{xk},\frac{p_{j+1}^{uk}-p_j^{uk}}{h_k}-\lambda^k w_j^{uk}-
\frac{2}{h_k}[\xi_j^k,\bar{u}_{j}^k],\frac{p_{j+1}^{bk}-p_j^{bk}}{h_k}-\lambda^k w_j^{bk}\right)\\
&\qquad\in D^*F\Big(\ox_j^k,\ou_j^k,\ob_j^k,\frac{\ox_{j+1}^k-\ox_j^k}{h_k}\Big)
\Big(\lambda^k(v_j^{xk}+\frac{1}{h_k}\theta_j^{xk})-p_{j+1}^{xk}\Big),\quad j=0,\ldots,k-1.
\end{split}
\end{equation}
Using the notation of Theorem~\ref{cod_pol} and employing the coderivative upper estimate \eqref{cod_inclusion}
therein, we deduce from \eqref{inclu2} the existence of vectors $\eta_j^k\in\R^m$ and $\gamma^k_j\in Q(\eta^k_j)$
for which all the relationships in \eqref{p_iteration}, \eqref{ee}, \eqref{vstar}, \eqref{orthogonal},
and \eqref{measuresign} are satisfied together with the conditions
\begin{eqnarray}\label{firstrel}
-(\bar{A}^k_j)^*\eta_j^k=\frac{\ox_{j+1}^k-\ox_j^k}{h_k},\;\mbox{ where }\;\eta_j^k\in N_{\R_{-}^m}(\bar{A}^k_j\bar{x}_j^k-\bar{b}_j^k),\quad j=1,\ldots,k-1.
\end{eqnarray}
It obviously follows from the conditions on $\eta^k_j$ in \eqref{firstrel} that $\eta_j^k\in\mathbb{R}^m_+$ and
the implications in \eqref{etajk-def} hold for $j=0,\ldots,k-1$. Defining $\eta_k^k:=\alpha_k$ with $\alpha_k$
taken from Theorem \ref{discrete}, we deduce from \eqref{complementary_d} that $\eta_k^k\in\mathbb{R}^m_+$ and
\eqref{etajk-def} holds for $j=k$ as well. The equations in \eqref{etajk} are consequences of those in \eqref{firstrel}
due to the definition of $\bar{A}^k_j$ given right before Lemma~\ref{coderest}. The transversality conditions in
\eqref{pxkk}--\eqref{alphak1} are  direct consequences of \eqref{transversality_end_d} due to the relationships
in \eqref{complementary_d} and \eqref{ksi}.

It remains to verify the enhanced nontriviality condition \eqref{nontriv1}.
To proceed, suppose that $\lambda^k=0$, $p^{uk}_0=0$, and $p^{bk}_0=0$.
Then we get from \eqref{psiu} that $p_j^{uk}=0$ and $p_j^{bk}=0$ for all $j=0,\ldots,k$.
This ensures that $\gamma_j^k=0$ for $j=0,\ldots,k-1$ by \eqref{vstar}, $p_k^{xk}=0$ by \eqref{pxkk},
and consequently $p_j^{xk}=0$ for $j=0,\ldots,k$ by \eqref{p_iteration}.
Furthermore, it follows from \eqref{ee} that $[\xi_j^k,\bar{u}_{j}^k]=0$ for $j=0,\ldots,k-1$ and
from \eqref{pk-1} that $[\xi_k^k,\bar{u}_{k}^k]=0$. By recalling the definition of $[\cdot,\cdot]$ in
\eqref{notat}, the latter conditions readily implies that $\xi_j^k=0$ for all $j=0,\ldots,k$.
This is due to the assumption of $\ou_i(t)\ne 0$ on $[0,T]$ made in the theorem, which implies
that $\ou^k_{ij}\ne 0$ for the discrete approximation due to the uniform convergence of
$\ou^k(t)\to\ou(t)$ in Theorem~\ref{conver}. Therefore $\lambda^k=0$, $\xi^k=0$, $\alpha^k=p_k^{bk}=0$ by \eqref{transversality_end_d},
and so $p_j^{xk}=0$ for $j=0,\ldots,k-1$.
This contradicts the nontriviality condition \eqref{nontriv} in Theorem~\ref{discrete} and thus verifies the enhanced one in \eqref{nontriv1}. $\h$

\section{Optimality Conditions for the Controlled Sweeping Process}\label{nec-con}
\setcounter{equation}{0}

In this section we proceed with the passage to the limit as $k\to\infty$ in the necessary optimality conditions of Theorem~\ref{R1}
for problems $(P^\t_k)$ and deriving in this way, with the help of Theorem~\ref{conver} and the developed tools of generalized differentiation,
necessary optimality conditions in the original optimal control problem(s) for the sweeping process formulated in Section~1.
Our major case is problem $(P^\t)$ with $0<\t<T$, but we also consider the situation when $\t=0$, which is the same as $\t=T$.
The results obtained for $(P^\t)$ are {\em explicit}, i.e., they involve only the problem data and the given local minimizers $\oz(\cdot)$
while not requiring calculations of any auxiliary objects as, e.g., coderivatives.

Keeping the assumptions above, we impose here some additional ones on the problem data, which seem to be reasonable for the controlled sweeping model under
consideration and are illustrated below by examples. The next theorem uses notation \eqref{notat} together with the symbol `$\co$' for the convex hull.

One more remark is needed before the formulation of our main result. Since it is derived by passing to the limit in the optimality
conditions for the discrete problems $(P^\t_k)$, the subdifferential construction used in Theorem~\ref{R1} in the case of the
nondifferentiable running cost $\ell(t,\cdot,\cdot)$ has to be {\em robust}, i.e., outer semicontinuous with respect to
perturbations of the reference point. As well known, this important property holds for our subdifferential \eqref{subdifferential};
see, e.g., \cite[p.\ 11]{m-book}. In the general nonautonomous setting under consideration, the robustness of the
subdifferential of $\ell$ with respect to the time parameter is also required, and we {\em postulate} it in what follows.
It does not seem to be restrictive (see the discussion in \cite{m95}) and can be completely avoided by considering the
{\em extended subdifferential} of $\ell$ as in \cite[Sec.\ 6.1.5]{m-book1}.

\begin{Theorem} {\bf (nondegenerate necessary optimality conditions for the controlled sweeping process).}\label{necopt}
Let $\oz(\cdot)=(\ox(\cdot),\ou(\cdot),\ob(\cdot))$ be a r.i.l.m.\ for problem $(P^{\tau})$ as $\t\in(0,T)$ under
the assumptions of Theorem~{\rm\ref{conver}}, and let LICQ hold on $[0,T]$. Suppose in addition that $\ph$ is locally Lipschitzian around $\ox(T)$
and the running cost $\ell$ is represented as
\begin{equation}\label{reprlagr}
\ell(t,z,\dot{z})=\ell_1(t,z)+\ell_2(\dot{x})+\ell_3(t,\dot{u},\dot{b}),
\end{equation}
where $\ell_1$, $\ell_2$, and $\ell_3$ are locally Lipschitzian around the given local minimizer $\oz(\cdot)$ with respect to all but time variables.
Suppose that $\ell_3$ is differentiable in $(\dot{u},\dot{b})$ on $\R^n\times\R$ with Lipschitz continuous partial derivatives and that there is a constant $L>0$ such that we have the estimate
\begin{equation}\label{boundgrad}
\begin{split}
\|\partial\ell_2(t,\dot{x})\|&\le L\|\dot{x}\|\;\mbox{ for all }\;t\in[0,T]\;\mbox{ and }\;\dot{x}\in\R^{n}\\
\|\nabla\ell_3(t,\dot{u},\dot{b})\|&\le L\big(\|\dot{u}\|+\|\dot{b}\|\big)\;\mbox{ for all }\;t\in[0,T]\;\mbox{ and }\;(\dot{u},\dot{b})\in\R^{n+1}.
\end{split}
\end{equation}
Suppose finally that $\oz$ satisfies the estimates
\begin{equation}\label{ouobbdd}
\begin{split}
\limsup_{k\to\infty}\left\|\frac{\ou(h_k)-\ou(0)}{h_k}\right\|\le L,&\qquad\limsup_{k\to\infty}\left\|\frac{\ob(h_k)-\ob(0)}{h_k}\right\|\le L,\\
\limsup_{k\to\infty}\left\|\frac{\ou(T)-\ou(T-h_k)}{h_k}\right\|\le L,&\qquad\limsup_{k\to\infty}\left\|\frac{\ob(T)-\ob(T-h_k)}{h_k}\right\|\le L,
\end{split}
\end{equation}
where $h_k$ is the step of the uniform discrete mesh taken from {\rm(\ref{partition})}.

Then there exist a multiplier $\lm\ge 0$, an adjoint arc $p(\cdot)=(p^x,p^u,p^b)\colon[0,T]\to\R^{n+nm+m}$ absolutely
continuous on $[0,T]$ with $\dot{p}(\cdot)\in L^2[0,T]$, signed vector measures $\gamma=(\gamma_1,\ldots,\gamma_m)\in C^*([0,T];\R^m)$
and $\xi=(\xi_1,\ldots,\xi_m)\in C^*([0,T];\R^{mn})$, as well as functions $\big(w(\cdot),v(\cdot)\big)\in L^\infty((0,T);\R^n)\times L^2((0,T);\R^n)$ with
\begin{eqnarray}\label{co}
\big(w(t),v(t)\big)\in{\rm co}\,\partial\ell\big(t,\oz(t),\dot\oz(t)\big)\;\mbox{ for a.e. }\;t\in[0,T]
\end{eqnarray}
such that the following conditions are satisfied:

$\bullet$ The {\sc primal-dual dynamic relationships}:
\begin{equation}\label{etajk-defl}
\big\langle\ou_{i}(t),\ox(t)\big\rangle<\ob_{i}(t)\Longrightarrow\eta_{i}(t)=0\;\mbox{ for a.e. }\;t\in[0,T],\quad \;i=1,\ldots,m,
\end{equation}
\begin{equation}\label{orthogonall}
\eta_{i}(t)>0\Longrightarrow\big\la\lm v^x(t)-q^{x}(t),\ou_i(t)\big\ra=0\;\mbox{ for a.e. }\;t\in[0,T],\quad \;i=1,\ldots,m,
\end{equation}
\begin{equation}\label{hamiltonx}
\dot p(t)=\lm w(t)+\big(0,\big[-\eta(t),\lm v^x(t)-q^x(t)\big],0\big)\;\mbox{ for a.e. }\;t\in[0,T],
\end{equation}
\begin{equation}\label{hamilton2ub}
q^u(t)=\lm\frac{\partial}{\partial\dot u}\ell_2\big(\dot\ou(t),\dot\ob(t)\big),\;
\mbox{ and }\;q^b(t)=\lm\frac{\partial}{\partial\dot b}\ell_2\big(\dot\ou(t),\dot\ob(t)\big)\;\mbox{ for a.e. }\;t\in[0,T],
\end{equation}
where $\eta(\cdot)=(\eta_1(\cdot),\ldots,\eta_m(\cdot))$ with the components $\eta_i(\cdot)\in L^{2}([0,T];\R_+)$, $i=1,\ldots,m$,
is a uniquely defined vector function determined by the representation
\begin{equation}\label{etajkl}
\dot\ox(t)=-\sum\limits_{i=1}^m\eta_{i}(t)\ou_{i}(t)\;\mbox{ for a.e. }\;t\in[0,T],
\end{equation}
and where $q\colon[0,T]\to\R^{n+nm+m}$ is a function of bounded variation on $[0,T]$ with its left-continuous
representative given, for all $t\in[0,T]$ excepting at most a countable subset, by
\begin{equation}\label{q}
q(t)=p(t)-\Big({\int_{[t,T]}}\sum\limits_{i=1}^m\ou_i(s)d\gamma_i(s),\int_{[t,T]}
\big[{\rm rep}_m(\ox(s)),d\gamma(s)\big]+2\int_{[t,T]}\big[\ou(s),d\xi(s)\big],-\int_{[t,T]}d\gamma(s)\Big).
\end{equation}

$\bullet$ The {\sc transversality conditions} at the right endpoint:
\begin{equation}\label{pxkkc}
-p^x(T)\in\lm\partial\ph\big(\ox(T)\big)+\sum\limits_{i=1}^m p^{b}_{i}(T)\ou_i(T),
\end{equation}
\begin{equation}\label{pk-1c}
p_i^u(T)+p_i^b(T)\ox(T)=\Big\la p_i^u(T)+p_i^b(T)\ox(T),\ou_i(T)\Big\ra\ou_i(T),\quad i=1,\ldots,m,
\end{equation}
\begin{equation}\label{alphak1c}
p^{b}(T)\in\R^m_+\mbox{ and }\;\big\langle\ou_{i}(T),\ox(T)\big\rangle<\ob_{i}(T)\Longrightarrow p^{b}_{i}(T)=0,\quad i=1,\ldots,m.
\end{equation}

$\bullet$ The {\sc measure nonatomicity conditions}:

$\quad$ {\bf (a)} If $t\in[0,T)$ and $\langle\ou_{i}(t),\ox(t)\rangle<\ob_{i}(t)$ for all $i=1,\ldots,m$, then there is a neighborhood $V_t$ of $t$ in
$[0,T]$ such that $\gamma(V)=0$ for any Borel subset $V$ of $V_t$.

$\quad$ {\bf (b)} If $t\in [0,\tau)\cup(T-\tau,T]$ and $1/2 <\|\ou_{i}(t)\|<3/2$ for all $i=1,\ldots,m$,
then there is a neighborhood $W_t$ of $t$ in $[0,\tau)\cup(T-\tau,T]$ such that $\xi(W)=0$ for any Borel subset $W$ of $W_t$.

$\bullet$ {\sc Nontriviality conditions}: We always have
\begin{equation}\label{nontriv2}
\lm +\|q(0)\|+\|p(T)\|\ne 0.
\end{equation}
Furthermore, the additional assumptions on $\la\ou_i(0),\ox(b)\ra<\ob_i(0)$ and $1/2<\|\ou_i(0)\|<3/2$
whenever $i=1,\ldots,m$ ensure the validity of the enhanced nontriviality condition $(\lm,p(T))\ne 0$.
\end{Theorem}
{\bf Proof.} First we construct all the functions with the claimed properties satisfying the {\em primal-dual dynamic relationships} of the theorem.
Fix any $\t\in(0,T)$ and for the given r.i.l.m.\ $\oz(\cdot)$ in $(P^{\tau})$ consider the discrete approximation problems $(P_k^{\tau})$ whose optimal
solutions $\oz^k=(\ox^k,\ou^k,\ob^k)$ exist by Proposition~\ref{ex-disc} with their piecewise linear extensions $\oz^k(\cdot)$, $0\le t\le T$,
converging to $\oz(\cdot)$ in the sense of Theorem~\ref{conver}. Our aim is to derive the claimed necessary optimality conditions
for $\oz(\cdot)$ by passing to the limit from those for $\oz^k(\cdot)$ obtained in Theorem~\ref{R1}. To proceed, for each $k\in\N$
denote by $w^k(\cdot)$ and $v^k(\cdot)$ the piecewise constant extensions to $[0,T]$ of the discrete functions $w^k_j$ and $v^k_j$,
respectively, satisfying \eqref{subl}. It follows from \eqref{subl} and \eqref{boundgrad} due to the strong $W^{1,2}$-convergence
of the sequence $\{\oz^k\}$ to $\oz$ and the local Lipschitz continuity of $\ell(t,\cdot,\cdot)$ that $\{(w^k(t),v^k(t))\}$ is weakly compact in
$L^2([0,T];\R^{2(n+nm+m)})=:L^2[0,T]$. Hence we have
\begin{equation*}
\big(w^k(\cdot),v^k(\cdot)\big)\rightarrow\big(w(\cdot),v(\cdot)\big))\;\mbox{ weakly in }\;L^2[0,T]\;\mbox{ as }\;k\to\infty
\end{equation*}
with some pair $\big(w(\cdot),v(\cdot)\big)\in L^2[0,T]$. Employing the aforementioned robustness property of the
subdifferential together with the well-known weak convergence result based on Mazur's theorem (see, e.g., \cite[Theorem~1.4.1]{AC})
allows us to deduce from \eqref{subl} that the convexified inclusion in \eqref{co} holds. Note also that $w(\cdot)$ belongs actually to $L^\infty[0,T]$ due to
its a.e.\ boundedness on $[0,T]$.

Further, based on \eqref{th} for all $k\in\N$ we define the functions
\begin{equation*}
\th^{xk}(t):=\frac{\th^{xk}_j}{h_k}\;\mbox{ for }\;t\in[t^k_j,t^k_{j+1}),\quad j=0,\ldots,k-1,
\end{equation*}
on $[0,T]$ and easily observe by the convexity of the integrand that
\begin{equation}\label{converthetax}
\begin{split}
\disp\int_0^T\|\th^{xk}(t)\|^2dt&=\disp\sum_{j=0}^{k-1}\frac{\|\th^{xk}_{j}\|^2}{h_k}\le\frac{4}{h_k}
\sum_{j=0}^{k-1} \left(\int_{t^k_j}^{t^k_{j+1}}\Big\|\dot{\ox}(t)-\frac{\ox^k_{j+1}-\ox^k_j}{h_k}\Big\|dt\right)^2\\
&\le 4\sum_{j=0}^{k-1}\int_{t^k_j}^{t^k_{j+1}}\Big\|\dot{\ox}(t)-\frac{\ox^k_{j+1}-\ox^k_j}{h_k}\Big\|^2dt=
4\disp\int_0^T\|\dot\ox(t)-\dot\ox^k(t)\|^2\,dt\to 0\;\mbox{ as }\;k\to\infty,
\end{split}
\end{equation}
where the convergence is due to Theorem~\ref{conver}. This implies that a subsequence of $\{\th^{xk}(t)\}$
converges to zero a.e.\ on $[0,T]$. The same conclusions hold for the similarly defined functions $\th^{uk}(t)$ and $\th^{bk}(t)$ on $[0,T]$.

It follows from \eqref{etajk} that for the piecewise linear interpolations of $\ox^k(\cdot)$ and $\ou^k(\cdot)$ on $[0,T]$ we have
\begin{equation}\label{he1}
\overset{.}{\bar{x}}\text{\thinspace}^{k}\left(t\right)=-\sum_{i=1}^{m}\eta_{ji}^{k}\bar{u}_{i}^{k}(t_{j}^{k})\;
\mbox{ for all }\;t\in(t_{j}^{k},t_{j+1}^{k}),\quad j=0,\ldots,k-1.
\end{equation}
Now extend $\eta_{j}^{k}$ to $[0,T]$ by $\eta^{k}\left(t\right):=\eta_{j}^{k}$ for $t\in\lbrack t_{j}^{k},t_{j+1}^{k})$ and define the functions
\begin{equation}\label{deftau}
\vt^{k}(t):=\max\big\{t_{j}^{k}\big|\;t_{j}^{k}\le t,\;0\le j\le k\big\}\;\mbox{ for all }\;
t\in\left[0,T\right],\;k\in\N.
\end{equation}
We clearly have $\vt^{k}(t)=t_{j}^{k}$ for all $t\in\lbrack t_{j}^{k},t_{j+1}^{k})$, $j=0,\ldots,k-1$,
and $\vt_k (t)\to t$ uniformly in $[0,T]$ as $k\to\infty$.
The uniform convergence of $\oz^k(\cdot)$ to $\oz(\cdot)$ on $[0,T]$ readily implies that
\begin{equation}\label{claim}
\bar{z}^{k}\big(\vt^{k}(t)\big)\to\bar{z}(t)\;\mbox{ uniformly on }\;[0,T]\;\mbox{ as }\;k\to\infty.
\end{equation}
This notation allows us to rewrite \eqref{he1} as
\begin{equation}\label{lincomb}
\overset{.}{\bar{x}}\text{\thinspace}^{k}\left(t\right)=-\sum_{i=1}^{m}\eta_{i}^{k}\left(t\right)
\bar{u}_{i}^{k}(\vt^{k}(t))\;\mbox{ for all }\;t\in\left[0,T\right]\backslash\big\{t_{0}^{k},\ldots,t_{k}^{k}\big\}.
\end{equation}

Consider further the subset of $[0,T]$ given by
\begin{equation*}
\mathcal{T}:=\left[0,T\right]\backslash\bigcup\limits_{k\in\N}\big\{t_{0}^{k},\ldots,t_{k}^{k}\big\}.
\end{equation*}
For any fixed $t\in\mathcal{T}$ denote by $J:=I(\bar{x}(t),\bar{u}(t),\bar{b}(t))$ the collection of active constraint indices from
\eqref{active_constr} and by $\tilde{u}^{k}(t)$ the matrix consisting of the rows $\bar{u}_{i}^{k}(\vt^{k}(t))$, $i\in J$,
while $\tilde{u}(t)$ stands for the matrix consisting of the rows $\bar{u}_{i}(t)$ as $i\in J$.
More precisely, this means that $\tilde{u}_{i}^{k}(t)=\bar{u}_{\phi (i)}^{k}(\vt^{k}(t))$ and
$\tilde{u}_{i}=\bar{u}_{\phi(i)}(t)$ for all $i\in\left\{1,\ldots,|J|\right\}$, where $|J|$ signifies the cardinality of $J$,
and where the mapping $\phi:\left\{1,\ldots,|J|\right\}\longrightarrow J$ is a bijection.
The assumed LICQ condition tells us that the rows of $\bar{u}_{i}(t)$ are linearly independent for $i\in J$,
and consequently we can build the generalized inverse matrix
\begin{equation*}
\hat{u}(t):=\big[\tilde{u}(t)\tilde{u}(t)^*\big]^{-1}\tilde{u}(t)\;\mbox{ for each }\;t\in{\cal T}.
\end{equation*}
It follows from (\ref{claim}) that $\bar{u}_{i}^{k}(\vt^{k}(t))\rightarrow_{k}\bar{u}_{i}(t)$ for $i\in J$, and so the generalized inverse
\begin{equation} \label{geninv}
\hat{u}^{k}(t):=\big[\tilde{u}^{k}(t)\tilde{u}^{k}(t)^*]^{-1}\tilde{u}^{k}(t)
\end{equation}
is well defined for all $k$ sufficiently large with $\hat{u}^{k}(t)\rightarrow\hat{u}(t)$ as $k\to\infty$. We have by the definition of $J$
that $\left\langle \bar{u}_{i}(t),\bar{x}(t)\right\rangle<\bar{b}_{i}(t)$ whenever $i\in\{1,\ldots,m\}\setminus J$, and hence (\ref{claim}) tells us that
\begin{equation*}
\left\langle\bar{u}_{i}^{k}(\vt^{k}(t)),\bar{x}^{k}(\vt^{k}(t))\right\rangle<\bar{b}_{i}^{k}\big(\vt^{k}(t)\big)\;\mbox{ for }\;
i\in J^c:=\big\{1,\ldots,m\}\setminus J
\end{equation*}
when $k$ is large. Since $\vt^{k}(t)=t_{j}^{k}$ for some $j\in\{0,\ldots,k\}$, we deduce that
$\left\langle\bar{u}_{ji}^{k},\bar{x}_{j}^{k}\right\rangle <\bar{b}_{ji}^{k}$ for all $i\in J^c$,
and so \eqref{etajk-def} yields $\eta_{ji}^{k}=0$ for this $j$ and all $i\in J^c$.
Remembering that $\eta_{ji}^{k}=\eta_{i}^{k}(t)$ by construction, we conclude that $\eta_{i}^{k}(t)=0$ for all $i\in J^c$ and large $k$.
This allows us to rewrite (\ref{lincomb}) as
\begin{equation}\label{lincomb2}
\overset{.}{\bar{x}}\text{\thinspace}^{k}\left(t\right)=-\sum_{i\in J}\eta_{i}^{k}\left(t\right)\bar{u}_{i}^{k}(\vt^{k}(t))=
-\tilde{u}^{k}(t)^*\tilde{\eta}^{k}\left(t\right)\;\mbox{ for large }\;k,
\end{equation}
where $\tilde{\eta}^{k}\left(t\right)$ collects the components $\eta_{i}^{k}\left(t\right)$ for $i\in J$, i.e.,
$\tilde{\eta}_{i}^{k}\left(t\right)=\eta_{\phi(i)}^{k}\left(t\right) $ for $i\in\left\{1,\ldots,|J|\right\}$ via the above bijection $\phi$.
Thus (\ref{geninv}) and \eqref{lincomb2} ensure the representation
\begin{equation}\label{repretak}
\tilde{\eta}^{k}\left(t\right)=-\hat{u}^{k}(t)\overset{.}{\bar{x}}^{k}\left(t\right)\;\mbox{ for large }\;k,
\end{equation}
and the passage to the limit implies that $\tilde{\eta}^{k}\left(t\right)\rightarrow-\hat{u}(t)\overset{.}{\bar{x}}\left(t\right)$
as $k\to\infty$.
Define now the required function $\eta(t)=(\eta_1(t),\ldots,\eta_m(t))$ for all $t\in J$ of full measure on $[0,T]$ by
\begin{eqnarray}\label{eta}
\eta_{\phi(i)}\left(t\right):=-\hat{u}_{i}(t)\overset{.}{\bar{x}}\left(t\right)\;\mbox{ for }\;i\in
\big\{1,\ldots,|J|\big\}\;\mbox{ and }\;\eta_{i}\left(t\right):=0\;\mbox{ for }\;i\in J^{c}
\end{eqnarray}
and observe from the constructions above that $\eta^{k}\left(t\right)\to\eta\left(t\right)$ as $k\to\infty$.
Since $\eta_{i}^{k}\left(t\right)=0$ for $i\in J^{c}$ and large $k$, we get \eqref{etajk-defl} by the definition of
$J$ and also $\eta_i(t)\ge 0$ for all $i=1,\ldots,m$ due to $\eta^k_j\in\R^m_+$ in Theorem~\ref{R1}.
It follows from (\ref{lincomb2}) by passing to the limit that \eqref{etajkl} holds.
The uniqueness of $\eta(t)$ in \eqref{etajkl} for a.e.\ $t\in[0,T]$ follows from the imposed LICQ condition.

To justify the claimed properties of $\eta(\cdot)$ in \eqref{etajk-defl}, it remains to show that $\eta(\cdot)\in L^2[0,T]$.
To see it, let us rearrange for each $t\in{\cal T}$ the active components of $\eta(t)$ by putting $\tilde\eta_i(t):=\eta_{\phi(i)}(t)$ for $i\in\{1,\ldots,|J|\}$
and $\tilde\eta_i(t)=0$ otherwise. Since $\tilde\eta_i(t)=-\hat u_i(t)\dot{\ox}(t)$ whenever $i\in\{1,\ldots,|J|\}$ by \eqref{eta} and the convergence
$\hat u^k(t)\to\hat u(t)$ is uniform on $[0,T]$ by \eqref{claim}, it follows that $\tilde\eta^k(\cdot)\to\tilde\eta(\cdot)$ and hence $\eta^k(\cdot)\to\eta(\cdot)$
{\em strongly} in $L^2[0,T]$. This not only verifies that $\eta(\cdot)\in L^2[0,T]$, but also allows us to get the estimate
\begin{equation}\label{boundetal1}
h_k\sum_{j=0}^k\|\eta^k_j\|^2=\int_0^T\|\eta^k(t)\|^2dt\le M
\end{equation}
with some constant $M>0$ independent of $k$. It immediately follows from \eqref{boundetal1} and \eqref{converthetax} that
\begin{equation}\label{qxbv}
\int_0^T\big\|\eta^k_i(t)\theta^{xk}(t)\big\|\,dt\to 0\;\mbox{ as }\;k\to\infty\;\mbox{ for all }\;i=1,\ldots,m.
\end{equation}

Next we use the notation of Theorem~\ref{R1} and define $q^k(\cdot)=(q^{xk}(\cdot),q^{uk}(\cdot),q^{bk}(\cdot))$
by extending $p_j^k$ piecewise linearly to $[0,T]$ with $q^{k}(t^k_j):=p^{k}_j$ for $j=0,\ldots,k$.
Construct $\gg^k(\cdot)$, $\xi^k(\cdot)$ on $[0,T]$ by
\begin{eqnarray}\label{gxt}
\gamma^k(t):=\gamma^k_j,\;\xi^k(t):=\frac{1}{h_k}\xi^k_j\;\mbox{ for }\;t\in[t^k_j,t^k_{j+1}),\quad j=0,\ldots,k-1
\end{eqnarray}
with $\gamma^k(T):=0$ and $\xi^k(T):=\xi^k_k$. Appealing to $\vt^k(\cdot)$ in \eqref{deftau},
equations \eqref{p_iteration}--\eqref{vstar} can be rewritten as
\begin{equation}\label{tranm}
\dot q^{xk}(t)-\lambda^k w^{xk}(t)=\sum\limits_{i=1}^m\gamma^k_i(t)\ou_i^k\big(\vt^k(t)\big),
\end{equation}
\begin{equation}\label{eem}
\dot q_i^{uk}(t)-\lambda^k w_i^{uk}(t)=2\xi^k_i(t)\ou^k_i\big(\vt^k(t))+\gamma^k_i(t)\ox^k(\vt^k(t)\big)-
\eta^k_i(t)\big(\lm^k v^{xk}(t)+\lm^k\th^{xk}(t)-q^{xk}(\vt^k_+(t))\big),
\end{equation}
\begin{equation}\label{vstarm}
\dot q^{bk}(t)-\lambda^k w^{bk}(t)=-\gamma^k(t)
\end{equation}
for every $t\in(t^k_j,t^k_{j+1})$, $j=0,\ldots,k-1$, and $i=1,\ldots,m$, where $\vt^k_+(t):=t^k_{j+1}$ for $t\in[t^k_j,t^k_{j+1})$.

Now we define $p^k(\cdot)=\Big(p^{xk}(\cdot),p^{uk}(\cdot),p^{bk}(\cdot)\Big)$ on $[0,T]$ by setting
\begin{equation}\label{pqi}
\begin{array}{ll}
p^{k}(t):=&q^{k}(t)+\disp\int_{t}^{T}\Big(\sum\limits_{i=1}^m\gamma^k_i(s)
                \ou_i^k\big(\vt^k(s)\big),\;2\big[\xi^k(s),\ou^k\big(\vt^k(s)\big)\big]\\
&+\disp\big[\gamma^k(s),{\rm rep}_m\big(\ox^k(\vt^k(s))\big)\big],\;-\gamma^k(s)\Big)ds
\end{array}
\end{equation}
for every $t\in [0,T]$. This gives us $p^k(T)=q^k(T)$ with the differential relation
\begin{equation}\label{pq}
\dot p^{k}(t)=\dot q^{k}(t)-\Big(\sum\limits_{i=1}^m\gamma^k_i(t)\ou_i^k\big(\vt^k(t)\big),\;2\big[\xi^k(t),\ou^k\big(\vt^k(t)\big)\big]+
\big[\gamma^k(t),{\rm rep}_m\big(\ox^k(\vt^k(t))\big)\big],\;-\gamma^k(t)\Big)
\end{equation}
holding for a.e.\ $t\in[0,T]$. Using this notation, equations \eqref{tranm}--\eqref{vstarm} can be rewritten as
\begin{equation}\label{tranmp}
\dot p^{xk}(t)=\lambda^k w^{xk}(t),
\end{equation}
\begin{equation}\label{eemp}
\dot p_i^{uk}(t)=\lambda^k w_i^{uk}(t)-\eta^k_i(t)\Big(\lm^k v^{xk}(t)+\lm^k\th^{xk}(t)-q^{xk}\big(\vt^k_+(t)\big)\Big),
\end{equation}
\begin{equation}\label{vstarmp}
\dot p^{bk}(t)=\lambda^k w^{bk}(t)
\end{equation}
for every $t\in(t^k_j,t^k_{j+1})$, $j=0,\ldots,k-1$, and $i=1,\ldots,m$. Define the vector measures $\gamma^k_{mes}$ and $\xi^k_{mes}$ by
\begin{equation}\label{gammaxi}
\int_B d\gamma^k_{mes}:=\int_B\gamma^k(t)\;\mbox{ and }\;\int_B d\xi^k_{mes}:=\int_B \xi^k(t)dt\;\mbox{ for every Borel subset }\;B\subset[0,T].
\end{equation}
From now on we drop for simplicity the index `mes' in the measure notation if no confusion arises.

Observe next that all the expressions in the statement of Theorem~\ref{R1} are positively homogeneous of degree
$1$ with respect to $\lambda^k$, $p^k$, $\gamma^k$, and $\xi^k$. Therefore the nontriviality condition \eqref{nontriv1}
allows us to normalize them by imposing the following relationships whenever $k\in\N$:
\begin{equation}\label{alphak}
\lm^k+\|q^{uk}(0)\|+\|q^{bk}(0)\|+\|p^k(T)\|+\int_0^T\|\gamma^{k}(t)\|dt+\int_0^T\|\xi^{k}(t)\|dt=1,
\end{equation}
which tell us that all the sequential terms in \eqref{alphak} are uniformly bounded. Passing below to subsequences of $k\to\infty$ if necessary, we can immediately
conclude that $\lambda^k\to\lambda$ for some $\lambda\ge 0$. Then the equality $p^k(T)=q^k(T)$, the uniform boundedness of the first integral terms in
\eqref{alphak} and of $\{w^k(\cdot)\}$ in $L^\infty[0,T]$ implies via \eqref{tranm} that the sequence $\{q^{xk}(\cdot)\}$ has uniformly bounded variations on $[0,T]$,
and so it is bounded in $L^\infty[0,T]$. Observe further that the right-hand sides of \eqref{tranmp} and \eqref{vstarmp} are obviously uniformly bounded in
$L^\infty[0,T]$. Concerning the right-hand sides of \eqref{eemp}, observe that the sequence  $\{\lambda\eta^k v^k\}$ is uniformly bounded in $L^2$, while the
remaining summands are uniformly bounded in $L^1$ due to \eqref{qxbv} and the uniform boundedness of $\{q^{xk}\}$. Since we also have their uniform integrability by
the arguments above, the classical Dunford-Pettis theorem on the weak compactness in $L^1$ allows us a subsequence of $\{\dot p^k(\cdot)\}$, which weakly
converges in $L^1$ to some function generating by the Newton-Leibniz formula an  absolutely continuous function $p(\cdot)$ such that $p^k(t)\to p(t)$ uniformly on
$[0,T]$. Moreover, the aforementioned Mazur's theorem gives us a subsequence of convex combinations of $\dot{p}^k(t)$ converging to $\dot p(t)$ a.e.\ pointwise on
$[0,T]$.

Using the uniform boundedness of $\int_0^T\|\gamma^{k}(t)\|dt$ and $\int_0^T\|\xi^{k}(t)\|dt$ by \eqref{alphak}, the relationships in \eqref{tranm}--\eqref{vstarm}
together with $p^k(T)=q^k(T)$ and \eqref{qxbv} ensures that $\{q^k(\cdot)\}$ is of uniformly bounded variation. This allows us to employ Helly's selection
theorem and find measures $\gamma\in C^*([0,T];\R^m)$, $\xi\in C^*([0,T];\R^{mn})$ and a function of bounded variation $q(\cdot)$ on $[0,T]$ such that a subsequence
of $\{q^k(\cdot)\}$ pointwise converges to $q(\cdot)$ while some of $\{(\gamma^k,\xi^k)\}$ weak$^*$ converges to $(\gamma,\xi)$ in $BV[0,T]$;
see \cite[Definition~3.11, Theorem~3.23, and Proposition~3.21]{AFP}. Thus having $q(T)=p(T)$ and combining it with the a.e.\ pointwise convergence of convex
combinations of $\dot p^k(\cdot)$ to $\dot p(\cdot)$ justify the possibility of passing to the limit in \eqref{tranmp}--\eqref{vstarmp} and to verify
\eqref{hamiltonx}. Combining \eqref{orthogonal} and \eqref{converthetax} gives us \eqref{orthogonall} while \eqref{hamilton2ub} follows from \eqref{psiu} and
\eqref{converthetax}.

Now we intend to prove the representation \eqref{q} for $q(\cdot)$ by passing to the limit in \eqref{pqi}.
It follows from the norm convergence of $\oz^k(\cdot)\to\oz(\cdot)$ in $W^{1,2}[0,T]$ that $\{\dot\ox^k(\cdot)\}$
is bounded in $L^2[0,T]$ and also that
\[
{\rm sup}_{t\in[0,T]}\big\|\ou^k\big(\vt^k(t)\big)-\ou^k(t)\big\|=\max_{0\le j\le k-1}\big\|\ou^k_{j+1}-\ou^k_j\big\|\to 0\;\mbox{ as }\;k\to\infty
\]
with the same for the $b$-components. Furthermore, for any fixed $i\in\{1,\ldots,m\}$ we have the estimate
\begin{equation*}
\begin{split}
\lefteqn{\!\!\!\!\!\!\!\!\!\!\!\!\!\!\!\Big\|\int_t^T\gamma^k_i(s)\ou_i^k\big(\vt^k(s)\big)ds-\int_t^T\ou_i(s)d\gamma_i(s)\Big\|}\\
&\;\le\Big\|\int_t^T\gamma^k_i(s)\ou_i^k\big(\vt^k(s)\big)ds-\int_t^T\gamma^k_i(s)\ou_i(s)ds\Big\|+\Big\|\int_t^T\gamma^k_i(s)\ou_i(s)ds-\int_t^T
\ou_i(s)d\gamma_i(s)\Big\|\\
&\;=\Big\|\int_t^T\gamma^k_i(s)\big(\ou_i^k(\vt^k(s))-\ou_i(s) \big)ds\Big\|+\Big\|\int_t^T\ou_i(s)\gamma^k_i(s)ds-\int_t^T\ou_i(s)d\gamma_i(s)\Big\|,
\end{split}
\end{equation*}
where the first summand vanishes as $k\to\infty$ on $[0,T]$ due to the uniform convergence of $\ou^k_i(\cdot)$ to $\ou_i(\cdot)$
and the uniform boundedness of $\int_0^T \|\gamma^{k}(t)\|dt$, while the second one vanishes for all $t\in[0,T]$ except
for at most countably many points because the measures $\gamma^k$ converge weak$^*$ to $\gamma$; see, e.g., \cite[p.\ 325]{v}. Hence
$$
\int_t^T\gamma^k_i(s)\ou_i^k\big(\vt^k(s)\big)ds\to\int_t^T\ou_i(s)d\gamma_i(s)\;\mbox{ as }k\to \infty
$$
except for at most countably many $t$. Proceeding in the same way with the rest of \eqref{pqi} gives us \eqref{q},
which is the left-continuous representative of the limiting function $q(\cdot)$ of bounded variation on $[0,T]$;
see, e.g., \cite[Theorem~3.8]{AFP}. Using this and the properties of $w(\cdot), v(\cdot)$ and $\eta(\cdot)$
established above we complete proving all the statements in the primal-dual dynamic relationships of the theorem.

Next we justify the {\em transversality conditions}, which is a much easier task. Indeed, the validity of \eqref{pxkkc} and \eqref{alphak1c} follows
by passing to the limit in \eqref{pxkk}
and \eqref{alphak1}, respectively, with taking into account that $\{p^k(\cdot)\}$ converges uniformly to
$p(\cdot)$ and that $p^k(T)=q^k(T)=p^k_k$. Then observe from \eqref{pk-1} that
\begin{equation*}
p^{uk}_{ki}+p^{bk}_{ki}\ox^k_k=\big\la p^{uk}_{ki}+p^{bk}_{ki}\ox^k_k,\ou^k_k\big\ra\ou^k_k\;\mbox{ for all }\;i=1,\ldots,m,
\end{equation*}
which gives us \eqref{pk-1c} by passing to the limit as $k\to\infty$.

Now we proceed with verifying the {\em measure nonatomicity conditions} of the theorem.
To check the one in (a), fix $t\in[0,T)$ with $\la\ou_{i}(t),\ox(t)\ra<\ob_{i}(t)$ for $i=1,\ldots,m$
and find a neighborhood $V_t$ of $t$ in $[0,T]$ such that for any $s\in V_t$ we have
$\la\ou_{i}(s),\ox(s)\ra<\ob_{i}(s)$, $i=1,\ldots,m$. This yields $\langle\ou^k_{i}(t_j^k),\ox^k(t_j^k)\big\rangle<\ob^k_{i}(t_j^k)$, $i=1,\ldots,m$,
when $k$ is sufficiently large and so $t_j^k\in V_t$. Thus it follows from \eqref{measuresign} that for any Borel subset
$V\subset V_t$ we have $\gamma^k(t)=0$ on $V$, which implies in turn that
$\|\gamma_{mes}^k\|(V)=\int_V d\|\gamma_{mes}^k\|=\int_V\|\gamma^k(t)\|dt=0$. Letting $k\to\infty$ shows that $\|\gamma\|(V)=0$.
The measure nonatomicity condition (b) for $\xi$ is justified similarly.

Our final step is to prove the {\em nontriviality conditions} starting with \eqref{nontriv2}.
Arguing by contradiction, suppose that $\lm=0$, $q(0)=0$, and $p(T)=0$ and hence get
$\lm^k\to 0$, $q^k(0)\to 0$, and $p^k(T)=p^k_k\to 0$ as $k\to\infty$. Substituting \eqref{psiu} into \eqref{vstar}
with the usage of \eqref{th}, \eqref{subl}, and \eqref{reprlagr}, we get
\begin{align*}
\int_0^T\|\gamma^k(t)\|dt&=\sum\limits_{j=0}^{k-1} h_k\|\gamma^k_j\|\le\;\sum\limits_{j=1}^{k-1}\Big\|p^{bk}_{j+1}-p^{bk}_{j}\Big\|+
\lm^k\sum \limits_{j=0}^{k-1} h_k\|w_j^{bk}\|+\|p^{bk}_{1}-p^{bk}_{0}\|\\
&\le\lm^k\sum\limits_{j=1}^{k-1}\Big\|\frac{\theta^{bk}_{j}-\theta^{bk}_{j-1}}{h_k}\Big\|+\lm^k\sum \limits_{j=1}^{k-1}\|v_j^{bk}-v_{j-1}^{bk}\|
+\lm^k\sum\limits_{j=0}^{k-1} h_k \|w_j^{bk}\|+\|p^{bk}_{1}\|+\|p^{bk}_{0}\|\\
&\le 2\lambda^k\Bigg(\sum_{j=0}^{k-2}\Big\|
\frac{\ob^k_{j+2}-2\ob^k_{j+1}+\ob^k_j}{h_k}\Big\|+\frac{1}{h_k}\Big|\int_{0}^{h^k}\dot{\ob}(t)\,dt\Big|
+\frac{1}{h_k}\Big|\int_{T-h^k}^{T}\dot{\ob}(t)\,dt\Big|\\
&\qquad+\Tilde L\sum\limits_{j=0}^{k-2}\Big(\Big\|\frac{\ob^k_{j+2}-2\ob^k_{j+1}+\ob^k_j}{h_k}\Big\|
+\Big\|\frac{\ou^k_{j+2}-2\ou^k_{j+1}+\ou^k_j}{h_k}\Big\|\Big)\Bigg)\\
&\qquad+\lm^k\sum\limits_{j=0}^k h_k\|w_j^{bk}\|+\|p^{bk}_{1}\|+\|p^{bk}_{0}\|,
\end{align*}
where $\tilde L$ is a Lipschitz constant of $\nabla\ell_2$. Deducing further from \eqref{psiu} that
\begin{equation*}
\|p^{bk}_{1}\|\le\lambda^k\|v^{bk}_0\|+\lambda^k h_k^{-1}|\theta^{bk}_0|,
\end{equation*}
let us show that $\|p^{bk}_{1}\|\to 0$ as $k\to\infty$. Indeed, we get by \eqref{subl} and the smoothness of $\ell_2$
that $v^{bk}_0=\frac{\partial\ell_2}{\partial \dot b}\big(0,\frac{\ou_1^k-\ou_0^k}{h_k},\frac{\ob_1^k-\ob_0^k}{h_k}\big)$.
Combining this with \eqref{initialubconstraint1} and \eqref{boundgrad}
ensures the boundedness of $\{v^{bk}_0\}$ and hence $\lambda^k\|v^{bk}_0\|\to 0$. It follows from \eqref{discretbound1},
\eqref{initialubconstraint1}, and \eqref{th} that $\{h_k^{-1}\theta^{bk}_0\}$ is bounded and thus $\lambda^k h_k^{-1}|\theta^{bk}_0|\to 0$,
which justifies the claim on $\|p^{bk}_{1}\|$. We also have by construction that $p^{bk}_{0}=q^k(0)\to 0$. Involving again \eqref{discretbound1} together with
\eqref{bvbounded} and the boundedness of $\{w^{bk}\}$ by \eqref{subl} gives us
$\int_0^T\|\gamma^k(t)\|dt\to 0$ as $k\to\infty$.

Considering now the functions $\xi^k(t)$ from \eqref{gxt}, we have
\begin{equation*}
\int_0^T\|\xi^k(t)\|dt=\sum\limits_{j=0}^{k-1}h_k\frac{1}{h_k}\|\xi^k_j\|=\sum\limits_{j=0}^{k-1}\|\xi^k_j\|,\quad k\in\N.
\end{equation*}
It follows from \eqref{ee} due to \eqref{contr-disc} that the estimate
\begin{equation}\label{uu1}
\begin{split}
\sum\limits_{j=0}^{k-1}\|\xi^k_{ji}\|&\le\sum\limits_{j=0}^{k-1}\Big\|p^{uk}_{(j+1)i}-p^{uk}_{ji}\Big\|+
\sum\limits_{j=0}^{k-1}\lm^k h_k\|w_{ji}^{uk}\|\\
&\quad+\sum\limits_{j=0}^{k-1}h_k|\eta^k_{ji}|\,\Big\|\lambda^kv^{xk}_j+\frac{\lambda^k\theta^{xk}_{j}}{h_k}-p^{xk}_{j+1}\Big\|+
\sum \limits_{j=0}^{k-1}h_k\|\gamma^k_{ji}\|\cdot\|\ox^k_j\|
\end{split}
\end{equation}
holds whenever $i=1,\ldots,m$. Furthermore, it follows from \eqref{psiu} that
\begin{equation}\label{uu2}
\sum\limits_{j=0}^{k-1}\Big\|p^{uk}_{j+1}-p^{uk}_{j}\Big\|\le\lm^k\sum\limits_{j=1}^k\Big\|\frac{\theta^{uk}_{j}-\theta^{uk}_{j-1}}{h_k}\Big\|
+\lm^k\sum\limits_{j=1}^k\|v_{j}^{uk}-v_{j-1}^{uk}\|+\|p^{uk}_{1}-p^{uk}_{0}\|.
\end{equation}
Recalling that $\|q^k(T)\|=\|p^k(T)\|\to 0$ and $\int_0^T\|\gamma^k(t)\|dt\to 0$, we get $q^x(\cdot)=0$ on $[0,T]$ by passing to the limit in \eqref{tranm} due to the weak$^*$ convergence of $\gamma^k(\cdot)$ in $BV[0,T]$.  Then \eqref{pqi} yields $p^x(\cdot)=0$, and so
\begin{eqnarray}\label{max-conver}
\max\limits_{t\in[0,T]}\|q^{xk}(t)\|=\max\limits_{j=0,\ldots,k}\|p^{xk}_j\|\to 0\;\mbox{ as }\;k\to\infty.
\end{eqnarray}
Combining \eqref{max-conver} with \eqref{th}, \eqref{converthetax}, \eqref{boundetal1},
and the $L^2$--boundedness of $\{v^{xk}\}$, which follows from \eqref{boundgrad} and the strong $W^{1,2}$--convergence of
Theorem \ref{conver}, tells us that the third summand in \eqref{uu1} vanishes as $k\to\infty$. Using the same arguments allowing us to prove that $\int_0^T\|\gamma^k(t)\|dt\to 0$ as $k\to\infty$, we get by \eqref{uu2} that the first summand in \eqref{uu1} vanishes as well, which
therefore verifies that $\int_0^T\|\xi^k(t)\|dt\to 0$. All of this leads us to the violation of \eqref{alphak}
and thus justifies the nontriviality condition \eqref{nontriv2}.

To verify finally the {\em enhanced} nontriviality condition under the additional assumptions made, suppose by contradiction that $(\lambda,p(T))=0$. By \eqref{hamiltonx}, $p^x(t)=0$, $p^b(t)=0$ for all $t\in[0,T]$. By \eqref{hamilton2ub}, $q^u(t)=0$, $q^b(t)=0$ for almost all $t \in [0,T]$. Combining those arguments and \eqref{q}, we get that also $q^x(t)=0$ for almost all $t \in [0,T]$. Using \eqref{hamiltonx} again yields that $p^u(t)=0$ for all $t \in [0,T]$.
Therefore, $p(t)=0$ for all $t \in [0,T]$ and $q(t)=0$ for almost all $t \in [0,T]$. By using the measure nonatomicity condition, we get also $q(0)=0$, hence contradicting the nontriviality condition \eqref{nontriv2} and thus completing the proof. $\h$

It is worth mentioning (as used in Example~\ref{ex3} below) that the differentiability assumption on $\ell$
with respect to $(\dot u,\dot b)$ can be replaced in the proof of Theorem~\ref{necopt} by the following:
there is $M>0$ such that for all the partitions $0<t_0<t_1<\ldots<t_k<T$ and
$(v^u_j,v^b_j)\in\partial_{\dot u,\dot b}\ell(t_j,\oz(t_j),\dot\oz(t_j))$, $j=0,\ldots,k$, we have
\begin{equation}\label{bvderivative}
\sum_{j=0}^{k-1}\|v^u_{j+1}-v^u_{j}\|\le M\;\mbox{ and }\;\sum_{j=0}^{k-1}\|v^b_{j+1}-v^b_{j}\|\le M.
\end{equation}
Indeed, \eqref{bvderivative} is exactly the condition employed above to justify nontriviality \eqref{nontriv2}.

\begin{Remark}{\bf (optimality conditions for problem $(P)$).}\label{necP} {\rm It is not hard to observe while following the limiting procedures developed in Theorems~\ref{conver} and \ref{necopt} that the passage to the limit as $\t\dn 0$ in the optimality conditions obtained for $(P^\t)$
in Theorem~\ref{necopt} leads us to necessary optimality conditions for intermediate local minimizers in problem $(P)$ with the validity of {\em all} the relationships \eqref{co}--\eqref{nontriv2} of this theorem {\em but} the second measure nonatomicity condition (b). However, the optimality conditions for $(P)$ derived in this way may {\em degenerate} in the sense that for any given {\em feasible} solution to $(P)$ we can find {\em some}
collection of dual elements satisfying the nontriviality condition \eqref{nontriv2} such that all the conditions \eqref{co}--\eqref{nontriv2} hold for them.
Indeed, this happens when
\begin{equation*}
\lm=0,\;p(\cdot)=0,\;\gg(\cdot)=0,\;\xi=\delta_{\{0\}}\;\mbox{ (Dirac measure at 0)}
\end{equation*}
and the adjoint arc $q(\cdot)$ of bounded variation on $[0,T]$ is constructed as follows:
\begin{equation*}
q(t):=\left\{\begin{array}{ll}
\big(0,-2\ou(0),0)&\mbox{ for }\;t=0,\\
0&\mbox{ for }\;t\in(0,T],
\end{array}\right.
\end{equation*}
where $\ou(\cdot)$ is the $u$-part of the given feasible solution $\oz(\cdot)$ to $(P)$.
Nevertheless, it is important to emphasize as illustrated by the examples in Section~\ref{sec-exa} that,
even in the degenerate case, the aforementioned necessary optimality conditions allow us to eliminate nonoptimal solutions and find optimal ones.}
\end{Remark}

Finally in this section, we consider yet another sweeping optimal control problem much related to $(P)$,
where the control actions $\ou_i(\cdot)$ in normal directions are {\em fixed} and the optimization is
provided by $b$-controls changing the position of the moving polyhedron. This problem can be modeled in the following form $(\Tilde P)$:
\[
\mbox{minimize }\;\Tilde J[x,b]:=\varphi\big(x(T)\big)+\int_0^T\Big(\ell_1\big(t,x(t),b(t),\dot{x}(t)\big)+\ell_2\big(\dot{b}(t)\big)\Big)dt
\]
subject to the constraints in \eqref{sw-con1} and \eqref{sw-con2}, where $u_i(\cdot)=\ou_i(\cdot)$, $i=1,\ldots,m$,
are fixed absolutely continuous functions on $[0,T]$. Since the equality constraints \eqref{sw-con3} or \eqref{delta}
are not imposed, there is no difference between problem $(\Tilde P)$ and its $\tau$-perturbations as before.

We have the following necessary optimality conditions for the new problem under consideration.

\begin{Theorem}{\bf (necessary conditions for problem with fixed normal directions).}\label{necoptufix}
Let $\oz(\cdot)=(\ox(\cdot),\ob(\cdot))$ be a given r.i.l.m.\ for problem $(\Tilde P)$, and let the LICQ
condition hold at $\oz(\cdot)$. Suppose that the assumptions of Theorem~{\rm\ref{necopt}} hold whenever appropriate.
Then there exist $\lm\ge 0$, an adjoint arc $p(\cdot)=(p^x,p^b)\colon[0,T]\to\R^{n+m}$
absolutely continuous on $[0,T]$, $L^\infty$-functions $(w(\cdot),v(\cdot))$ satisfying $(w(t),v(t))\in{\rm co}\,\partial\ell(t,\oz(t),\dot\oz(t))$
for a.e.\ $t\in[0,T]$ with $\ell=\ell_1+\ell_2$, and a measure $\gamma\in C^*([0,T];\R^m)$ such that for all $i=1,\ldots$
we have the optimality relationships \eqref{etajk-defl}, \eqref{orthogonall}, \eqref{hamilton2ub}, \eqref{etajkl}, \eqref{pxkkc}, \eqref{alphak1c}, and \eqref{nontriv2} holding together with the first measure nonatomicity condition $(a)$. Moreover, \eqref{hamiltonx} and \eqref{q} read as
\begin{equation*}
\dot p(t)=\lm w(t)\;\mbox{ for a.e. }\;t\in[0,T],
\end{equation*}
\begin{equation*}
q(t)=p(t)-\Big({\int_{[t,T]}}\sum\limits_{i=1}^m\ou_i(s)d\gamma_i(s),-\int_{[t,T]}d\gamma(s)\Big)\;\mbox{ for a.e. }\;t\in[0,T].
\end{equation*}
If finally $\langle\ou_{i}(0),\ox(0)\rangle<\ob_{i}(0)$ as $i=1,\ldots,m$, we have the enhanced nontriviality $(\lambda,p(T))\ne 0$.
\end{Theorem}
{\bf Proof.} Following the proof of Theorem~\ref{discrete} shows that \eqref{th}--\eqref{psiu} hold with
$\theta_j^{uk}=0$, $w_j^{uk}=0$, $v_j^{uk}=0$ for all $j=0,\ldots,k-1$, $\xi^k=0$, and $p^{uk}=0$. Inclusion \eqref{euler1} reads now as
\begin{equation}\label{euler2}
\begin{split}
&\left(\frac{p_{j+1}^{xk}-p_j^{xk}}{h_k}-\lambda^k w_j^{xk},\frac{p_{j+1}^{bk}-p_j^{bk}}{h_k}-\lambda^k w_j^{bk},p_{j+1}^{xk}-\lambda^k\Big(v_j^{xk}+\frac{1}{h_k}\theta_j^{xk}\Big)\right)\\&\in N\Big(\Big(\ox_j^k,\ou_j^k,\ob_j^k,\frac{\ox_{j+1}^k-\ox_j^k}{h_k}\Big);\gph F(\cdot,\ou(t_j^k),\cdot)\Big),\quad j=0,\ldots,k-1.
\end{split}
\end{equation}
Note that $F$ satisfies the qualification condition in \cite[Corollary~3.17]{m-book1}, which allows us to deduced that
\begin{equation*}
D^*F_u(\ox_j^k,\ou(t_j^k),\ob_j^k)(v)\subset\text{proj}_{\R^n\times\R^m}D^*F(\ox_j^k,\ou(t_j^k),\ob_j^k)(v),
\end{equation*}
where $v:=-p_{j+1}^{xk}+\lambda^k\Big(v_j^{xk}+\frac{1}{h_k}\theta_j^{xk}\Big)$. Employing the coderivative estimate \eqref{cod_inclusion} gives us \eqref{nontriv1}--\eqref{p_iteration}, \eqref{vstar}--\eqref{pxkk}, and \eqref{alphak1}. Then the proof is completed by using the same arguments as in Theorem~\ref{necopt}. $\h$

\section{Examples and Applications}\label{sec-exa}
\setcounter{equation}{0}

We split this section into six examples, which are of a different scale. The first one describes an application of
the obtained results to a class of elastoplasticity problems, which can be modeled via the sweeping process over
controlled polyhedral moving sets. The second example addresses a particular sweeping process known as
the play-and-stop operator, which has various applications to practical models in physics, mechanics, engineering, etc.
The other examples illustrate special features of the established necessary conditions in determining optimal solutions
to the controlled sweeping process in one- or two-dimensional settings.

\begin{Example}{\bf (quasistatic elastoplasticity with hardening).}\label{elastoplasticity}
{\rm We refer the reader to the book \cite[Chapters~2--4]{HR} for models of this type (with no control) and mechanical
processes they describe with the notation therein; see also some related models in \cite{kre}.
This example is particularly inspired by models in {\em quasistatic} small-strain elastoplasticity with hardening.
Note that an optimization problem of {\em static} plasticity with linear kinematic hardening was studied in \cite{hermey},
where the external forces are taken as static controls. Here we adopt an essentially different {\em dynamic} approach,
which seems to be more realistic from the viewpoint of mechanical applications. Besides allowing the natural
{\em time evolution}, we also treat the underlying {\em yield criterion} as a control action.
This leads us to the following optimal control model of {\em dynamic elastoplasticity optimization}:
to design an elastoplastic material by (dynamically) adjusting its yield criterion in order to minimize an appropriate cost.

To proceed in more detail, consider a body in $\mathbb{R}^3$ whose displacement from the initial position is
$\mathbf{u}(t,x)$. The {\em strain} $\boldsymbol{\ve}$ is the symmetric part of the gradient of $u$, i.e.,
\[
\boldsymbol{\ve}=\frac12\Big(\nabla\mathbf{u}+\nabla^*\mathbf{u}\Big).
\]
It can be decomposed into the sum of the plastic strain $\mathbf{p}$ and of the elastic strain $\mathbf{e}$ by
$\mathbf{e}=\boldsymbol{\ve}-\mathbf{p}$. The {\em stress} $\boldsymbol{\sigma}$ depends on the elastic strain as
$\boldsymbol{\sigma}=\mathbf{C}\mathbf{e}$, where $\mathbf{C}$ is the {\em elasticity tensor} and $\boldsymbol{\sigma}$
satisfies the equilibrium equation $\mathrm{div}\,\boldsymbol{\sigma}+\mathbf{d}=0$ on an open set $\Omega$ with
smooth boundary that contains all the possible positions of the body together with the boundary condition
$\boldsymbol{\sigma}\cdot\mathbf{n}=\mathbf{c}$ on $\partial\Omega$. Here $\mathbf{n}$ denotes the external normal to
$\Omega$ while $\mathbf{d}$ and $\mathbf{c}$ represent the external forces that are taken as control actions.
We assume that $\boldsymbol{\sigma}=\boldsymbol{\sigma}(t)$ and $\mathbf{p}$ are {\em independent of $x$}.
This corresponds to the so-called {\em pseudo-rigid body}; see, e.g., \cite{SD}.
Of course, a more realistic model requires dependence on $x$, but this would lead us to
considering the sweeping process with an infinite-dimensional state space, which is beyond the scope of this paper.

If the material undergoes a linear kinematic hardening, then the ``plastic flow law'' is given by
\begin{equation}\label{swepelasto}
\dot{\mathbf{p}}\in N(\boldsymbol{\sigma}-k\mathbf{p};K)
\end{equation}
(see, e.g., \cite[pp.\ 89--90]{HR}), where $k$ is a positive constant, and where $K$ is a compact convex
subset of $\R^{3}$ called the ``region of admissible stresses." There is a number of interesting practical
models of this type with a polyhedral region of admissible stresses; e.g., it is a {\em hexagon} in the
model with the Tresca yield criterion described in \cite[p.\ 63]{HR}). Denoting $q:=-k\mathbf{p}$, we have
\[
K:=\big\{q\in\R^{3}\big|\;\langle\sigma+q,u_i\rangle\le b_i,\;i=1,\ldots,m\big\},
\]
which induces the polyhedral moving set
\[
C(t):=\big\{q\in\R^{3}\big|\;\langle q,u_i(t)\rangle\le b_i(t)-\langle\sigma(t),u_i(t)\rangle,\;i=1,\ldots ,m\big\}.
\]
This allows us to reformulate model \eqref{swepelasto} as the {\em controlled sweeping process}
\[
-\dot{q}(t)\in N\big(q(t);C(t)\big)
\]
over the moving polyhedron $C(t)$ with the control functions $u_i(t)$, $b_i(t)$, and $\sigma(t)$.
The theory developed in this paper can be readily applied to optimize the class of models under
consideration with respect to general cost functions depending on the state and control variables as well as their velocities.
Observe that our necessary conditions do not involve $\mathbf{d}$ and $\mathbf{c}$ directly, but only $\boldsymbol{\sigma}$.}
\end{Example}

The next example concerns a particular model, which appears in the literature is several contexts.

\begin{Example} {\bf (play-and-stop operator).}\label{play} {\rm This name is associated with the sweeping process given by
\begin{equation}\label{plstop}
-\dot{x}(t)\in N\big(x(t);b(t)-Z\big),\;x(0)\in b(0)-Z,
\end{equation}
where $x\in\R^n$ (in general $x$ belongs to a Hilbert space), where $Z$ is a closed and convex set
(polyhedron in our case), and where $b\colon[0,T]\to\R^n$ is absolutely continuous.
We refer the reader to, e.g., \cite[Section~7]{smb} and the bibliographies therein for more details on such operators and their applications.

To describe the possibility of applying our results, consider for simplicity the case when $Z$ is the symmetric rectangle centered at the origin
\[
Z:=\big\{(x_1,x_2)\in\R^2\big|\;|x_1|\le\bb_1,\;|x_2|\le\bb_2\}\;\mbox{ with }\;\bb_1,\bb_2>0
\]
and the control is provided by $b(t)=(b_1(t),b_2(t))$ for $t\in[0,T]$ under the fixed constant $u$-components
$u_1:=(1,0)$, $u_2:=(0,1)$, $u_3:=-u_1$, $u_4:=-u_2$. Then we have
\begin{eqnarray*}
C(t)= \big\{x\in\R^2\big|\;\langle x,u_1\rangle\le\bb_1+b_1(t),\langle x,u_2\rangle\le\bb_2+b_2(t),\langle x,u_3\rangle\le\bb_1+b_1(t),
\langle x,u_4\rangle\le\bb_2+b_2(t)\big\}
\end{eqnarray*}
and are in a position to apply the necessary optimality conditions of Theorem~\ref{necoptufix} to the optimal
control problem described by \eqref{eq:MP} and \eqref{plstop} with the moving set $C(t)$ given above.}
\end{Example}

Now we present several examples illustrating some characteristic features of the necessary optimality conditions derived in
Theorem~\ref{necopt} and also showing how to use these conditions to determine intermediate local minimizes.
Note that in the examples below the running cost is {\em convex} in velocity variables, and so there are no
difference between intermediate and relaxed intermediate local minimizers.\vspace*{0.03in}

The following simple one-dimensional example (as well as the more involved subsequent ones) illustrates the
procedure of solving problems $(P^\t)$ whenever $0\le\t\le T$ by using Theorem~\ref{necopt}, even in the
case of the possible degeneracy for $(P)$ as discussed in Remark~\ref{necP}.

\begin{Example}{\bf (calculating optimal controls in one-dimensional problems).}\label{Ex2} {\rm Let $(P)$ be given by
\begin{eqnarray}\label{ex2a}
n=m=T=1,\;x_0=0,\;\varphi(x):=\frac{(x-1)^2}{2},\;\mbox{ and }\;\ell(t,x,u,b,\dot x,\dot u,\dot b):=\frac{1}{2}\dot b^2.
\end{eqnarray}
It follows from the structure of $(P)$ that we can put $\ou(t)=-1$ on $[0,1]$ and thus consider the minimization of the cost functional \eqref{eq:MP} with
data \eqref{ex2a} subject to the dynamic constraint
\begin{eqnarray}\label{ex2b}
-\dot x(t)\in N\big(x(t);C(t)\big),\;\mbox{ where }\;C(t):=\big\{x\in\R\big|\;-x(t)\le b(t)\big\}\;\mbox{ for a.e. }\;t\in[0,1].
\end{eqnarray}
It is easy to see that the variational problem in \eqref{ex2a} and \eqref{ex2b} admits an optimal solution; it also follows from the general theory due to the convexity and coercivity of the the integrand $\ell$ in \eqref{ex2a} with respect to velocity. Thus we can apply the necessary optimality conditions of Theorem~\ref{necopt} to the problem in \eqref{ex2a}, \eqref{ex2b} and determine in this way its local solution. Employing \eqref{co}--\eqref{nontriv2} with taking into account that \eqref{pk-1c} carries no information in this case give us the following relationships valid for a.e.\ $t\in[0,1]$:
\begin{eqnarray*}
\begin{array}{ll}
(1)\;\;w=0,\;v=(0,0,\dot{\ob});\quad (2)\;\;-\ox(t)<\ob(t)\Longrightarrow\eta(t)=0;\\
(3)\;\;\eta(t)>0\Longrightarrow q^x(t)=0;\quad (4)\;\;\dot{\ox}(t)=\eta(t);\\
(5)\;(\dot{p}^x,\dot{p}^u,\dot{p}^b)(t)=\big(0,\eta(t)q^x(t),0\big);\quad (6)\;\;q^u(t)=0,\;\;q^b(t)=\lambda\dot{\ob}(t);\\
(7)\;(q^x,q^u,q^b)(t)=(p^x,p^u,p^b)(t)-\disp\Big(-\int_t^1d\gamma,\int_t^1\ox(s)d\gamma-2\int_t^1d\xi,-\int_t^1d\gamma\Big);\\
(8)\;p^b(1)\ge 0,\;-\ox(1)<b(1)\Longrightarrow p^b(1)=0;\quad (9)-p^x(1)=\lambda\big(\ox(1)-1\big)-p^b(1);\\
(10)\;\lambda+\|q(0)\|+\|p(1)\|\ne 0.
\end{array}
\end{eqnarray*}
We consider first that case where $\dot{\ox}(t)\ne 0$ for a.e.\ $t\in[0,1]$. Then it is evident that $\dot{\ox}(t)=-\dot{\ob}(t)>0$ 
and hence it follows from
(4) that $q^x(t)=0$ for a.e.\ $t$. Observe furthermore that (5) implies that $p(\cdot)$ is constant, which ensures by (7) that 
$\int_{[t,1]}d\gamma$ is
constant on $(0,T]$ as well. The latter means that either $\gamma$ is zero or it is a Dirac measure concentrated at $t=1$. 
In both cases we have that $q^b(t)$
is constant by (5) and (7), and so $\dot{\ob}(\cdot)$ is constant by (6) provided that $\lm\ne 0$; otherwise we do not have 
enough information to proceed. Assuming $\lambda\ne 0$ yields in this case that there is only one feasible trajectory 
satisfying the necessary optimality conditions; namely $(\ox(t),\ou(t),\ob(t))=(t/2,-1,-t/2)$ with the cost value of $1/4$. 
The case of $\lambda=0$, where no information can be deduced on $\ox$, cannot be ruled out. 
To examine finally the opposite case of $\ox(\cdot)=\ob(\cdot)=0$, we see by the same arguments as above that 
$q(\cdot)$ is constant on $(0,T]$ with $q^x(t)=\lambda-p^b(1)$. This choice satisfies necessary optimality conditions with 
the cost value of $1/2>1/4$, and thus we found a reasonable candidate
to be an optimal solution to this problem.}
\end{Example}

Observe that the problem in Example~\ref{Ex2} can be also treated by necessary optimality conditions of the Pontryagin Maximum 
Principle from conventional control theory by taking into account that the state constraint therein is active, i.e., $x(t)=-b(t)$, 
and so we can consider 
$\dot b(t)$ as the new control. In fact, PMP allows us to show in this setting that $(\ox(t),\ou(t),\ob(t))=(t/2,-1,-t/2)$ is a global minimizer.
However, it is not the case in the following modification of the previous example, where the moving constraint is not active, 
and we cannot reduce the sweeping process to a conventional control system.

\begin{Example} {\bf (necessary conditions for the controlled sweeping process versus PMP).}\label{Ex4}
{\rm The only difference of this example in comparison with Example~\ref{Ex2} is that the running cost is given now by
\begin{eqnarray*}
\ell(t,x,u,b,\dot x,\dot u,\dot b):=\frac12\big((b-1)^2+|\dot b|^2\big).
\end{eqnarray*}
The trivial choice of $b(t)=1$ for all $t\in[0,T]$, gives us the value $1/2$ of the cost function \eqref{eq:MP}. 
If instead the moving constraint is active, then obviously $u(t)\equiv -1$ and $x(t)=-b(t)\ge 0$ on $[0,T]$, which shows that
\[
\frac{(x(1)-1)^2}{2}+\frac12\int_0^1\big((-x(t)-1)^2+|\dot{b}(t)|^2\big)dt>\frac{(x(1)-1)^2}{2}+ \frac12\ge\frac12.
\]
It is easy to see that the trivial solution above satisfies the necessary conditions in 
Theorem~\ref{necopt} (take $q^x(\cdot)=p^x(\cdot)=\lambda=1$
and let the other dual elements vanish). Since in this case the sweeping state constraint is not active, we cannot employ PMP 
as in Example~\ref{Ex2}.}
\end{Example}

The next example shows how to exclude nonoptimal solutions by using necessary optimality conditions from Theorem~\ref{necoptufix}. 
Observe that in this example the measure $\gg$ has an atom and the corresponding adjoint arc $q(\cdot)$ is discontinuous 
inside the time interval.

\begin{Example} {\bf (excluding nonoptimal solutions).}\label{ex3} 
{\rm Consider problem $(\tilde{P})$ with $\bar{u}\equiv -1$ on $[0,1]$,
$n,m,T$ as in \eqref{ex2a}, $x_0=1/5$, $\varphi(x):=(x-1)^2$, and the running cost given by
\begin{equation*}
\quad\disp\ell(t,x,b,\dot x,\dot b):=\big(b+t-s_0(t)\big)^2+\alpha|\dot b+4t-2|\;\mbox{ with }\;\alpha\ge 0,\;
s_0(t):=\begin{cases}
\Big(t-\frac{1}{5}\Big)^2&\quad{\rm if}\;t<\frac{1}{5},\\
0&\quad{\rm if}\;\frac{1}{5}\le t\le\frac{4}{5},\\
\Big(t-\frac{4}{5}\Big)\Big(t+\frac15\Big)&\quad{\rm otherwise.}
\end{cases}
\end{equation*}
It is easy to see that the couple $(\ox(t),\ob(t)\big)=\big(v_0(t),-t+s_0(t))$ on $[0,T]$ is optimal for this problem if $\al=0$, where
\[
v_0(t):=
\begin{cases}
\frac15 & \text{ if } t<\frac15,\\
t  & \text{ if } \frac15\le t\le \frac45,\\
\frac45 &\text{ otherwise}.
\end{cases}
\]
Let us show it is not (at least in the intermediate local sense) if $\al>0$. It can done it by applying the necessary optimality 
conditions of Theorem~\ref{necopt} while observing that all its assumptions are satisfied with the validity of \eqref{bvderivative} 
for $\ell_2$ nondifferentiable in $\dot b$. Furthermore, due to $\ox(0)\ou=-1/5<1/25=\ob(0)$ we can employ the enhanced 
nontriviality condition $(\lm,p(1))\ne 0$.

To proceed, deduce from \eqref{etajkl} and \eqref{orthogonall} that $q^x(t)=0$ for a.e.\ $t\in[1/5,4/5]$. Moreover, since $\eta\equiv 0$ both in
$[0,1/5]$ and in $[4/5,1]$ thanks to \eqref{etajk-defl}, we deduce from \eqref{hamiltonx} that $\dot p(\cdot)=0$, and so 
$p(t)=p(1)$ for all $t\in[0,1]$. Assuming by contradiction that $\lm \ne 0$, we get from \eqref{hamilton2ub} 
that $q^b(t)=\alpha\lambda$ for $t>3/4$ and $q^b(t)=-\alpha\lambda$ for $t<3/4$. Then \eqref{q} tells us that
\begin{equation*}
q^x(t)=\int_{[t,1]}d\gamma(s)\;\mbox{ and }\;q^b(t)=-\int_{[t,1]}d\gamma(s)\;\mbox{ for a.e. }\;t\in[0,1].
\end{equation*}
Therefore, on one hand the measure $\gamma$ is zero on $(1/5,4/5)$, while on the other it must have a nonzero mass at $t=3/4$. 
This is a contradiction, which shows that $\lm=0$.  at the same time we get $p(1)=0$ by the transversality conditions 
\eqref{pxkkc}--\eqref{alphak1c} due to
$\ox(1)\ou=-4/5<-1+6/25=\ob(1)$. This contradicts the enhanced nontriviality and thus verifies that 
$(\ox(\cdot),\ob(\cdot))$ is not optimal for $(\tilde{P})$ if $\alpha>0$.}
\end{Example}

Finally, we present a two-dimensional example that can be analyzed on the basis of Theorem~\ref{necoptufix}.

\begin{Example} {\bf (controlled sweeping process in two dimensions).}\label{Ex5} {\rm Let the data of $(P)$ be:
\begin{eqnarray}\label{ex5a}
n=m=2,\;x_0=(1,1),\;T=1,\;\varphi(x)=\frac{\|x\|^2}{2},\;\mbox{ and }\;\ell(t,x,u,b,\dot x,\dot u,\dot b):=\frac{1}{2}\big(\dot{b}_1^2+\dot{b}_2^2\big).
\end{eqnarray}
Consider the version $(\Tilde P)$ of this problem with the fixed normal vectors $u_1\equiv(1,0)$ and
$u_2\equiv(0,1)$ and apply Theorem~\ref{necoptufix} to determine optimal solutions $\ob(t)=(\ob_1(t),\ob_2(t))$
and  $\ox(t)=(\ox_1(t),\ox_2(t))$ on $[0,1]$. The necessary optimality conditions of Theorem~\ref{necoptufix} give us the relationships on $[0,1]$:

\begin{itemize}
\item[(1)] $w(\cdot)=0$, $v^x(\cdot)=0$, $v^b(\cdot)=\big(\dot{b}_1(\cdot),\dot{b}_2(\cdot)\big)$;\quad(2)$\;\dot{\ox}_i(t) \neq 0\Longrightarrow
q^x_i(t)=0$,\;$i=1,2$.
\item[(3)] $p^b(\cdot)$ is constant with nonnegative components, and $-p_i^x(\cdot)=\lambda\ox(1)+p_i^b(\cdot)u_i$, $i=1,2$ is also constant.
\item[(4)] $q^x(t)=p^x-\gamma([t,1])$,\quad $q^b(t)=\lambda\dot{\ob}(t)=p^b+\gamma([t,1])$ for a.e.\ $t\in[0,T]$.
\item[(5)] $\lambda+\|q(0)\|+\|p(1)\|\ne 0$ with $\lm\ge 0$.
\end{itemize}
Observe first that the trivial solution with $\ox(t)=(1,1)$ and $\dot{\ob}(t)=0$ on $[0,1]$ satisfies necessary conditions (take $p_1^x=p_2^x=-1$, $p_1^b=p_2^b=\gamma_1=\gamma_2=0$, and $\lambda=1$). In this case the cost value is $1$. If the $i$-th constraint is {\em pushing} (i.e., $\dot{\ox}_i(t)<0$ on a set of positive measure), it follows from (4) that $\gamma_i([t,1])$ is constant on that set and also $\dot{\ob}_i$ is constant for $i=1,2$ provided that $\lm\ne 0$, which is supposed to hold. We consider only the simplified case where pushing occurs on at most one interval. There are the following three possibilities in this case:

(a) Both constraints are pushing with constant speed at the same time in the interval $[0,\vt]$, where $0<\vt\le 1$ is to be determined.

(b) The constraints are pushing alternatively (with constant speed); by symmetry we may assume that they push for the same time,
say first $u_1$ in the interval $[0,\vt]$ and then $u_2$ in the interval $[\vt,2\vt]$, where $0<\vt\le 1/2$ is to be determined.

(c) Only one constraint is pushing in the interval $[0,\vt]$, again with constant speed; by symmetry we may assume that
the first one is active.

To proceed further, denote the constant speed of the $i$-th moving constraint by $\beta_i<0$. Then in case (a) the cost value is calculated by
\[
\frac12(\vt^2+\vt)(\beta_1^2+\beta_2^2)+\vt(\beta_1+\beta_2)+1,
\]
and it is subject to minimization over $\beta_1,\beta_2<0$ with $0<\vt\le 1$. Straightforward calculations show that
\[
\beta_1=\beta_2=\frac{-1}{1+\vt},
\]
and in this case the $\vt$-component of the gradient of the cost function is negative, and so $\vt=1$.
Thus $\beta_1=\beta_2=-\frac12$, which gives us the cost value $\frac12$.

In case (b) we have the same cost value as in (a) while $\vt\in(0,1/2]$. The same calculations tell us that the optimal cost with this strategy is obtained for $\vt=\frac12$ and $\beta_1=\beta_2=-\frac12$ and its value is $\frac{11}{16}$. In case (c) the cost value is $\frac34$ with the choice of $\vt=1$,
$\beta_1=-\frac12$, and $\beta_2=0$. Combining all the above allows us to conclude that the strategy in case (a) is the most appropriate when $\lm\ne 0$. Finally, observe that if $\lambda=0$, we do not have enough information to proceed.}
\end{Example}

\noindent {\bf Acknowledgements.} The authors are gratefully indebted to Franco Cardin and Toma\'s Roubi\v{c}ek for several bibliographical suggestions and commentaries on elastoplasticity models. We also thank Luk\'a\v s Adam and Tan Cao for their useful remarks and discussions on the original presentation.

\end{document}